%% file: parabolic-generic-transversality.tex
\documentclass[a4paper,12pt]{article}
\usepackage{amssymb,amsfonts,latexsym,amsmath,amsthm,amscd,graphpap}
\usepackage{enumerate,graphicx,xcolor,mathabx}
\usepackage[top=3cm,bottom=3cm,left=3cm,right=3cm]{geometry}

\newtheorem{theorem}{Theorem}[section]
\newtheorem{proposition}[theorem]{Proposition} 
\newtheorem{defi}[theorem]{Definition}

\newtheorem{coro}[theorem]{Corollary}
\newtheorem{conj}[theorem]{Conjecture}
\newenvironment{demo}{ \noindent \textbf{Proof:}}{\hfill$\square$\\
\vspace{0.4cm}}

\newcounter{compte}
\newenvironment{enum-i}{\begin{list}{\roman{compte})} {\usecounter{compte}
\topsep=1mm \itemsep=0.2mm \leftmargin=5mm  }}{\end{list}}
\newenvironment{enum-1}{\begin{list}{\arabic{compte})} {\usecounter{compte}
\topsep=1mm \itemsep=0.2mm \leftmargin=5mm  }}{\end{list}}
\newenvironment{enum-a}{\begin{list}{\alph{compte})} {\usecounter{compte}
\topsep=1mm \itemsep=0.2mm \leftmargin=5mm  }}{\end{list}}
\newenvironment{remarks}{ \noindent \emph{\textbf{Remarks:}} \begin{enum-1}}
{\end{enum-1} \bigskip }

\newcommand{\RR}{\mathbb{R}}
\newcommand{\NN}{\mathbb{N}}
\newcommand{\CC}{\mathbb{C}}
\newcommand{\ZZ}{\mathbb{Z}}

\newcommand{\Cc}{\mathcal{C}}
\newcommand{\Gc}{\mathcal{G}}

\newcommand{\Mc}{\mathcal{M}}
\newcommand{\Nc}{\mathcal{N}}
\newcommand{\Oc}{\mathcal{O}}
\newcommand{\Uc}{\mathcal{U}}
\newcommand{\Vc}{\mathcal{V}}
\newcommand{\Wc}{\mathcal{W}}

\newcommand{\Cg}{\mathfrak{C}}
\newcommand{\Gg}{\mathfrak{G}}

\newcommand{\grad}{{\nabla}}

\renewcommand{\Re}{\operatorname{Re}}

\newcommand{\no}{n$^{\text{o}}$}
\newcommand{\pc}{ \usefont{T1}{cmtl}{m}{n} \selectfont}

\numberwithin{equation}{section}

\begin{document}

\title{\bf Generic transversality of heteroclinic and homoclinic orbits for 
scalar parabolic equations}

\author{Pavel {\sc Brunovsk\'y}\footnote{Department of Applied Mathematics and Statistics, Comenius University Bratislava, Bratislava 84248, Slovakia.}, Romain {\sc Joly}\footnote{Universit\'e Grenoble Alpes, CNRS, 
Institut Fourier, F-38000 Grenoble, France, email: {\pc 
romain.joly@univ-grenoble-alpes.fr}} and Genevi\`eve 
{\sc Raugel}\footnote{Universit\'e Paris-Sud \& CNRS,
Laboratoire de Math\'ematiques d'Orsay, 91405 Orsay cedex, France.}}

\date{June 2019}

\maketitle
\vspace{1cm}

\begin{abstract}
In this paper, we consider the scalar reaction-diffusion
equations 
$$\partial_t u= \Delta u+f(x,u, \nabla u)$$ 
on a bounded domain $\Omega\subset\RR^d$ of class $\Cc^{2,\gamma}$.
We show that the heteroclinic and homoclinic orbits connecting hyperbolic 
equilibria and hyperbolic periodic orbits are transverse, generically with 
respect to $f$. 
One of the main ingredients of the proof is an accurate study of the singular 
nodal set of solutions of linear parabolic equations.
Our main result is a first step for proving the genericity of Kupka-Smale 
property, the generic hyperbolicity of periodic orbits remaining unproved. 
\\[2mm]
\vspace{1cm}
\end{abstract}

\section{Introduction}

Let $d\geq 2$ and let $\Omega\subset \RR^d$ be a bounded domain 
of class $\Cc^{2, \gamma}$, where $0< \gamma \leq 1$. Let $p > d$ be fixed, let 
$X=L^p(\Omega)$ and let 
$$\Delta_D~:~D(-\Delta_D)= W^{1,p}_0(\Omega) \cap 
W^{2,p}(\Omega)~\longrightarrow~X=L^p(\Omega)$$ 
be the Laplacian operator with homogeneous 
Dirichlet boundary conditions. Let $\alpha \in (1/2+d/2p,1)$, so that 
$X^\alpha=D((-\Delta_D)^\alpha)\hookrightarrow
W^{2\alpha,p}(\Omega)$ is compactly embedded in $\Cc^1(\overline{\Omega})$. 

We consider the scalar parabolic equation
\begin{equation}\label{eq}
\left\{\begin{array}{ll}\partial_t u(x,t)=\Delta_D u(x,t)+f(x,u(x,t),\nabla
u(x,t)),&\quad (x,t)\in \Omega\times (0,+\infty)\\ 
u(x,t)=0,& \quad (x,t)\in\partial\Omega\times (0,+\infty)\\ 
u(x,0)=u_0(x) \in X^\alpha, 
\end{array}\right.
\end{equation}
where $f\in\Cc^2(\overline\Omega\times\RR\times\RR^d,\RR)$ and $u(x,t)\in\RR$.

The local existence and uniqueness of classical solutions $u(t) \in
\Cc^0([0,\tau), X^{\alpha})$ of Equation \eqref{eq}, as well as the continuous 
dependence of the solutions with respect to
 the initial data $u_0$ in $X^{\alpha}$, are well known (see \cite{Henry} for
example and  Section \ref{section-basics} for more details). Thus, Eq. 
\eqref{eq} generates a local dynamical system $S(t) \equiv S_f(t)$ on 
$X^\alpha$. This dynamical system contains all the features of a classical 
finite-dimensional system: equilibrium points and periodic orbits, stable and 
unstable manifolds\ldots{} We recall the 
definition of these objects, the definition of hyperbolicity and of 
transversality in Section \ref{section-dyn}. There, we also present their 
construction in our framework. Notice that 
the realizations results of \cite{Dancer-Polacik} and \cite{Polacik5} show the
possible existence of very complicated dynamics for \eqref{eq}, such as chaotic 
dynamics, as soon as $d\geq 2$.

In what follows, for any $r \geq 2$, we denote by $\Cg^r$ the
space $\Cc^r(\overline\Omega\times\RR\times\RR^d,\RR)$ endowed with the Whitney
topology, which is a Baire space (see Appendix \ref{section-whitney} for 
definitions, including the one of generic subset). 
In fact, our result still holds if we embed $\Cc^r$ with another reasonable topology, 
but the Whitney one is the most classical. 
See \cite{GG} and Appendix \ref{section-whitney} below for more details. 

\vspace{3mm}

Our main result is as follows.
\begin{theorem} \label{th-main} {\bf Generic transversality of connecting 
orbits}\\
Let $r\geq 2$ and let $f_0\in \Cg^r$. Let $\Cc^-_0$ and $\Cc^+_0$ be two 
critical elements of the flow of \eqref{eq}, i.e. $\Cc^\pm_0$ are equilibrium 
points or periodic orbits, $\Cc^-_0=\Cc^+_0$ being possible. 

Assume that both $\Cc^-_0$ and $\Cc^+_0$ are hyperbolic. Then, there exists a 
neighborhood $\Oc$ of $f_0$ in $\Cg^r$ and a generic set $\Gg\subset\Oc$ such 
that:
\begin{enum-i}
\item there exist two families $\Cc^-(f)$ and $\Cc^+(f)$ of critical 
elements (either equilibrium points or periodic orbits) of the flow of 
\eqref{eq}, depending smoothly of $f\in\Oc$, such that 
$\Cc^\pm(f_0)=\Cc^\pm_0$ and $\Cc^\pm(f)$ is hyperbolic for any $f\in\Oc$.
\item for any $f$ in the generic set $\Gg\subset\Oc$, the unstable manifold 
    $W^u(\Cc^-(f))$ and the stable manifold $W^s(\Cc^+(f))$ intersect 
transversally, i.e. $W^u(\Cc^-(f))\pitchfork W^s(\Cc^+(f))$.
\end{enum-i}
\end{theorem}

\begin{figure}[htp]
\begin{center}
\resizebox{13cm}{!}{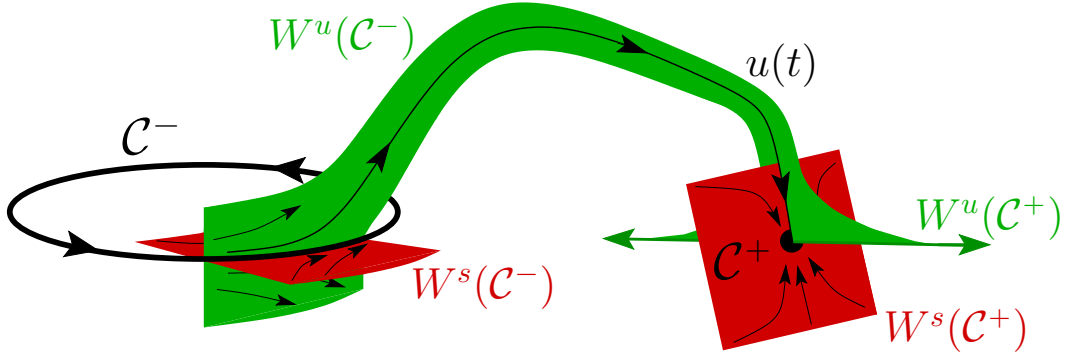}
\end{center}
\caption{\it A typical transversal heteroclinic orbit connecting a periodic orbit $\Cc^-$ and an equilibrium point $\Cc^+$. If $\Cc^\pm$ are hyperbolic, they admit stable and unstable manifolds. Theorem \ref{th-main} states that, the transversality of $u(t)$ in this picture is a generic situation in the parabolic equation \eqref{eq}. 
Here $\Cc^-$ is a periodic orbit and $\Cc^+$ is an equilibrium point. This situation is robust to perturbation
and yields several important qualitative properties of the dynamics. See the third part 
of this introduction for the historical background and Section \ref{section-dyn} for precise
definitions. \label{fig-intro}}
\end{figure}

Theorem \ref{th-main} states the generic transversality of connecting orbits, 
i.e. heteroclinic and homoclinic orbits, between hpyerbolic critical elements (either 
equilibrium points or periodic orbits). See Figure \ref{fig-intro} for an 
illustration of a typical transversal connecting orbit.
This is a first step to obtain the 
genericity of Kupka-Smale property. Below in this introduction, we recall the 
historical background and previous results. We discuss about the missing 
ingredients to obtain the genericity of the whole Kupka-Smale property in 
Appendix \ref{section-hyperbolicity}.

Notice that we do not need to assume global existence of solutions in Theorem 
\ref{th-main}. Indeed, we consider closed and connecting orbits, which are by 
definition solutions $u(t) \in X^{\alpha}$ of \eqref{eq}, which are defined for 
any time $t \in \RR$ and are
also uniformly bounded for $t \in \RR$. So, we do not really care about
solutions of Eq. \eqref{eq}, which do not exist globally. If one wants that all
solutions of \eqref{eq} exist for $0 \leq  t \leq \infty$,  one has to 
introduce additional hypotheses on $f$ (see \cite{Pol} for instance).

We also enhance that our result may apply to settings different from 
\eqref{eq}. Typically, we can choose different boundary conditions or consider 
systems of parabolic equations. We discuss this kind of straightforward 
generalizations in Section \ref{section-beyond}.

\vspace{5mm}

{\noindent {\textbf{Observability of trajectories, unique 
continuation and singular nodal sets.}}}\\[1mm]
As in the classical case of generic transversality in ODEs, the proof of 
Theorem \ref{th-main} consists in finding suitable perturbation of the 
non-linearity $f$ for breaking the non-transversal orbits.  
Of course, even if the general patterns and the spirit of the proofs stay the 
same, working with PDE's instead of ODE's gives rise to several more or less 
delicate technical problems. For example, for proving generic properties, 
instead of using Thom's transversality theorem (as in \cite{Peix66}), we will 
apply a Sard-Smale theorem stated in Appendix \ref{section-Sard-Smale}. Here, we 
want to emphasize that, in the case of PDE's, the main new difficulty
arises in the construction of appropriate perturbations. When one wants to prove
that a property is dense in the set of ODE's of the form $\dot y(t)=g(y(t))$, 
for each $g$, one has to construct a particular perturbation $\varepsilon h$ 
with small $\varepsilon$ such that the flow of $\dot y(t)=(g+\varepsilon 
h)(y(t))$ satisfies the desired property. The 
vector field $h$ of the perturbation can be
chosen freely and localized, so that his support intersects the 
trajectory of $y(t)$ only in the neighborhood of $y(t_0)$. In the case of 
PDE's, we have to construct a perturbation $h$ of the non-linearity such that the 
flow of $\partial_t u(x,t)=\Delta u(x,t) + (f+\varepsilon h)(x,u(x,t),\grad 
u(x,t))$ 
satisfies the desired property. Therefore, the perturbation $h$ of the PDE's is
of the form
\begin{equation}\label{eq-pertu}
u(\cdot)\in X^\alpha ~\longmapsto~ h(\cdot,u(\cdot),\grad u(\cdot))
\end{equation}
Since two distinct functions $u_1$ and $u_2$ can take the same value 
 $(u_1(x_0), \nabla u_1(x_0)) = (u_2(x_0), \nabla u_2(x_0))$  at a given 
$x_0\in\Omega$, the
perturbations of the form \eqref{eq-pertu} are in general ``non local" in
$X^\alpha$. Given
a particular trajectory $u(t)$ and a time $t_0$, our strategy consists in
constructing a perturbation \eqref{eq-pertu}, whose support, even if it is
large, intersects $u(x,t)$ only around $(x_0, t_0)$, which allows to consider
\eqref{eq-pertu} as a local perturbation. However, this construction is not
straightforward  and requires deep properties of the PDE. This problem is 
close to observability questions: how much information on a solution 
$u(t)$ can we get from the observation at one point $x_0$ of $u(x_0,t)$ and 
$\grad u(x_0,t)$? 

To be able to prove Theorem \ref{th-main}, we will prove in Section 
\ref{section-1to1} results of the following type.
\begin{theorem}\label{th-1to1-intro1}{\bf Injectivity properties of connecting 
orbits}\\
Let $f\in\Cc^\infty(\Omega\times\RR\times\RR^d,\RR)$. Let $u(t)$ be a 
heteroclinic or homoclinic orbit connecting two critical elements. Then there 
exists a dense open set of points $(x_0,t_0)\in\Omega\times\RR$ such that the 
curve $t\mapsto (u(x_0,t),\grad u(x_0,t))$ is one to one at $t_0$ in the 
sense that:
\begin{enum-i}
\item $(\partial_t u(x_0,t_0),\grad \partial_t u(x_0,t_0))\neq 0$,
\item for all $t\in \RR$, $(u(x_0,t),\grad 
(x_0,t))=(u(x_0,t_0),\grad (x_0,t_0))$ $\Longrightarrow$ $t=t_0$.
\end{enum-i}
\end{theorem}

The above result is a key property to be able to construct a suitable 
perturbation of the non-linearity $f$ in the proof of Theorem \ref{th-main}.
The following result is similar: it shows that 
the period of a periodic orbit of the parabolic equation may be observed 
very locally. This result is not required in the proof of our main theorem, but it may 
be interesting by itself and could be a key step to prove the generic 
hyperbolicity of periodic orbits (see the discussion of Appendix 
\ref{section-hyperbolicity}).
\begin{theorem}\label{th-1to1-intro2}{\bf Pointwise observability of the 
period of periodic orbits}\\
Let $f\in\Cc^\infty(\overline\Omega\times \RR\times\RR^d,\RR)$. Let $p(t)$ be a
periodic solution of \eqref{eq} with minimal period $\omega>0$. Then there
exists a dense open set of points $(x_0,t_0)\in \Omega\times\RR$ such that 
$$(p(x_0,t),\nabla p(x_0,t)) = (p(x_0,t_0),\nabla 
p(x_0,t_0))~~~\Longrightarrow~~~ t \in t_0+\ZZ \omega~.$$
\end{theorem}

Notice that in dimension $d=1$, the above results are true for all $(x_0,t_0)$ 
and not only for a dense subset (see \cite{JR}).

To obtain these injectivity properties of $(x,t)\longmapsto(x,u(x,t), \nabla 
u(x,t))$, where $u(t)=S_f(t)u_0$ is a bounded complete trajectory of \eqref{eq}, 
we set 
$$v(x,t,\tau)=u(x,t) - u(x,t+\tau)~,$$
and remark by using the equation \eqref{eq} that $v(x,t)$ is the solution of a 
linear parabolic equation with parameter of the form
\begin{equation}\label{eq-intro-nodal}
\partial_t v(x,t,\tau)=\Delta v(x,t,\tau) + a(x,t,\tau)v(x,t,\tau)+b(x,t,\tau).
\nabla_x v(x,t,\tau)~,
\end{equation}
in the domain $\Omega$ of $\RR^d$.  The non-injectivity points of the image 
of $(x,u(x,t), \nabla u(x,t))$, $(x,t)\in \Omega \times \RR$, 
are described by the nodal singular set of \eqref{eq-intro-nodal}, that is, 
the set of points $(x,t,\tau)$ where 
$v(x,t,\tau)$ and $\nabla_x v (x,t,\tau)$ both vanish. The singular nodal set of
solutions of the parabolic equations, with coefficients independent of the
parameter $\tau$, have already been studied in \cite{HanLin} and in \cite{Chen}
for example. Here, generalizing an argument of \cite{HardtSimon} and applying
unique continuations results (recalled in Section \ref{section-basics}), we 
prove the following theorem, see Section \ref{section-nodalset}. 
\begin{theorem}\label{th-Intro-nodalset} {\bf Singular nodal sets for 
parabolic PDEs with parameter}\\
Let $I$ and $J$ be open intervals of $\RR$. Let
$a\in\Cc^\infty(\Omega\times
I\times J,\RR)$ and $b\in\Cc^\infty(\Omega\times
I\times J,\RR^d)$ be bounded coefficients. Let $v$ be a strong solution
of \eqref{eq-intro-nodal} with Dirichlet boundary conditions. Let $r\geq 1$
and assume that $v$ is of class $\Cc^r$ with respect to $\tau$ and of class
$\Cc^\infty$ with respect to $x$ and $t$. Assume moreover that the 
null solution is not part of the family, that is that, there are no time
$t\in I$ and parameter $\tau\in J$ such that $v(.,t,\tau)\equiv 0$.

Then, the set
\begin{equation*}
 \{ (x_0,t_0) \in \Omega \times I ~ | ~ \nexists\, \tau  \in J 
\text{ such that } (v(x_0,t_0,\tau),\nabla v(x_0,t_0,\tau) )=(0,0)\}
\end{equation*}
is generic in $\Omega\times I$. In other words, the projection of all the 
singular nodal sets of the family of solutions $v(\cdot,\cdot,\tau)$ is 
negligible in $\Omega\times I$. 
\end{theorem}

\vspace{5mm}

{\noindent {\textbf{Historical background: the Morse-Smale and 
Kupka-Smale properties.}}}\\[1mm]
The transversality of unstable and stable manifolds stated in Theorem 
\ref{th-main} is related to the local stability of the qualitative dynamics. 
In the modeling of phenomena in physics or biology, we often work on 
approximate systems: some phenomena are neglected, only
approximate values of the parameters are known, or we work with a 
discretized version of the system for simulation by computer\ldots{} 
Therefore, it 
is important to know if such small approximations may qualitatively change the 
dynamics or not. Unfortunately, when perturbing general dynamical systems, 
drastic changes in the local or global dynamics can occur due for example to 
bifurcation phenomena. Thus, the common hope is that these bifurcations are 
rare, that is, that the systems, whose dynamics are robust under perturbations, 
are dense or generic. Here, we obtain the generic transversality of heteroclinic 
and homoclinic orbits between critical elements. 
Roughly, Theorem \ref{th-main} says that if we consider two hyperbolic 
closed orbits of the flow of the parabolic equation \eqref{eq} and if we 
observe a connecting orbit between them, then, ``almost surely'' this 
connection still remains after small perturbations of the system (numerical 
computation, changes of the parameters\ldots). 

Such stability questions have been extensively studied in the case of vector
fields or iterations of maps. In 1937, Andronov and Pontrjagin  introduced the 
fundamental notion of structurally stable vectors fields  (``syst\`{e}mes 
grossiers" or ``coarse
systems"), that is, vector fields $X_0$ which have a neighborhood 
$V_0$ in the $\Cc^1$-topology such that any vector field $X$ in $V_0$ is
topologically equivalent to $X_0$. In 1959 (\cite{Smale0}), Smale defined the
class of nowadays called Morse-Smale dynamical systems on compact 
$n-$dimensional manifolds, that is, systems for which the 
non-wandering set consists only in a finite number of hyperbolic equilibria and 
hyperbolic periodic orbits and for which the  intersections of the stable and 
unstable manifolds of equilibria and periodic orbits are all transversal.
Peixoto (\cite{Peix62}) proved that Morse-Smale vector fields are dense 
and have structurally stable qualitative dynamics in compact orientable
two-dimensional manifolds. In 1968, Palis and Smale (\cite{Palis}, 
\cite{Palis-Smale}) proved the structural stability of the Morse-Smale dynamical 
systems in any dimension. However, the density of Morse-Smale systems fails in 
dimension higher than two, due to ``Smale horseshoe''. In 1963, Smale 
(\cite{Smale63}) and also Kupka (\cite{Kupka}) introduced the Kupka-Smale vector 
fields, that is, the vector fields for which all the equilibria and periodic 
orbits are 
hyperbolic and the intersections of the stable and unstable manifolds of 
equilibria and periodic orbits are all transversal. They both show the density 
of such systems in any dimension (see also \cite{Peix66}).
The qualitative dynamics of Kupka-Smale systems are locally stable: periodic 
orbits, the local dynamics around them and their connections move smoothly when 
a parameter of the equation is changing. 

For the partial differential equations (PDE's in short), the history of
structural stability and of local stability is more recent. 
Notice that a trajectory of the dynamical system $S(t)$ generated by such a PDE 
is of the form $t\mapsto S(t)u_0=u(\cdot,t)$, where $u(x,t)$ is the solution of 
the PDE with initial data $u_0(x)$. In particular, the trajectory moves in a 
functions space (often a Sobolev space), which is infinite-dimensional. As a 
generalization of \cite{Palis} and \cite{Palis-Smale}, \cite{HMO84} and  
\cite{Oli} proved that Morse-Smale and Kupka-Smale properties are still 
meaningful in infinite-dimensional systems for 
the problem of stability of the qualitative dynamics. Therefore, there is a 
great interest in obtaining generalizations of the above mentioned 
finite-dimensional generic results. Notice that, if we want to 
get a meaningful genericity result, we have to allow perturbations only in the 
same class of PDE's. Typically, the parameter with respect to which the 
genericity is obtained is the non-linearity $f$.

The first example of transversality of unstable and stable manifolds for PDE's 
is due to Henry (\cite{Henry85}) in 1985 for the reaction-diffusion equation in 
the segment
\begin{equation}\label{eq-dim1}
\partial_t u=u_{xx} +f(x,u,u_x), ~~~(x,t)\in
(0,1) \times (0,+\infty)
\end{equation}
with Dirichlet, Neumann or Robin boundary conditions. More strikingly, he 
obtained the noteworthy property that the stable and unstable manifolds of two 
hyperbolic equilibria of \eqref{eq-dim1} always intersect transversally. A key 
ingredient for proving this automatic transversality is the use of 
the non-increase of the ``Sturm number" or ``zero number"
\cite{Sturm} of the solutions of the corresponding linearized parabolic 
equations. In addition to this automatic transversality, the gradient structure 
proved in \cite{Zelen} shows the genericity of Morse-Smale property for the 
flow of \eqref{eq-dim1} with separated boundary conditions. 

If we consider \eqref{eq-dim1} with periodic boundary conditions, that is the 
parabolic equation on the circle $S^1$
\begin{equation}\label{eq-S1}
\partial_t u=u_{xx} +f(x,u,u_x), ~~~(x,t)\in
S^1 \times (0,+\infty)
\end{equation}
then the gradient structure fails but 
the flow of \eqref{eq-S1} still has particular properties equivalent to the 
ones of two-dimensional ODEs, such as the Poincar\'{e}-Bendixson
property proved in \cite{FMP} (the reader interested in the correspondence 
between the dynamics of \eqref{eq-dim1} and the ones of low-dimensional ODEs 
may consider the review paper \cite{JR-vf}). In 2008, still using the powerful 
tool of the ``zero number", Czaja and Rocha (\cite{Czaja-Rocha}) proved that, 
for the parabolic equations on the circle \eqref{eq-S1}, the stable and unstable 
manifolds of hyperbolic periodic orbits always intersect transversally.  In 
2010, the second and third authors completed the results of Czaja and Rocha. 
More precisely, they proved in \cite{JR} that the equilibria and periodic orbits are 
hyperbolic, generically with respect to the 
nonlinearity $f$. They also proved that the stable and unstable
manifolds of hyperbolic critical elements ${\mathcal C}^-$ and 
${\mathcal C}^+$  intersect transversally, unless both critical elements
${\mathcal C}^-$ and ${\mathcal C}^+$ are equilibria of same Morse index and
moreover that, generically with respect to $f$, such connecting orbits between
equilibria with the same Morse index (\cite{JR2}) do not exist. Finally, the 
Poincar\'{e}-Bendixson theorem of \cite{FMP} yields that, generically with 
respect to $f$, the equation \eqref{eq-S1} is Morse-Smale
(see \cite{JR2}).

Concerning spatial dimension higher than $d=1$, the generic transversality of 
stable and unstable manifolds has been shown in 1997 by the first author and P. 
Pol\'a\v{c}ik (\cite{Bruno-Pola}) in the case $f\equiv f(x,u)$, that is, for the 
equation
\begin{equation}\label{eq-grad}
\partial_t u=\Delta u +f(x,u), ~~~~ (x,t)\in
\Omega \times (0,+\infty)
\end{equation}
with $\Omega\subset \RR^d$, $d\geq 2$. As a consequence, since \eqref{eq-grad} 
is a gradient system, they deduce that, under additional dissipative 
conditions on the non-linearity, the Morse-Smale property holds for 
the flow \eqref{eq-grad} generically with respect to $f\in \Cc^2$.
It is noteworthy, as shown by Pol\'a\v{c}ik
(\cite{Polacik99}), that this generic transversality property is not true  if 
one considers homogeneous functions 
$f(x,u) \equiv f(u)$ only. 

We also mention that generic transversality properties have been shown by the 
authors for various gradient damped wave equations, see \cite{Bruno-Raugel} 
and \cite{Joly}.

Due to the realization results of Dancer and 
Pol\'a\v{c}ik, \cite{Dancer-Polacik} and \cite{Polacik5}, we know that the 
dynamics of the flow of the general parabolic equation \eqref{eq} in dimension 
$d\geq 2$ may be as complicated as chaotic flows. We may only hope to prove 
the genericity of the Kupka-Smale property and not of the Morse-Smale one. Notice that the 
flow of \eqref{eq} is not gradient (periodic orbits may exist) and the very 
particular and helpful ``zero number property" of spatial dimension $d=1$ fails.
In the present paper, we prove the generic transversality property. The generic 
hyperbolicity of equilibrium points is already proved in \cite{JR} in any 
space dimension. Thus, the generic hyperbolicity of periodic 
orbits is the only remaining step to obtain the genericity of the Kupka-Smale 
property. 

Some years ago, in a preliminary draft of this paper, we were convinced to have 
proved the genericity of the Kupka-Smale property. However, 
Maxime Percy du Sert pointed to us a gap in the proof of generic hyperbolicity 
of periodic orbits. We did not manage to fill it. Recently, two of the three authors 
passed away and we decided to publish the results as obtained together. In particular, 
we prove the generic transversality only (unlike claimed in \cite{JR-vf}). In Appendix \ref{section-hyperbolicity}, we quickly discuss our ideas to obtain the generic hyperbolicity of periodic orbits and 
indicate where the gap remains.

\vspace{5mm}

{\noindent {\textbf{Plan of the article.}}}\\[1mm]
In Section \ref{section-basics}, we recall the classical existence and 
uniqueness properties of the solutions of the scalar parabolic equation and the 
corresponding linear and linear adjoint equations. We also review unique 
continuation properties, which are fundamental in this paper. 
In Section \ref{section-dyn}, we remind some basic definitions
such as hyperbolicity of critical elements and we state the main properties of 
the dynamical system $S_f(t)$, namely the existence of $\Cc^1$ immersed 
finite-codimensional (resp. finite-dimensional) stable (resp. unstable) 
manifolds of hyperbolic critical elements. 
Section \ref{section-nodalset} is devoted to the study of 
the singular nodal sets and to the proof of Theorem 
\ref{th-Intro-nodalset}. 
In Section \ref{section-1to1}, we show that Theorem 
\ref{th-Intro-nodalset} leads to one-to-one properties such as Theorems 
\ref{th-1to1-intro1} and \ref{th-1to1-intro2}.
Using these tools, in Section \ref{section-transverse}, we prove Theorem 
\ref{th-main}, i.e. we show the generic transversality of heteroclinic and 
homoclinic orbits of the parabolic equation \eqref{eq}. 
Section \ref{section-beyond} contains discussions about some generalizations of 
Theorem \ref{th-main}. 
We conclude by two appendices recalling the basic facts about the Whitney 
topology 
and Sard-Smale theorems, which will be used in this paper, and one appendix 
discussing the still open problem of generic hyperbolicity of periodic orbits 
of \eqref{eq}.

\vspace{5mm}

{\noindent {\textbf{Dedication:}}} Very sadly, both Pavol Brunovsk\'y and Genevi\`eve Raugel passed away before the publication of this article, respectively in december 2018 and in may 2019. They were still working actively on the manuscript and the present version is exactly the one which have been completed by them. This article is dedicated to their memories.

\vspace{3mm}

{\noindent {\textbf{Acknowledgement:}}} The last two authors have been funded by the research project {\it ISDEEC} ANR-16-CE40-0013. 


\section{Some basic results on parabolic PDEs}\label{section-basics}

\subsection{Local existence and regularity results of the parabolic equation  
\eqref{eq}}

The solutions of the scalar parabolic equation \eqref{eq} exist  locally and are
unique, see for example \cite{Pazy} or \cite{Henry}. In the whole paper, 
$\alpha$ belongs to the open interval $(\frac{1}{2} + \frac{d}{2p},1)$.
We recall that we use the notation $f\in\Cc^r(E,\RR)$ to indicate the 
regularity of $f$, i.e. to say that the function $f:E\rightarrow \RR$ is of 
class $\Cc^r$. Where a topology is required (smooth dependences on $f$ etc.), 
the notation $\Cg^r(E,\RR)$ refers to the space $\Cc^r(E,\RR)$ endowed with the 
Whitney topology (see Appendix \ref{section-whitney}).

\begin{proposition}\label{Pr-exist}
Let $r\geq 1$ and $f\in \Cc^r(\overline\Omega\times\RR\times\RR^d,\RR)$.
\begin{enum-i}
\item For any $u_0 \in X^{\alpha}$, there exists a maximal time $T (u_0)>0$
such that  \eqref{eq} has a unique classical solution
$S_f(t)u_0=u(t)\in 
\Cc^0([0,T],X^\alpha) \cap \Cc^1((0,T],X^{\beta}) \cap \Cc^0((0,T],
D(-\Delta_D))$, for any $0 \leq \beta <1$ and for any $0 <T \leq T(u_0)$.  If  
$T(u_0)$ is finite, then $\|u(t)\|_{X^\alpha}$ goes to $+\infty$ 
when $t < T(u_0)$ tends to $T(u_0)$. 

\noindent Moreover, $t\mapsto \partial_t u(t)$ is locally H\"older 
continuous from $(0,T]$ into $X^\beta$, for $0 \leq \beta <1$. In particular, 
$u(\cdot) \equiv S_f(\cdot) u_0$ belongs to the space $\Cc^0((0,T], 
W^{3,p}(\Omega)) \cap 
\Cc^1((0,T], W^{s,p}(\Omega))$, for any $s<2$, and thus 
belongs to the spaces $\Cc^0((0,T], \Cc^2(\overline\Omega)) \cap 
\Cc^1((0,T], \Cc^1(\overline\Omega))$ and $\Cc^1(\overline\Omega \times 
[\tau,T],\RR)$, 
for any $0 <\tau <T$. If, in addition,  the first derivatives $D_u f(x, \cdot, 
\cdot)$ and 
$D_{\nabla u}f(x, \cdot, \cdot)$ are Lipschitz-continuous on the bounded sets of
$\overline\Omega \times \mathbb{R} \times \mathbb{R}^d$, then  
$u(\cdot)$ belongs to $\Cc^1((0,T], W^{2,p}(\Omega))
\cap \Cc^2((0,T], W^{s,p}(\Omega))$, for any $s<2$ and hence
$u(\cdot)$ also belongs to $\Cc^2(\overline\Omega \times [\tau,T],\RR)$,
for any $0 <\tau <T$.

\item  For any $u_0 \in X^{\alpha}$, for any $T<T(u_0)$, there exist a
neighborhood 
$\Uc \equiv \Uc(T)$ of $u_0$ in $X^\alpha$ and a neighborhood $\Vc \equiv
\Vc(T)$ of $f$
 in $\Cg^1$ such that, for any $v_0\in \Uc$ and any $g \in \Vc$, 
 $v(t) \equiv S_g(t)v_0$ is well defined on $[0,T]$, depends
continuously on $v_0 \in X^{\alpha}$ and $g \in \Cg^1$, and there exists a 
positive number 
$R \equiv R(T, \Uc ,\Vc)$ such that  $((S_g(t)v_0)(x), (\nabla S_g(t)v_0)(x))$ 
belongs to the ball 
$B_{\RR^{d+1}}(0, R)$, for all $(t,v_0,g,x) \in [0,T] \times \Uc \times \Vc 
\times \overline\Omega$.

\item Moreover, for any $u_0 \in X^{\alpha}$, for any $T<T(u_0)$, the map 
$(t,u_0) \in (0,T] \times \Uc \mapsto S_{f}(t) u_0 \in X^{\alpha}$ is of class 
$\Cc^r$ and, in  particular, 
$S_f(t)$ is a local semigroup of class $\Cc^r$. In addition, there exists a
neighborhood $\mathcal{W}$ of 
$f$ in the space $\Cg^r(\overline\Omega\times [-2R,2R] \times [-2R,2R]^d,\RR)$  
such that the map
$(t,u_0,g) \in (0,T] \times \Uc \times \mathcal{W} \mapsto S_{g}(t) u_0 \in 
X^{\alpha}$ is of 
class $\Cc^r$.
\end{enum-i}
\end{proposition}

\begin{remarks}
\item The statement (i) is a direct consequence of the existence and
regularity results 
given in \cite[Chapter 3]{Henry} and of elliptic regularity properties. We only
want to emphasize that, since the solution $u(\cdot) \equiv S_f(\cdot) u_0$
belongs to $\Cc^0([0,T], X^\alpha)$ and that 
$X^{\alpha}$ is continuously embedded in $\Cc^1(\overline\Omega)$, $u(\cdot)$  
automatically  belongs to the space $\Cc^0([0,T), \Cc^1(\overline\Omega))$. 
Since $u(\cdot)$ is a classical solution and belongs to $\Cc^0((0,T], 
W^{2,p}(\Omega)) \cap  \Cc^1((0,T], W^{1,p}(\Omega))$, 
$f(x,u, \nabla u) -\partial_t u$ is in the space $\Cc^0((0,T],  
W^{1,p}(\Omega))$ and 
the regularity 
properties of the elliptic equation
$$
\Delta_D u = \partial_t u - f(x,u, \nabla u)~,
$$
imply that $u(\cdot)$ belongs to the space
 $\Cc^0((0,T], W^{3,p}(\Omega)) \subset  \Cc^0((0,T],\Cc^2(\overline\Omega))$.
\item Statements (ii) and (iii) are also easy consequences of 
\cite[Theorem 3.4.4 and Corollary 3.4.5]{Henry}. We want to point out that, for 
any $u_0 \in X^\alpha$ and any $0 <T < T(u_0)$, there exists $R_0 >0$ such that 
$(u(x, t), \nabla u(x,t))$, for all
$(x,t )\in \overline \Omega \times [0,T]$ is bounded in $\RR^{d+1}$ by a 
positive number $R_0 \equiv R_0(u_0, T)$. Since $g(x,u(x,t), \nabla u(x,t))$ 
depends only on the values of $x$, $u(x,t)$ and 
$\nabla u(x,t)$, we can show, by applying the continuity results of 
\cite[Section 3.4]{Henry}, that,
for any $R >R_0$, for any $0 < \varepsilon < (R-R_0)/2$, there exists a 
positive number $\eta$ 
such that, for any $g(\cdot, \cdot,\cdot) \in \Cc^r(\overline\Omega\times 
[-R,R] \times [-R,R]^d,\RR)$, $\eta$-close to $f$ in the classical norm of 
$\Cg^r(\overline\Omega\times [-R,R] \times [-R,R]^d,\RR)$,
$((S_g(t)u_0)(x), (\nabla S_g(t)u_0)(x))$ belongs to the ball 
$B_{\RR^{d+1}}(0, R_0 +\varepsilon)$, for all $(x,t) \in  \overline\Omega 
\times [0,T]$.
\item Notice that the statement (ii) of Proposition \ref{Pr-exist} implies that
the maximal time
$T(u_0)$ is a lower-semi-continuous function of the initial data $u_0$
\end{remarks}

As we have already seen, the parabolic equation has a smoothing effect at any 
finite positive time.
 If the boundary of the domain $\Omega$ was of class $\Cc^{\infty}$ and 
$f$ belonged to $\Cc^\infty(\overline\Omega\times\RR\times\RR^d,\RR)$, the 
solutions of 
Eq.\eqref{eq} would be in $\Cc^{\infty}(\overline\Omega \times [\tau,T],\RR)$ 
for any $0 < \tau < T < T(u_0)$. However, if $f \in 
\Cc^\infty(\overline\Omega\times \RR \times\RR^d,\RR)$, we can still show that 
the solutions are regular in the interior of $\Omega$, even if $\Omega$ is of 
class $\Cc^{2,\alpha}$ only. 

In the whole paper, we say  that $u(t): t \in \RR \mapsto u(t)$ is a {\it 
bounded complete solution (or trajectory)} of \eqref{eq} if it is a solution of 
\eqref{eq}, defined for any $t \in \RR$ and
 bounded in $X^{\alpha}$, uniformly with respect to $t \in \RR$. 
 
Since we are only interested in the regularity of the bounded complete 
solutions of \eqref{eq}, we will  state a $\Cc^{\infty}$-regularity result for 
such solutions.
\begin{proposition} \label{Pr-Cinfini} Assume that $f$ belongs to 
$\Cc^\infty(\overline\Omega\times \RR \times \RR^d,\RR)$. Then, any bounded 
complete
solution $u(t)$ of \eqref{eq} belongs to $\Cc^\infty(\Omega\times\RR,\RR)$.
More precisely, for any open set $\mathcal{O}$, such that $\overline 
{\mathcal{O}} \subset \Omega$, for any $R >0$, any $m \in \mathbb{N}$, any $k 
\in \mathbb{N}$, and any $q \in [1,\infty]$, there exists a positive constant 
$K(\mathcal{O},R,m,k,q)$, such that any bounded complete solution $u(t)$, with 
$\sup_{t \in\RR}~\|u(t)\|_{X^\alpha}\leq R$, satisfies 
\begin{equation}
\label{estimksq}
\sup_{t \in \mathbb{R}} \left\|\frac{d^k 
u}{dt^k}(t)\right\|_{W^{m,q}(\mathcal{O})} \leq 
K(\mathcal{O},R,m,k,q)~.
\end{equation}
\end{proposition}

\begin{demo} We will not give all the details of the proof, but will indicate 
only the main arguments. The proof consists in a recursion argument with 
respect to $k$ and $m$.
Let $u(t)$ be a bounded complete solution of \eqref{eq} satisfying 
$\sup_{t \in \mathbb{R}} \| u(t)\|_{X^{\alpha}} \leq R$.

\vspace{2mm}

\noindent {\underline{First ste}p\underline{:}} Since $f$ belongs to 
$\Cc^\infty(\overline\Omega\times \RR \times \RR^d,\RR)$, by \cite[Corollary 
3.4.6]{Henry}, 
the function $t \in \mathbb{R} \mapsto u(t) \in X^{\alpha}$ is of class 
$\Cc^k$, for 
any $k \in \mathbb{N}$ and $\frac{d^k u}{dt^k}(t) \in \Cc^0(\mathbb{R}, 
X^{\alpha} \cap W^{2,p}(\Omega)) \cap \Cc^1(\mathbb{R}, X^{\beta})$, for any 
$\beta <1$, is a classical solution of the equation
\begin{equation}
\label{eqdtk}
\frac{d}{dt}(\frac{d^{k}u}{dt^{k}}) = \Delta \frac{d^{k} u}{dt^{k}} + 
\frac{d^{k}}{dt^{k}}(f(x, u, \nabla u))~.
\end{equation}
We notice that the term $\frac{d^{k}}{dt^{k}}(f(x, u, \nabla u))$ can be 
computed by using 
the Faa Di Bruno formula \cite{FaadiBruno} and its generalization
\cite{FaadiBrunoGen} as follows.
We introduce the $(d+1)$-dimensional vector 
$w(x,t) = (u,\grad u)(x,t)$, that is $w_1=u$ and $w_{i+1}=\partial_{x_i} 
u$. Using the generalized Faa Di Bruno formula (\cite{FaadiBrunoGen}), we can 
write,
\begin{align}
\frac{d^k}{dt^k}(f(x,u(x,t), & \nabla u(x,t))) = \sum_{m_j=1, |m|=1} 
D^m_{w}f(x,w(x,t)) \frac{d^k}{dt^k} (w_j)(x,t) \nonumber\\
 &~~~ + \sum_{2 \leq |m|\leq k} D^m_{w}f(x,w(x,t)) \sum_{p(k,m)} k! 
\Pi_{j=1}^{k} 
 \frac{\bigl[\frac{d^{\ell_j}}{dt^{\ell_j}}w\bigr]^{n_j}}{(n_j ! ) [\ell_{j} 
!]^{|n_j|}} \nonumber\\
& \equiv \sum_{m_j=1, |m|=1} D^m_{w}f(x,w(x,t)) \frac{d^k}{dt^k} (w_j)(x,t) + 
g_k(x,t) \label{2FDBk}
\end{align}
where $p(k,m) = \{(n_1,\ldots,n_k; \ell_1, \ldots, \ell_k)\,|\,\exists s\in 
\ldbrack 1,k\rdbrack,~n_i = \ell_i =0 \text{ for } 1 \leq i \leq n-s \}$
and $g_k$ contains only derivatives with respect to $t$ of order less or equal 
to $k-1$.

We notice that the estimate \eqref{estimksq} for $k=0$, $m=2$ and $q=p$ is a 
direct consequence of the hypothesis and of Proposition \ref{Pr-exist}.
Using \eqref{2FDBk}, the fact that $W^{1,p}(\Omega)$ is an 
algebra and the bound 
$\sup_{t \in\RR}\|u(t)\|_{X^\alpha}\leq R$, one shows by recursion on  
$k$ that
\begin{equation}
\label{estimdtk-1}
\sup_{t \in \mathbb{R}} \|\frac{d^k u}{dt^k}(t)\|_{W^{2,p}(\Omega)} \leq 
C_2(R,k)~,
\end{equation}
where $C_2(R,k)$ is a positive constant depending only on $R$, $k$ (and of 
$f$). Like in the remarks following Proposition 
\ref{Pr-exist}, the elliptic regularity properties allow also to deduce from
Eq.\eqref{eqdtk} and from the estimate \eqref{estimdtk-1} that, 
\begin{equation}
\label{estimdtk}
\sup_{t \in \mathbb{R}} \|\frac{d^k u}{dt^k}(t)\|_{W^{3,p}(\Omega)} \leq 
C_3(R,k)~,
\end{equation}
where $C_3(R,k)$ is a positive constant depending only on $R$, $k$ (and of 
$f$). 

\vspace{2mm}

\noindent {\underline{Second ste}p\underline{:}} One easily shows, by recursion 
on $n \in \mathbb{N}$ (and also $k$) that, 
\begin{equation}
\label{estimdtkn}
\sup_{t \in \mathbb{R}} \|\frac{d^k u}{dt^k}(t)\|_{W^{3+n,p}(\mathcal{O})} \leq 
C_{3+n}(\mathcal{O}, R,k)~.
\end{equation}
Indeed, let $O_j$, $j=1,2, \ldots,n+1$, be a sequence of regular open sets such 
that 
$\overline{\mathcal{O}} \subset O_{n+1} \subset \overline{O_{n+1}} 
\subset O_{n} \subset \ldots \subset O_{j+1} \subset \overline O_{j+1} 
\subset O_{j} \ldots \subset O_1 \subset \overline{O_1} \subset \Omega$ and 
$\varphi_j$, $j=1,2, \ldots$, be  a 
corresponding sequence of regular functions
such that $\varphi_j(x) \in [0,1]$, $x \in \overline \Omega$, and $\varphi_j(x) 
\equiv 0$, for $x \in 
\overline{\Omega}  \setminus O_{j}$ and $\varphi_j(x) \equiv 1$, for $x \in  
O_{j+1}$. We recall that, by the remarks following Proposition \ref{Pr-exist}, 
one already knows that the estimates \eqref{estimdtk}
hold for any $k \in \mathbb{N}$. We remark that $\varphi_1u$ is a solution of 
the elliptic equation
\begin{equation}
\label{eqO1}
\Delta (\varphi_1 u)= \varphi_1 \frac{du}{dt} + u \Delta \varphi_1 + 2 \nabla u 
\cdot \nabla \varphi_1
- \varphi_1 f(x,u,\nabla u)~
\end{equation}
where $\varphi_1 \frac{du}{dt} + u \Delta \varphi_1 + 2\nabla u \cdot \nabla 
\varphi_1
- \varphi_1 f(x,u,\nabla u)$ belongs to $W^{3-1,p}(O_1) \cap W^{1,p}_0(O_1)$. By 
the elliptic regularity 
results, $\varphi_1u$ belongs to $W^{3+1,p}(O_1)$ and 
\begin{equation}
\label{est01O1}
\sup_{t \in \mathbb{R}} \| \varphi_1 u(t)\|_{W^{3+1,p}(O_1)} \leq  C_{3+1}( O_1, 
R,0,\varphi_1)~,
\end{equation}
where $C_{3+1}( O_1, R,0,\varphi_1)$ is a positive constant depending only on 
$O_1$, $R$, 
$\varphi_1$. Likewise, writing the elliptic equality satisfied by  $\Delta 
(\varphi_1 (\frac{d^k}{dt^k}u))$ 
and using the equalities \eqref{eqdtk} and \eqref{2FDBk}, one shows, by 
recursion on $k$, that
$\frac{d^k}{dt^k}(\varphi_1u)$ belongs to $W^{3+1,p}(O_1)$ and 
\begin{equation}
\label{estk1O1}
\sup_{t \in \mathbb{R}} \| \frac{d^k}{dt^k}(\varphi_1 u)(t)\|_{W^{3+1,p}(O_1)} 
\leq  
C_{3+1}(O_1,R,k,\varphi_1)~,
\end{equation}
where $C_{3+1}(O_1,R,k,\varphi_1)$ is a positive constant depending only on 
$O_1$, $R$, $k$ and 
$\varphi_1$. We notice that $\frac{d^k}{dt^k}(\varphi_1u)(x) =
\frac{d^k}{dt^k}u(x)$, for any $x \in O_2$. 

We next assume that $\frac{d^k}{dt^k}(\varphi_j u)$ belongs to $W^{3+j,p}(O_j)$ and 
that the estimates
\eqref{est01O1} and \eqref{estk1O1} hold with $1$ replaced by $j$. Remarking that
$\varphi_{j+1}u$ is a solution of the elliptic equation
\begin{equation}
\label{eqOj+1}
\Delta (\varphi_{j+1} u)=\varphi_{j+1} \frac{du}{dt} + u \Delta \varphi_{j+1}
 + 2\nabla u \cdot \nabla \varphi_{j+1} - \varphi_{j+1} f(x,u,\nabla u)~
\end{equation}
where $\varphi_{j+1} \frac{du}{dt} + u \Delta \varphi_{j+1} + 2\nabla u \cdot 
\nabla \varphi_{j+1}- \varphi_{j+1} f(x,u,\nabla u)$ belongs to $W^{3+j-1,p}(O_{j+1}) \cap 
W^{1,p}_0(O_{j+1})$, we at once
show that $\varphi_{j+1} u$ belongs to $W^{3+j+1,p(O_{j+1}) \cap 
W^{1,p}_0(O_{j+1})}$ and that the estimate \eqref{est01O1} holds with $1$ 
replaced by $j+1$. Likewise, one shows by recursion on $k$ that
$\frac{d^k}{dt^k}(\varphi_{j+1}u)$ belongs to $W^{3+j+1,p}(O_{j+1})$ and that 
the estimate \eqref{estk1O1} holds with $1$ replaced by $j +1$. Thus, we have 
proved by recursion on $n$ and $k$ that
$\frac{d^k}{dt^k}(u)$ belongs to $W^{3+n,p}(\mathcal{O})$  
and that the estimates \eqref{estimdtkn} are satisfied. 

The general estimate \eqref{estimksq} is a direct consequence of the estimates
\eqref{estimdtkn} and the classical Sobolev embedding theorem.
\end{demo}

\subsection{The linear and linear adjoint equations}

Let $0 \leq s <T$ and let $a(\cdot)\in\Cc^1([0,T],L^{\infty}(\Omega))$  and 
$b(\cdot)\in\Cc^1 ([0,T],W^{1,\infty}(\Omega)^d)$. We consider solutions
$v$ of the linear parabolic equation 
\begin{equation}
\label{lineariz}
\begin{split}
v_t(x,t) = &\Delta_D v(x,t) + a(x,t)v(x,t) +b(x,t).\nabla v(x,t)~, \quad t>s,  x 
\in \Omega, \cr
v(x,s) = & v_s ~.
\end{split}
\end{equation}
In what follows, we denote $A(t)$ the operator 
$$
A(t) = \Delta_D + a(x,t) . + b(x,t).\nabla~.
$$
Equation \eqref{lineariz} arises either when one linearizes the parabolic
equation \eqref{eq} along a solution
$u$, in which case we have
\begin{equation}\label{eq-ab1}
\left\{\begin{array}{l} a(x,t)=f'_u(x,u(x,t),\nabla u(x,t)) \\
b(x,t)=f'_{\nabla u}(x,u(x,t),\nabla u(x,t))\end{array}\right.
\end{equation}
or when one considers the difference $v(t)=u_2(t)-u_1(t)$ between two solutions
$u_1$ and $u_2$ of \eqref{eq}, in which case we have
\begin{equation}
\label{eq-ab2}
\left\{\begin{array}{l} a(x,t)=\int_0^1 f'_u(x,(\theta u_2
+(1-\theta) u_1)(x,t),\nabla (\theta u_2
+(1-\theta) u_1)(x,t))d\theta \\
b(x,t)=\int_0^1 f'_{\nabla u}(x,(\theta u_2
+(1-\theta) u_1)(x,t),\nabla (\theta u_2
+(1-\theta) u_1)(x,t)) d\theta   \end{array}\right.
\end{equation}
Notice that, since $f$ belongs to 
$\Cc^2(\overline\Omega\times\RR\times\RR^d,\RR)$,
 due to Proposition \ref{Pr-exist}, in both cases the
coefficients of \eqref{lineariz} belong to $\Cc^1((0,T],W^{1,\infty}(\Omega))$. 
Since in what follows, 
we are mainly applying the results of this section to  bounded complete 
trajectories, we can consider, without loss of generality, that the 
coefficients of \eqref{lineariz} belong to 
$\Cc^1([0,T],W^{1,\infty}(\Omega))$.  

\begin{proposition}\label{prop-lin}
Let $r\in [1,\infty)$  and let $v_s\in L^r(\Omega)$. Equation
\eqref{lineariz} has a unique solution $v(t) \equiv U(t,s) v_s \in
\Cc^0([s,T], L^r(\Omega))\cap \Cc^1((s,T], L ^r(\Omega))
\cap \Cc^0((s,T],W^{2,r}(\Omega)\cap W^{1,r}_0(\Omega))$ satisfying
$v(s)=v_s$. Moreover, $v: t \in (s,T] \mapsto v(t) \in X^{\alpha}$ is H\"older 
continuous 
 and belongs to $\Cc^1((s,T], L ^q(\Omega))
\cap \Cc^0((s,T],W^{2,q}(\Omega)\cap W^{1,q}_0(\Omega))$ for any $q\in
[1,+\infty]$. In particular $v\in \Cc^0((s,T],\Cc^1(\overline\Omega))$.
\end{proposition}

\begin{demo}
For the existence, uniqueness and regularity of the solution of $v(t) \equiv 
U(t,s) v_s$  in 
$\Cc^0([s,T], L^r(\Omega))\cap \Cc^1((s,T], L ^r(\Omega))
\cap \Cc^0((s,T],W^{2,r}(\Omega)\cap W^{1,r}_0(\Omega))$, we refer to 
\cite[Theorem
7.1.3]{Henry}. To prove that $v(t)$ belongs to any space $L^q(\Omega)$ (and thus
to $X^\alpha$), we will use a bootstrap argument. Assume that $v_s$ belongs to 
$L^r(\Omega)$
and set $r=r_0$. By \cite[Theorem 7.1.3]{Henry}, $v(s+\delta)\in 
W^{2,r_0}(\Omega)$ for any
$\delta>0$. If $d -2r_0 \leq 0$, then, $v(s+\delta)\in W^{2,r_0}(\Omega) 
\subset L^q(\Omega)$, 
for any positive number $q \geq 1$, by the classical Sobolev embedding. If, $d 
-2r_0 > 0$, again by the Sobolev embedding theorem, $v(s+\delta)\in 
W^{2,r_0}(\Omega) \subset L^{r_1}(\Omega)$, for 
$r_1=dr_0/(d-2r_0) = r_0 + 2r_0^2/(d -2r_0)$. We again apply \cite[Theorem
7.1.3]{Henry} to deduce that 
$v(s+2\delta)\in W^{2,r_1}(\Omega)$, for any $\delta>0$.  Again, if $d -2r_1 > 
0$, we obtain that 
$v(t+ 2\delta)\in W^{2,r_2}(\Omega) \subset L^{r_2}(\Omega)$, for $r_2
=dr_1/(d-2r_1) \geq r_1 + 2r_1^2/(d -2r_1) \geq r_0 + 2r_0^2/(d -2r_0) +
2r_1^2/(d -2r_1)$. Clearly, since the increment $r\mapsto 2r/(d-2r)$ is
increasing until $d-2r\leq 0$, after a finite number of steps, we obtain that
$v(t) \in L^q(\Omega)$.
\end{demo}

Proposition \ref{prop-lin} tells that Equation \eqref{lineariz} generates a
family of evolution operators $U(t,s)$ on $L^p(\Omega)$, which is extended
to $L^r(\Omega)$ for any $r\geq 1$.

Let now $1 <p < + \infty $, which implies that $X= L^p(\Omega)$ is reflexive. Denote by 
$p^*$ the conjugate exponent of $p$, that is, $p^*=p /(p-1)$; consider the 
adjoint space 
$X^*=( L ^p(\Omega))^* = L ^{p^*}(\Omega)$ of $X$ and the 
{\sl adjoint evolution operator} $U(t,s)^*: X^* \to X^*$. Let $T >0$; for  
$\psi_T \in L ^{p^*}(\Omega)$, we define the function $\psi :  s \in [0,T] 
\mapsto \psi (s) = U(T,s)^*\psi_T$. In general, $\psi(s)$ is only 
a {\sl weak$^*$ solution} of the equation 
\begin{equation}\label{eq-adjoint}
\partial_s\psi(x,s) = - \Delta_D \psi(x,s) - a(x,s)\psi(x,s) +
\hbox{div}(b(x,s)\psi(x,s))
\end{equation}
with $(x,s)\in\Omega\times (0,T)$ and with final data $\psi(T)=\psi_T$ in 
the weak-$\ast$ sense. More precisely, $s \in [0,T) \mapsto \psi(s) \in X^*$ is 
locally H\"older continuous, for each $\phi\in X$, $\langle \phi, \psi(s)\rangle \rightarrow \langle 
\phi, \psi_T \rangle$ when $s \rightarrow T^-$ and, for each 
$\phi \in D(A^*)$, $(\phi, \psi(s))$ is differentiable on $[0,T)$ with
$\partial_t(\phi,\psi(s)) = (A(s)\phi,\psi(s))$.

Usually, $\psi (s) = U(T,s)^*\psi_T$ is only a solution of \eqref{eq-adjoint} 
in a weak sense. But here, since 
$a(\cdot)\in\Cc^1([0,T],L^{\infty}(\Omega))$  and 
$b(\cdot)\in\Cc^1 ([0,T],W^{1,\infty}(\Omega)^d)$, $\psi(s)$ is a strong 
solution of \eqref{eq-adjoint}, as we shall see in the proposition below. 
Notice that \eqref{eq-adjoint} is a parabolic equation
solved backwards in time.

\begin{proposition}\label{prop-adjoint}
$~$
\begin{enum-1}
\item With the above notations, $\psi(s)=U(T,s)^*\psi_T$ belongs to
$\Cc^1([0,T),X^*) \cap\Cc^0([0,T),$ $W^{2,p^*}(\Omega)\cap
W^{1,p^*}_0(\Omega))$.
Moreover, it satisfies \eqref{eq-adjoint} in the strong sense and
$\psi(s)$ belongs to $\Cc^1([0,T),L^q(\Omega))\cap
\Cc^0([0,T),W^{2,q}(\Omega)\cap W^{1,q}_0(\Omega))$ for any $q\geq 1$.
\item Let $\tilde\psi_{T}\in X^*$. For any $0 < \eta < T$,
 $\tilde\psi_{T-\eta}=U(T,T-\eta)^*((-\Delta_D)^\alpha)^*\tilde\psi_{T}$ is well
defined in $X^*$. Hence, for 
 $s < T- \eta$, $\tilde \psi(s)= U(T-\eta,s)^*
\tilde\psi_{T-\eta}=U(T,s)^*((-\Delta_D)^\alpha)^*\tilde\psi_{T}$ 
belongs to  $\Cc^1([0,T-\eta),X^*)\cap \Cc^0([0,T-\eta),W^{2,p^*}(\Omega)\cap
W^{1,p^*}_0(\Omega))$ and a strong  solution of \eqref{eq-adjoint}.
\end{enum-1}
\end{proposition}
\begin{demo}
The first part of the proposition is a direct consequence of
\cite[Theorem 7.3.1]{Henry} 
on the existence and regularity of solutions for the adjoint equation and on 
the fact that the coefficients have the regularity
$a(\cdot)\in\Cc^1([0,T],L^{\infty}(\Omega))$ and 
$b(\cdot)\in\Cc^1 ([0,T],W^{1,\infty}(\Omega)^d)$. The fact that
$\psi(s)$ belongs to any $L^q(\Omega)$ is proved by recursion as in Proposition
\ref{prop-lin}. 

To show the second part of the proposition, let $\tilde\psi_{T}\in
X^*$ and let
$\varphi\in X= L ^p(\Omega)$. By Proposition \ref{prop-lin},
$U(T,T-\eta)\varphi$
belongs to $X^\alpha=D((-\Delta_D)^\alpha)$ and thus $\langle \tilde\psi_{T} |
(-\Delta_D)^\alpha U(T,T-\eta)\varphi\rangle_{ L ^{p^*}, L ^p}$ is well defined.
Therefore, $U(T,T-\eta)^*((-\Delta_D)^\alpha)^*\tilde\psi_{T}$ is well defined
and belongs to $ L ^{p^*}(\Omega)$. To finish, we apply \cite[Theorem
7.3.1]{Henry} (or the first part of the proposition) to the initial data
$\psi_{T}=U(T,T-\eta)^*((-\Delta_D)^\alpha)^*\tilde\psi_{T}$.
\end{demo}

\subsection{Unique continuation properties}

In this section, we recall some important unique continuation properties
satisfied by the linear parabolic equation \eqref{lineariz}. We enhance
that these properties will apply to solutions $v(t)\in X^\alpha$ of 
\eqref{lineariz}
with coefficients given by \eqref{eq-ab1} or \eqref{eq-ab2}. Hence, we may 
apply it to 
 the difference of two solutions of the nonlinear parabolic equation
\eqref{eq}. In particular, the
 unique continuation properties below will have fundamental consequences on the
properties of the dynamics of \eqref{eq}, such as the injectivity of the
flow. 

The following result is a direct consequence of the backward uniqueness
property stated in \cite[Theorem II.1]{Bardos-Tartar}. 
\begin{proposition}\label{prop-backward-uniqueness}
$~$
\begin{enum-1}
\item Let $T>0$. Let $a(x,t)\in L^\infty(\Omega\times (0,T))$ and let
$b(x,t)\in L^\infty(\Omega\times (0,T))^d$. Let $v(t)\in
L^2((0,T), H^1_0(\Omega))$ be a solution of the linear parabolic
equation \eqref{lineariz}.
Then, $v(T)\equiv 0$ in $\Omega$ if and only if $v$ vanishes identically in
$(0,T)\times\Omega$.

\item Likewise, assume that $a(x,t)\in L^\infty(\Omega\times (0,T))$, that
$b(x,t)\in L^\infty(\Omega\times (0,T))^d$ and that $D_{x_i} b(x,t) \in
L^\infty(\Omega\times (0,T))^d$, $0 \leq i \leq d$.
 Let $\psi(t)\in L^2((0,T), H^1_0(\Omega))$ be a solution of the adjoint linear 
equation \eqref{eq-adjoint}.
Then, $\psi(0)\equiv 0$ in $\Omega$ if and only if $\psi$ vanishes identically 
in
$(0,T)\times\Omega$. 
\end{enum-1}
\end{proposition}

Let now $u_1$ and $u_2$ be two solutions on the time interval $[0,T]$ of the
equation \eqref{eq}. 
We already remarked that $v(t) = u_2 (t) - u_1(t)$ satisfies the linear 
equation \eqref{lineariz} with the coefficients $a$ and $b$ given by 
\eqref{eq-ab2}. By Proposition \ref{Pr-exist}, the coefficients $a$, $b$ and 
the function $v(t)$ satisfy the regularity assumptions of the above proposition 
\ref{prop-backward-uniqueness}. Thus, if $u_1(T)=u_2(T)$, then $u_1\equiv u_2$ 
on $[0,T]$. This leads to state the following corollary.

\begin{coro}\label{coro-backward-unique-S}
Let $T>0$. Let $u_1(t)$ and $u_2(t)$ be two solutions on the time interval 
$[0,T]$ of the equation \eqref{eq}. If $u_1(T)=u_2(T)$, then $u_1(t) = u_2(t)$, 
for any $t \in [0,T]$. In other terms, the local dynamical system $S_f(t)$ 
generated by \eqref{eq} has the backward uniqueness property.
\end{coro}

The following result is proved in \cite{Saut-Scheurer} and shows that the
set of the zeros of the solutions of the linear parabolic equation is a closed
set with empty interior.
\begin{proposition}\label{prop-Saut-Scheurer}
Let $T>0$, $a$ and $b$ be as in Proposition \ref{prop-backward-uniqueness}. We 
assume that
$v(x,t) \in L^2((0,T),H^2(\Omega)\cap H^1_0(\Omega))$ is a solution of the 
linear parabolic
equation \eqref{lineariz}.
If $v(x,t)$ vanishes on an open non-empty subset of $\Omega\times (0,T)$, then
$v(x,t)$ identically vanishes on $\Omega\times (0,T)$.
\end{proposition}

A similar result has been obtained for the strong solutions of 
the adjoint equation in \cite[Corollary 2.12]{Fabre}.
\begin{proposition}\label{prop-Fabre}
Let $T>0$. Let $a(x,s)\in L^\infty(\Omega\times (0,T))$ and let
$b(x,s)\in L^\infty(\Omega\times (0,T))^d$. Let
$\psi(s)\in L^2((0,T),H^2(\Omega)\cap H^1_0(\Omega))$ be a solution
of the adjoint equation \eqref{eq-adjoint}. If $\psi(x,t)$ vanishes on an open
non-empty subset of $\Omega\times (0,T)$, then $\psi(x,t)$ identically vanishes
on $\Omega\times (0,T)$.
\end{proposition}

In the particular case of smooth solutions of \eqref{lineariz} (typically if
one considers global bounded solutions and a smooth non-linearity $f$), we will
need stronger properties on the zeros of the solutions in Section \ref{section-nodalset}.

We say that $v$ vanishes to infinite order in both the space
and time variables at $(x_0,t_0)$ if, for any $k \geq 1$, there is a
constant $C_k>0$, such that, for any $(x,t) \in \Omega \times
[-T,0]$,
\begin{equation}\label{ordinf}
|v(x,t)| \leq C_k( |x-x_0|^2 + |t-t_0|)^{k/2}~.
\end{equation}
We shall often apply the following unique continuation result of Escauriaza and
Fern\'{a}ndez \cite{EscFer03}.
 
\begin{proposition} \label{uniqcont} 
Assume that $v \in \Cc^0((-T,0],\Cc^2(\overline{\Omega})) \cap \Cc^1((-T,0],
\Cc^1(\overline{\Omega}))$ is a solution of \eqref{lineariz} and satisfies either homogeneous
Dirichlet or homogeneous Neumann boundary conditions. Suppose that $v$
vanishes to infinite order at $(x_0,0)$ in both the space and time
variables in the sense of \eqref{ordinf}. Assume moreover that 
there exists a positive constant $K$ such that for any  $(x,t) \in \Omega \times
(-T,0]$,
\begin{equation}\label{inequal}
|v_t(x,t) -\Delta v(x,t) | \leq K( |\nabla v(x,t) | +|v(x,t)| )~.
\end{equation}
Then, $v(x,0)$ vanishes for any $x \in \Omega$ and therefore $v(x,t)$
identically vanishes in $\Omega \times [-T,0]$.
\end{proposition}

We say that $v$ vanishes to infinite order in space at $(x_0,t_0)$ if, for any
$k \geq 1$, there is a constant $C_k>0$, such that
\begin{equation}\label{ordinf2}
|v(x,t_0)| \leq C_k |x-x_0|^k ~.
\end{equation}
{}From Proposition \ref{uniqcont} and \cite[Theorem 1]{AlesVes03}, we deduce the
following unique continuation result for solutions $v \in
\Cc^0((-T,0],\Cc^2(\overline{\Omega})) \cap 
\Cc^1((-T,0], \Cc^1(\overline{\Omega}))$ of \eqref{inequal}, which vanish to
infinite order in space. The
following result can also be deduced from Proposition \ref{uniqcont}, a simple
computation and, a recursion argument when $v(x,t)$ is a $\Cc^{\infty}$-function
in the variables $(x,t)$. Indeed,  if for example $v(x,t_0)$ vanishes 
to order $2$ (resp. $4$) in space at $(x_0,t_0)$, then, due to the equation
\eqref{lineariz}, $v_t(x,t_0)$ vanishes to order $0$ (resp. $2$) in space at
$(x_0,t_0)$. Moreover, if $v(x,t_0)$ vanishes to order $4$ in space at
$(x_0,t_0)$, deriving the equation \eqref{lineariz} with respect to $t$, one
shows that $v_{tt}(x,t)$ vanishes at order $0$ in space. Finally, continuing the
recursion argument on $k$ and on the derivatives with respect to $t$, one shows
that $v$ vanishes to infinite order at 
$(x_0,t_0)$ in both the space and time variables in the sense of \eqref{ordinf}

\begin{proposition} \label{uniqcontB} 
Assume that $v \in \Cc^0((-T,0],\Cc^2(\overline{\Omega})) \cap \Cc^1((-T,0],
\Cc^1(\overline{\Omega}))$ satisfies the inequality \eqref{inequal} and either
homogeneous Dirichlet or homogeneous Neumann
boundary conditions. Suppose also that $v$ vanishes to infinite order in space
at $(x_0,0)$, for some $x_0\in\Omega$. Then, $v(x,0)$ vanishes for any $x \in
\Omega$ and therefore $v(x,t)$ identically vanishes in $\Omega \times [-T,0]$. 
\end{proposition}


\section{The local infinite-dimensional dynamical system 
$S_f(t)$}\label{section-dyn}

In this section, we recall some basic properties of the local dynamical system 
$S_f(t)$ generated by
the parabolic  equation \eqref{eq} on $X^{\alpha}$ (if the dependence on $f$ is 
clear, we simply write 
$S(t)$). As we have seen in the introduction, the hyperbolicity of the critical 
elements  (that is, the equilibrium points and periodic orbits) and the 
transversality of the stable and unstable manifolds play a primordial role. 
Thus, we will focus on recalling the definitions and main properties of these 
objects.

\subsection{Critical elements and hyperbolicity}
Let $e\in X^\alpha$ be an equilibrium point of \eqref{eq}. The linearization 
$(D_uS(t)e)$ of
the dynamical system $S(t)$ at $e$ is given by the linear semigroup $e^{L_et}$
on $X^\alpha$, where $L_e: D(\Delta_D) \mapsto L^p(\Omega)$ is the linear 
operator defined 
by 
$$
L_e v=\Delta_D v + f'_u(x,e(x),\nabla e(x))v+f'_{\grad u}(x,e(x),\nabla
e(x)).\nabla v~.
$$
The operator $-L_e$ is a sectorial operator and a Fredholm operator with compact
resolvent. Therefore, the spectrum of $L_{e}$  consists of a sequence of 
isolated
eigenvalues of finite multiplicity, the norms of which converge to infinity.
Since the resolvent of $L_e: X \to X$ is compact, the linear 
$C_0$-semigroup $e^{L_et}$ from $X$ into $X$ is compact and its
 spectrum consists of a sequence of isolated
eigenvalues of finite multiplicity converging to $0$. By \cite[Chapter 2, 
Theorem 2.4]{Pazy},
 $\mu$ is an eigenvalue of $e^{L_e}$ if and only if $\mu =e^{\lambda}$, where 
$\lambda$ is an eigenvalue of $L_e$.
\begin{defi}
The equilibrium point $e$ is said {\bf simple} if $1$ does not belong to the 
spectrum of 
$e^{L_e}$.
The equilibrium point $e$ is {\bf hyperbolic} if $e^{L_e}$ has no
spectrum on the unit circle $S^1 \equiv \{z\in\CC~|~|z|=1\}$. 

In the case of the equation \eqref{eq}, we may equivalently say that the 
equilibrium point $e$ is simple if and only if 
$0$ is not an eigenvalue of $L_e$ and that it is hyperbolic if and only 
if $L_e$ has no eigenvalue with zero real part. 

The {\bf Morse index} $i(e)$ is the (finite) number of eigenvalues of $e^{L_e}$ 
of norm strictly larger 
than $1$ (counted with their multiplicities) or equivalently the number of 
eigenvalues of $L_e$ with positive real part.
\end{defi}

\vspace{5mm}

Let $p(t)$ be a periodic solution of the scalar parabolic equation \eqref{eq}
with period $\omega>0$. This periodic solution describes the {\bf periodic 
orbit}
$\Gamma = \{ p(t) \, | \, t \in [0,\omega) \}$.
The linearization of the dynamical system $S(t)$ along
$p(t)$ is given by the evolution operator $\Pi_{f,p}(t,s): v_s \in X^{\alpha} 
\mapsto v(t) \in X^{\alpha}$, $t \geq s$,  where $v(\tau)$ solves the 
non-autonomous equation
\begin{equation}\label{lin-periodic}
\left\{ \begin{array}{l} \partial_\tau v(x,\tau)=\Delta v(x,\tau)+
f'_u(x,p,\grad p) v(x,\tau) +  f'_{\grad u}(x,p,
\grad p) \nabla v(x,\tau)\\ 
v(x,s)=v_s(x)~.\end{array} \right.
\end{equation}
The operator $\Pi_{f,p}(\omega,0)$ is called the (corresponding) {\bf period 
map}. One remarks that
$\Pi_{f,p}(t+\omega,t) = \Pi_{f,p}(t+m\omega,t+ (m-1)\omega)$ for any $t \geq 
0$ and any $m�\in \mathbb{N}$.
Notice that $\partial_t p(t)$ is a solution of \eqref{lin-periodic} and thus 
that $1$
is an eigenvalue of $\Pi_{f,p}(\omega,0)$ with eigenvector $\partial_t p(0)$. We
emphasize that, due to the smoothing properties in finite positive time of the
parabolic equation \eqref{lin-periodic}, the operator 
 $\Pi_{f,p}(t,s): X^{\alpha} \to X^{\alpha}$, $t > s$, is compact. Therefore, 
the spectrum of 
$\Pi_{f,p}(t +\omega,t)$ consists of a sequence of isolated eigenvalues of 
finite multiplicity, converging to $0$. 
As for the linearized operator $e^{L_e}$ at  the equilibrium point $e$, $0$ is
the only point where the
spectrum of $\Pi_{f,p}(t+\omega,t)$  accumulates. Actually, by the backward 
uniqueness property,
$0$ is not an eigenvalue neither of $e^{L_e}$, nor of $\Pi_{f,p}(t+\omega,t)$. 
By \cite[Lemma 7.2.2]{Henry}, the spectrum  $\sigma(\Pi_{f,p}(t +\omega,t))$ of 
$\Pi_{f,p}(t+\omega,t)$ is independent of $t \in [0,+\infty)$. For this reason, 
the following definition makes sense.

To simplify the notation, when there is no confusion, we will simply 
write $\Pi(t,s)$ instead of $\Pi_{f,p}(t,s)$. 

\begin{defi}
A periodic solution $p(t)$ of period $\omega$ is {\bf simple} or {\bf 
non-degenerate} if the number $1$ is a simple (isolated) eigenvalue of
$\Pi_{f,p}(\omega,0)$.

The periodic solution $p(t)$ is {\bf hyperbolic} if $\Pi_{f,p}(\omega,0)$ has no
spectrum on the unit circle $S^1$ except the eigenvalue one,
which is simple and isolated. 

Since $\Pi_{f,p}(\omega,0)$ is a compact operator, the periodic solution 
$p(t)$ is  hyperbolic if and only if $1$ is a simple, isolated eigenvalue of 
$\Pi_{f,p}(\omega,0)$
and is the only eigenvalue on the unit circle. 

The {\bf Morse index} $i(p)$ of $p(\cdot)$, or the {\bf Morse index} 
$i(\Gamma)$ of  $\Gamma$, is the (finite) number of eigenvalues of 
$\Pi_{f,p}(\omega,0)$ of norm 
strictly larger than $1$ (counted with their multiplicities).
\end{defi}

In what follows, we will sometimes say that the periodic orbit
$\Gamma = \{ p(t) \, | \, t \in [0,\omega) \}$ is simple (resp. hyperbolic), 
instead of saying that $p(t)$ is simple (resp. hyperbolic).

\vspace{5mm}

A first important consequence of the simplicity property  is the persistence of 
equilibrium points and periodic orbits under perturbations.
\begin{theorem}\label{th-stab-hyperbo}
Let $r\geq 2$ be given and let $f_0\in\Cg^r$.
\begin{enum-1}
\item Let $e_0$ be a simple equilibrium point of \eqref{eq} with $f=f_0$. There
exist a neighborhood $\Nc$ of $f_0$ in $\Cg^r$ and a neighborhood $\Uc$ of
$e_0$ in $X^\alpha$ such that, for any $f\in\Nc$, there exists a unique
equilibrium point $e(f)$ in $\Uc$. This equilibrium depends 
continuously on $f \in \Cg^r$. In addition, the eigenvalues of $L_{e(f)}$ 
continuously depend 
on $f \in \Cg^r$. 

Moreover, if $e_0$ is hyperbolic, the 
neighborhoods $\Nc$ and $\Uc$ can be chosen small enough so that  
$e(f)$ is also hyperbolic and so that the Morse index $i(e)$ is equal to 
$i(e_0)$.

\item Let $p_0(t)$ be a simple periodic solution with period (resp. minimal
period)
$\omega_0$ of \eqref{eq} for $f=f_0$. There exist a neighborhood $\Nc$ of $f_0$
in $\Cg^r$, a 
positive number $\eta$ and a neighborhood $\Uc$ of $\Gamma_0= \{p_0(t)\, | \, t
\in [0, \omega_0) \}$ in $X^\alpha$ such that, for any $f \in \Nc$, there exists
a unique periodic orbit $\Gamma(f) = \{p(f)(t)\, | \, t \in [0, \omega(f)) \}$
in $\Uc$, of period (resp. minimal period) $\omega(f)$ with $|
\omega (f) - \omega_0| \leq \eta$. The  period $\omega(f)$ and the periodic
orbit $\Gamma(f)$ continuously depend on $f$. In addition, the eigenvalues of
$\Pi_{f,p(f)}(\omega(f),0)$ continuously depend on $f \in \Cg^r$. 

Moreover, if $f_0$ is hyperbolic, the neighborhoods $\Nc$ and $\Uc$ and
$\eta>0$ can be chosen small enough so that the periodic solution $p(f)(t)$ is
hyperbolic and so that the Morse index 
$i(p(f))$ is equal to the Morse-index $i(p_0)$.
\end{enum-1}
\end{theorem}

\begin{demo}
The first statement about the persistence of simple equilibria $e_0$ is very
classical. Assume that 
$\|e_0 \|_{L^{\infty}} \leq m$ and  $\|\nabla e_0 \|_{L^{\infty}} \leq m$. Then,
applying the implicit function theorem or the fixed point theorem of strict
contraction (see the proof \cite[Lemma 4.c.2]{Bruno-Pola}), one shows that there
exist a neighborhood $\mathbf{N}_0$ of $f_0$ in $\Cc^r(\overline{\Omega} 
\times [-2m,2m] \times [-2m,2m]^d)$ and a neighborhood $\mathcal{U}$ of $e_0$
in $X^\alpha$ such that for any $f\in \mathbf{N}_0$, there exists a unique
equilibrium point $e(f)$ in $\Uc$. This equilibrium depends continuously of $f
\in \mathbf{N}_0$ and, moreover, all the other properties of the first statement
hold. 
Using the restriction mapping $R$ of Section 2.1, we conclude that there exists 
a neighborhood 
$\Nc$ of $f_0$ in $\Cg^r$ such that, for any $f\in \mathcal{N}$, there exists a 
unique equilibrium point $e(f)$ in $\Uc$ and that all the other properties of 
the first statement hold.

Let $p_0(t)$ be a simple periodic solution of period $\omega_0 >0$ 
of \eqref{eq} for $f=f_0$. Assume that $\sup_{t \in [0,\omega_0)} \|p_0(t)
\|_{L^{\infty}} \leq m$ and  
$\sup_{t \in [0,\omega_0)} \| \nabla p_0(t) \|_{L^{\infty}} \leq m$. 
The statement of the persistence of a simple periodic solution $p_f(t)$ near 
$p_0(t)$ with period  
$\omega_f$ close to $\omega_0$ and also of the uniqueness (up to a time
translation) of this periodic solution, if $f$ belongs to a small enough
neighborhood of $f_0$ in 
$\Cc^r(\overline{\Omega} \times [-2m,2m] \times [-2m,2m]^d)$, is a direct
consequence of 
\cite[Theorem 8.3.2]{Henry}; it is proved by using the method of Poincar\'{e} 
sections and the implicit function theorem or the fixed point theorem of strict 
contraction (for further results in the case where the perturbations are less 
regular, see also \cite{HaRau10} and \cite{HaRau12}).  One concludes like in 
the proof of the statement 1) by using the restriction mapping  $R$ of Section 
2.1.

The continuous dependence of the eigenvalues of $L_{e(f)}$ or of
$\Pi_{f,p(f)}(\omega(f),0)$ with respect to $f \in \Cg^r$ is a consequence of
the proof of the continuity results of Kato (see \cite[Theorems IX.24, IV.31,
IV.3.18]{Kato}) and of the properties of the restriction mapping $R$.
Detailed proofs of continuity of the point spectrum can also be found in
\cite[Section 3]{HaRau92}.
\end{demo}  

Notice that a periodic solution $p(t)$ of period $\omega$ can be simple, 
whereas the same periodic solution $p(t)$, considered as periodic 
solution of period $n\omega$ can be non-simple. This is the case when the 
spectrum of $\Pi_{f,p}(\omega,0)$ contains a $n$-th root of $1$. Thus, in the 
statement 2) of Theorem \ref{th-stab-hyperbo}, when 
$p_0(t)$ is a simple periodic solution of period $\omega_0$ of \eqref{eq} for 
$f=f_0$, we do not know if $\Gamma(f) = \{p(f)(t)\, | \, t \in [0, \omega(f)) 
\}$ is the unique periodic orbit 
of \eqref{eq} in the neighborhood $\Uc$ of $\Gamma_0$ if $f$ belongs to $\Nc$.
Indeed,if the spectrum of $\Pi_{f_0,p_0}(\omega_0,0)$ contains a $n$-th root of 
unity, then it is possible that new periodic orbits of period close to 
$n\omega_0$ are created (in the case where $n=2$, it is the famous 
``period-doubling bifurcation''). 

Of course, when $p_0(t)$ is hyperbolic, no such new periodic solutions can be 
created and $\Gamma(f)$ is still isolated in the set of periodic orbits. 
Hyperbolicity is a notion independent of the chosen period.

\subsection{Stable and unstable manifolds}
We recall that a {\bf critical element} means either an equilibrium point or a 
periodic orbit of \eqref{eq}.
\begin{defi}\label{defi-Wu-Ws} Let $\mathcal{C}$ be a critical element of 
\eqref{eq}.
The {\bf global stable and unstable sets} of $\mathcal{C}$ are respectively 
defined as 
\begin{equation*}
\begin{split}
& W^s(\mathcal{C}) = \{u_0 \in X^{\alpha}\, | \, S_f(t)u_0
\xrightarrow[t \rightarrow 
  +\infty]{} \mathcal{C}\}~, \cr
& W^u(\mathcal{C}) = \{u_0 \in X^{\alpha} \, | \, \, \forall \, t\leq
0,~S_f(t)u_0\hbox{ is well defined and }
S_f(t)u_0 \xrightarrow[t \rightarrow 
  - \infty]{} \mathcal{C}\}~.
  \end{split}
\end{equation*}
Likewise, if $U_{\mathcal{C}}$ is a neighborhood of $\mathcal{C}$ in 
$X^{\alpha}$, we introduce the 
 {\bf local stable and unstable sets} of $\mathcal{C}$ defined as 
\begin{equation*}
\begin{split}
& W^s(\mathcal{C},U_{ \mathcal{C}}) \equiv  W^s_{loc}(\mathcal{C})
\equiv \{u_0 \in U_{ \mathcal{C}}\, | \, S_f(t)u_0 \in U_{ \mathcal{C}},
 t \geq 0\}~, \cr
& W^u(\mathcal{C}, U_{ \mathcal{C}}) \equiv  W^u_{loc}(\mathcal{C}) \equiv
 \{u_0 \in U_{ \mathcal{C}} \, | \, \forall \, t\leq 0,~S_f(t)u_0\text{ is
well defined and stays in }U_{ \mathcal{C}}\}.
  \end{split}
\end{equation*}
\end{defi}
If we need to specify the dependence with respect to the non-linearity $f$, we
will denote these manifolds as $W^s(\mathcal{C},U_{ \mathcal{C}},f)$ and 
$W^u(\mathcal{C}, U_{ \mathcal{C}},f)$ or 
as $W^s_{loc}(\mathcal{C},f)$ and $W^u_{loc}(\mathcal{C},f)$.

Let $e_0$ be an equilibrium point of \eqref{eq} and let $(D_{u}S(t)e_0)= 
e^{L_{e_0}t}$ be the 
corresponding linearized operator around $e_0$. We denote by $P_u$ (resp. $P_s$)
the projection
in $X^{\alpha}$ onto the space
generated by the (generalized) eigenfunctions of $e^{L_{e_0}}$
corresponding to the eigenvalues with
modulus strictly larger than $1$ (resp. with modulus strictly smaller than $1$).
Let $X_u^{\alpha} = P_u(X^{\alpha})$ and $X_s^{\alpha} = P_s(X^{\alpha})$. We 
have seen that, in the case of the parabolic equation \eqref{eq}, the Morse 
index of every hyperbolic equilibrium point is finite, which implies that 
$P_u(X) = P_u(X^{\alpha})$.

The following theorem states the existence of the local stable and unstable manifolds
near hyperbolic equilibrium points. The result is very classical. In the case of a vector field 
on a finite-dimensional compact manifold, we refer the reader to  
\cite{AbRobbin}, \cite{Palis-de-Melo}, \cite{HPS77} for example, and in the 
infinite dimensional case, we refer to \cite{Henry}, 
\cite{HMO84}, \cite{HVL}, \cite{ChHaBT}, \cite{Ruelle}. 
\begin{theorem}\label{th-eWuL-WsL}
Let $f_0$ be given in $\Cg^r$, $r \geq 2$, and let $e_0$ be a hyperbolic 
equilibrium point of 
$S_{f_0}(t)$. Then there is a neighborhood
$U_0$ of $e_0$ such that the local unstable manifold
$W^u(e_0,U_0)$ (resp. the local stable manifold $W^s(e_0,U_0)$)
is a $\Cc^r$-submanifold of dimension $i(e_0)$ (resp. codimension $i(e_0)$),
which is tangent to 
$X_u^{\alpha}$ (resp. $X_s^{\alpha}$) at $e_0$. 

More precisely, there exist a neighborhood $U_0$ of $e_0$ in $X^{\alpha}$, two 
mappings 
$h_u (f_0)\equiv h_u^0: P_u X^{\alpha}\to P_sX^{\alpha}$ and  $h_s(f_0) \equiv
h_s^0 : P_s X^\alpha  \to P_u X^{\alpha}$ of class $\Cc^r$ such that $h_u^0(0)=
0$, $Dh_u^0(0)= 0$, $h_s^0(0) =0$, $Dh_s^0(0) =0$ and 
\begin{align}
W^u_{loc}(e_0,f_0) &\equiv W^u(e_0,U_0, f_0)\nonumber \\
&= \{ v \in U_0\, | \, v = e_0 + P_u(v-e_0) + h_u^0(P_u(v-e_0)) \} \nonumber \\
W^s_{loc}(e_0,f_0) & \equiv W^s(e_0,U_0,f_0) \nonumber \\ 
&= \{ v \in U_0 \, | \, v = e_0 + P_s(v-e_0) + h_s^0(P_s(v-e_0)) 
\}~.\label{graphes}
\end{align}
Furthermore, the convergence rates to the origin are exponential.
More precisely, there are positive constants $k_1$, $k_2$ and constants  
 $0 < \gamma_2 < 1 <\gamma_1$, such that,
\begin{equation}
\label{exprate}
\begin{split}
\|S_{f_0}(t)x\|_{X} & \, \leq \, k_1 \gamma_1^t \, ,\quad \forall  \, t \leq 0 
\, , 
\quad  \forall  \, x \in W^u(e_0,U_0)~,\cr 
\|S_{f_0}(t)x\|_{X} &\, \leq \, k_2 \gamma_2^t \, , \quad \forall \, t \geq 0\, 
, 
\quad  \forall \, x \in W^s(e_0,U_0)~.\cr
\end{split}
\end{equation}
In addition, the local stable and unstable manifolds ``continuously" depend of 
the nonlinear map $f$. More precisely, there exists $\rho >0$ and, for any
$\varepsilon >0$, there is a neighborhood $\mathcal{N}$ of $f_0$ in $\Cg^r$
such that, for any $f \in \mathcal{N}$, $S_{f}(t)$ has a unique equilibrium
point $e(f)$ in the ball $B_{X^{\alpha}}(e_0,\rho)$ of center $e_0$ and radius 
$\rho$ in $X^{\alpha}$, and $\| e(f) -e_0\|_{X^{\alpha}} \leq \varepsilon$. 
Moreover, the corresponding local unstable and local stable manifolds of $e(f)$ 
are given by
\begin{align*}
W^u_{loc}(e(f),f) &\equiv W^u(e(f),U_0,f) \\
&= \{ v \in U_0\, | \, v = e(f) + P_u(v-e(f)) + h_u(f)(P_u(v-e(f))) \} \\
W^s_{loc}(e(f),f)  &\equiv W^s_{loc}(e(f),U_0, f)\\
&= \{ v \in U_0 \, | \, v = e(f)+ P_s(v-e(f)) + h_s(f)(P_s(v-e(f))) \}~,
\end{align*}
where $h_u (f): P_u X^{\alpha}\to P_sX^{\alpha}$ and  $h_s(f) : P_s X^\alpha 
\to P_u X^{\alpha}$ are maps of class $\Cc^r$ such that $h_u (f)(0)=0$,
$h_s(f)(0) =0$ and $\| h_u(f) - h_u^0\|_{\Cc^r} \leq \varepsilon$ and $\| h_s(f)
- h_s^0\|_{\Cc^r} \leq \varepsilon$. Finally, for any $f \in \mathcal{N}$, the
above constants $k_i$, $\gamma_i$ are independent of $f$.
\end{theorem}
\begin{demo}
We refer to \cite[Theorems 5.2.1. and 5.2.2]{Henry} for the 
existence of the local stable and unstable manifolds in the case of a hyperbolic equilibrium
point of a parabolic equation. To obtain the last part of the Theorem, that is the smooth 
dependence with respect to $f$, we simply use a fixed point theorem with parameter. 
Indeed, the proof of Theorem 5.2.1 of \cite{Henry} consists in constructing the mappings $h_u$ and $h_s$
as fixed points of suitable contraction mappings. These maps depend smoothly on $f$ and thus 
remain contractions mappings for $f$ close to $f_0$ and their fixed points $h_u(f)$ and $h_s(f)$ 
depend smoothly on $f$. Notice that in general $Dh_u(f)(0)$ and $Dh_s(f)(0)$ do not vanish,
but are only small of order $\varepsilon$.
\end{demo}

Let $p(x,t)$ be a hyperbolic periodic solution of \eqref{eq} of
minimal period $\omega >0$, let $\Gamma=\{p(t)\, | \, t\in [0,\omega)\}$ be the 
associated orbit
and let $\Pi(t,0): X^{\alpha} \to X^{\alpha}$, be the associated evolution 
operator 
defined by the linearized equation \eqref{lin-periodic}.
We denote $\mu_i$, $i \in \NN$, the eigenvalues of the period map 
$\Pi(\omega,0)$. Since $p(x,t)$ is a hyperbolic periodic solution, the 
intersection of
the spectrum of $\Pi(\omega,0)$ with the unit circle $S^1$ of $\mathbb{C}$ 
reduces to
the eigenvalue $1$, which is a simple (isolated) eigenvalue. We recall
that, if $p(a)$, $a \in [0,\omega)$, is another point of the periodic
orbit, the spectrum of $D_u(S_{f}(\omega,0) p(a))$ coincides with the
one of $\Pi(\omega,0)$ whereas the corresponding eigenfunctions depend on
the point $p(a)$.

We  denote  $P_u(a)$ (resp. $P_c(a)$, resp. $P_s(a)$) the projection
in $X^{\alpha}$ onto the space
generated by the (generalized) eigenfunctions of $D_u(S_{f}(\omega,0)p(a))$
corresponding to the eigenvalues with
modulus strictly larger than $1$ (resp. equal to $1$, resp. with
modulus strictly smaller than $1$).

Since a hyperbolic periodic orbit is a particular case of a normally
hyperbolic $\Cc^1$ manifold, we may apply, for example, the existence
results of  \cite{Bates-Lu98}, \cite{HPS70}, \cite{HPS77} or
\cite[Theorem 14.2 and Remark 14.3]{Ruelle} and thus, we may state the 
following theorem.
Other methods of proofs are also given in  \cite{AbRobbin}, \cite{HPS77},
\cite{HMO84}, \cite{HVL} and \cite{Palis-de-Melo}.

\begin{theorem}\label{th-Wuper}
Let $f_0$ be given in $\Cg^r$, $r \geq 2$, and let $\Gamma_0=\{p_0(t)\, | \, 
t\in [0,\omega_0)\}$ be a hyperbolic periodic orbit of
Eq.  \eqref{eq} of minimal period $\omega_0>0$. 
\begin{enum-1}
\item There exists a small
neighborhood $U_{\Gamma_0}$ of $\Gamma_0$ in $X^{\alpha}$ such that the
local unstable and stable sets
\begin{equation*}
\begin{split}
&W^u_{loc}(\Gamma_0) \equiv W^u(\Gamma_0,U_{\Gamma_0}) = \{u_0 \in X^{\alpha}
\, | \, S_{f_0}(t)u_0 \in U_{\Gamma_0}\, , \, \forall t \leq 0\} \cr
&W^s_{loc}(\Gamma_0) \equiv W^s(\Gamma_0,U_{\Gamma_0}) = \{u_0 \in X^{\alpha}
\, | \, S_{f_0}(t)u_0 \in U_{\Gamma_0}\, , \, \forall t \geq 0\} 
\end{split}
\end{equation*}
are (embedded) $\Cc^1$-submanifolds of $X^{\alpha}$ of dimension $i(\Gamma_0)+1$
and
codimension $i(\Gamma_0)$ respectively. 

\item Moreover,
$W^{s}_{loc}(\Gamma_0)$ and $W^u_{loc}(\Gamma_0)$ are fibrated by the local
strongly stable (resp. unstable) manifolds at each point $p_0(a) \in \Gamma_0$, 
that is,
$$
W^{s}_{loc}(\Gamma_0) = \cup_{a \in [0,\omega_0)} W^{ss}_{loc}(p_0(a))~, \quad
W^{u}_{loc}(\Gamma_0) = \cup_{a \in [0,\omega_0)} W^{su}_{loc}(p_0(a))~,
$$
where there exist positive constants ${\tilde r}_0$, $\kappa_0$ and 
$\kappa^*_0$ such
that
\begin{equation}
\label{Dperaux1}
\begin{split}
W^{ss}_{loc}(p_0(a))= &\{u_0 \in X^{\alpha} \, | \, \|S_{f_0}(t) u_0 -
p_0(a+t)\|_{X^{\alpha}} < {\tilde r}_0~, ~\forall t \geq 0 \, ,\cr
& ~\lim_{t \to \infty} e^{\kappa_0 t} \|S_{f_0}(t) u_0 - p_0(a+t) 
\|_{X^{\alpha}}
=0\} ~, \cr
W^{su}_{loc}(p_0(a))=& \{u_0 \in X^{\alpha} \, | \, \|S_{f_0}(t) u_0 -
p_0(a +t)\|_{X^{\alpha}} < {\tilde r}_0 ~, ~ \forall t \leq 0 \, , \cr
&~\lim_{t \to  -\infty} e^{\kappa^*_0 t} \|S_{f_0}(t) u_0 - p_0(a+t) 
\|_{X^{\alpha}}
=0\}~.
\end{split}
\end{equation}
For any $a \in [0, \omega_0)$, $W^{su}_{loc}(p_0(a))$ (resp.
$W^{ss}_{loc}(p_0(a))$) is a $\Cc^{r}$-submanifold of $X^{\alpha}$ of dimension
$i(\Gamma)$
(resp. of codimension $i(\Gamma)+1$)
tangent at $p_0(a)$ to $P_u(a) X^{\alpha}$ (resp. $P_s(a) X^{\alpha}$). 

\item Finally, the local stable and unstable manifolds of the
periodic orbit continuously depend on the nonlinear map $f \in \Cg^r$.
\end{enum-1}
\end{theorem}

We have seen that the local stable and unstable manifolds are 
$\mathcal{C}^r$ graphs over $P_sX^{\alpha}$ and $P_uX^{\alpha}$ respectively.
In general, the global stable and unstable manifolds are not embedded
submanifolds of $X^{\alpha}$. 

Adapting the proof of  \cite[Theorem 6.1.9]{Henry}, one easily shows the 
following result.
\begin{theorem}\label{th-pWuL-WsL} Let $f \in \Cg^r$, $r \geq 2$, be given. 
\begin{enum-1}
\item Let $e_0$ be a hyperbolic equilibrium point of \eqref{eq}.
Then, the global unstable set $W^u(e_0)$ (resp. global stable set $W^s(e_0)$) 
is an injectively immersed invariant manifold of class $\mathcal{C}^r$ in 
$X^{\alpha}$ of dimension  (resp. of codimension) $i(e_0)$.

\item Likewise, let $\Gamma_0 =\{p_0(t) \, | \, t \in [0,\omega_0] \}$ be a
hyperbolic periodic orbit of minimal period $\omega_0 >0$. Then, 
the global unstable set $W^u(\Gamma_0)$ (resp. global stable set
$W^s(\Gamma_0)$) is an injectively immersed invariant manifold of class
$\mathcal{C}^r$ in $X^{\alpha}$ of dimension  $i(\Gamma_0) +1$ (resp. of
codimension $i(\Gamma_0)$). 
\end{enum-1}
\end{theorem}

\begin{demo} We will give the proof in the case of a hyperbolic equilibrium 
$e_0$, since the proof is very similar
 in the case of a hyperbolic periodic orbit.
 
{\noindent \it Proof for the unstable manifold:}  For every $m \in \NN$, we introduce the open
set
 $$
 U_0(m)= \{ x \in U_0 \, | \, S_f(t)x \hbox{ is well defined}, 0 \leq t \leq m 
\}~,
 $$
 where $U_0$ is the neighborhood of $e_0$, in which the local stable and 
unstable manifolds are given as graphs (see Theorem \ref{th-eWuL-WsL}). By 
Proposition \ref{Pr-exist}, $U_0(m)$ is an open subset of $U_0$ and thus 
$W^u_{loc}(e_0) \cap U_0(m)$ is an open subset of  $W^u_{loc}(e_0)$. We readily 
check that 
\begin{equation}
\label{global-unstable}
W^u(e_0) = \cup_{m=0}^{+ \infty} S_f(m) (W^u_{loc}(e_0) \cap U_0(m))~.
\end{equation}
Moreover, since $W^u_{loc}(e_0)$ is negatively invariant, we  have, for any $m 
\in \NN$,
$$
S_f(m) (W^u_{loc}(e_0) \cap U_0(m)) \subset S_f(m +1 ) (W^u_{loc}(e_0) \cap 
U_0(m +1))~.
$$
By Corollary \ref{coro-backward-unique-S}, $S_f(m)$ is an injective map from 
$U_0(m)$ into $X^{\alpha}$.
Moreover, by Proposition \ref{prop-backward-uniqueness}, for any $x \in U_0(m)$,
$D_uS_f(t)x$ is an injective map from $X^{\alpha}$ into itself, thus
$S_f(m)_{|U_0(m)}$ is an injective $\Cc^r$-immersion.
By Theorem \ref{th-eWuL-WsL}, $W^u_{loc}(e_0)$ is the image of an injective
$\Cc^r$-map $H_u$ from the open ball $B_{\RR^k}(0,1)$ of center $0$ and radius
$1$ of $\RR^k$ into $X^{\alpha}$, where 
$k=i(e_0)$. Moreover, the derivative $DH_u(y)$  has rank $k$ at each point $y 
\in B_{\RR^k}(0,1)$.
We recall that  $H^{-1}_u(W^u_{loc}(e_0) \cap U_0(m))$ is an open subset
$V(k,m)$ of $B_{\RR^k}(0,1)$. It follows that $S_f(m) W^u_{loc}(e_0) \cap
U_0(m))$ is the image of the injective $\Cc^r$-immersion $S_f(m) \circ H_u :
V(k,m) \to  X^{\alpha}$ and thus is a $\Cc^r$-submanifold of dimension $k$.
Since the invariance is obvious, Statement 1) is proved. 

{\noindent \it Proof for the stable manifold:} We first remark that 
\begin{equation}
\label{global-stable}
W^s(e_0) = \cup_{m=0}^{+ \infty} S_f(m)^{-1}(W^s_{loc}(e_0))~.
\end{equation}
Moreover, since $W^s_{loc}(e_0)$ is positively invariant, we  have, for any $m 
\in \NN$,
$$
S_f(m)^{-1}( W^s_{loc}(e_0)) \subset S_f(m +1)^{-1}( W^s_{loc}(e_0))~.
$$
As a consequence of the property \eqref{graphes} in Theorem \ref{th-eWuL-WsL}, 
where $h_s^0$ is a 
$\Cc^r$-map of $P_sX^{\alpha}$ into the $k$-dimensional space $P_uX^{\alpha}$
and where $Dh_s^0(0)=0$, 
$W^s_{loc}(e_0)$ is actually represented as the set $\{ v \in U_0 \, |\, 
g(v)=0\}$, where $g: x \in U_0
\mapsto g(x) \in \RR^k$ is a map of class $\Cc^r$ and $Dg(v)$ has constant rank
$k$ at every point 
$v \in g^{-1}(0)$. By \cite[Theorem 7.3.3]{Henry}, $DS_f(m)u$ has dense range 
at every point 
$u \in X^{\alpha}$ at which $S_f(m)u$ exists if $(DS_f(m)u)^*$ is injective. By 
Proposition \ref{prop-backward-uniqueness}, the adjoint equation 
\eqref{eq-adjoint} also satisfies the backward uniqueness property.
Thus $DS_f(m)u$ has dense range at every point $u \in S_f(-m)W^s_{loc}(e_0)$, 
which implies that, at every point 
$u \in (g \circ S_f(m))^{-1}(0)$, $D(g(S_f(m))u)$ has rank $k$. In other terms, 
 the mapping 
$v \to g(S_f(m)v)$ is a submersion of constant rank $k$  at every point $u \in(g
\circ S_f(m))^{-1}(0)$.  By a theorem on Page 12 of \cite{Marsden} for example,
$(g \circ S_f(m))^{-1}(0)$ is a $\Cc^r$-submanifold of $X^{\alpha}$ of
codimension $k$. Thus, since $S_{f(m)}$ is injective, $W^s(e_0)$ is an
injectively immersed manifold of codimension $k$. Since the invariance is
obvious, Statement 2) is proved. 
\end{demo}

\subsection{Transversality of connecting orbits}
We use here the above concepts of stable and unstable manifolds of hyperbolic equilibrium points
or periodic orbits. The definitions related to Theorem \ref{th-main} are as follows.
\begin{defi}
Let $\mathcal{C}^{\pm}$ be two hyperbolic critical elements. We say that
$W^u(\mathcal{C}^-)$ and $W^s(\mathcal{C}^+)$ {\bf intersect transversally} (or 
are {\bf transverse}) and we denote it by
$$
W^u(\mathcal{C}^-) \pitchfork W^s(\mathcal{C}^+)~,
$$
if, at each intersection point $u_0 \in W^u(\mathcal{C}^-) \cap 
W^s(\mathcal{C}^+)$, 
$T_{u_0}W^u(\mathcal{C}^-)$ splits, that is, contains
a closed complement of $T_{u_0}W^s(\mathcal{C}^+)$ in $X^{\alpha}$.
\end{defi}
It is important to notice that, in this paper, the complement of $T_{u_0}W^s$ 
in $X^\alpha$ is always closed since $T_{u_0}W^u(\mathcal{C}^-)$ is 
finite-dimensional. Also note that, by definition, manifolds which do not 
intersect are transverse.
\begin{defi}
Let $\Cc^-\neq \Cc^+$ be two different hyperbolic critical elements.
A trajectory $u(t)$ of $S(t)$ is a {\bf heteroclinic orbit}  connecting $\Cc^-$ 
to $\Cc^+$ if $u(t)\in W^u(\mathcal{C}^-) \cap W^s(\mathcal{C}^+)$.

Let $\Cc$ be a hyperbolic critical element.
A trajectory $u(t)$ of $S(t)$ is a {\bf homoclinic orbit} to $\Cc$ if $u(t)\in 
W^u(\mathcal{C}) \cap W^s(\mathcal{C})$.

A heteroclinic or homoclinic orbit is {\bf transverse} if the above 
intersections of stable and unstable manifolds are transverse.
\end{defi}


\section{Singular nodal sets for linear parabolic equations
with parameter}\label{section-nodalset}

In this section, we consider a general linear parabolic equation with parameter
\begin{equation}\label{eq-nodalset}
\partial_t v(x,t,\tau)=\Delta v(x,t,\tau) + 
a(x,t,\tau)v(x,t,\tau)+b(x,t,\tau).
\nabla_x v(x,t,\tau)~,
\end{equation}
in a domain $\Omega$ of $\RR^d$.

We are interested in the singular nodal set of $v$, that is the points
$(x,t,\tau)$ where $v$ and $\nabla_x v$ both vanish. To this end, we use
techniques coming from \cite{HardtSimon}. The singular nodal set of
solutions of the parabolic equations, with coefficients independent of the
parameter $\tau$, has already been studied in \cite{HanLin} and in \cite{Chen}. 
Notice that we
assume that $v$ is smooth in the variables $(x,t) \in \Omega \times \RR$, but 
this is not a restriction since this property holds in the applications, that 
we have in mind (see Section \ref{section-1to1}).

\begin{theorem}\label{th-nodalset}
Let $I$ and $J$ be open intervals of $\RR$. Let
$a\in\Cc^\infty(\Omega\times
I\times J,\RR)$ and $b\in\Cc^\infty(\Omega\times
I\times J,\RR^d)$ be bounded coefficients. Let $v$ be a strong solution
of \eqref{eq-nodalset} with Dirichlet boundary conditions. Let $r\geq 1$
and assume that $v$ is of class $\Cc^r$ with respect to $\tau$ and of class
$\Cc^\infty$ with respect to $x$ and $t$. Assume moreover that there are no time
$t\in I$ and no parameter $\tau\in J$ such that $v(.,t,\tau)
\equiv 0$. Then,
\begin{enum-1}
\item $M= \{(x,t,\tau) \in \Omega\times I\times J \, | \,
v(x,t,\tau)=0~,~\nabla_x v(x,t,\tau) =0\}$ 
is contained in a countable union of $\Cc^r-$manifolds of dimension $d$,
\begin{itemize}
\item  either parametrized by $t$, $\tau$ and $d-2$ components of $x$,
\item  or parametrized by $\tau$ and $d-1$ components of $x$.
\end{itemize}
\item the set
\begin{equation*}
(TNS)= \{ (x_0,t_0) \in \Omega \times I \, |\,\nexists \tau  \in J \text{ such 
that } (v(x_0,t_0,\tau),\nabla v(x_0,t_0,\tau) )=(0,0)\}
\end{equation*}
is generic in $\Omega\times I$.
\end{enum-1}
\end{theorem}

\begin{demo}
We introduce the set
\begin{equation*}
\begin{split}
M_q=\{(x,t,\tau) \in &\Omega \times I \times J \text{ such that for all } 
|\alpha| \leq q~,~~D^{\alpha}_x v(x,t,\tau)=0 ,\cr &\text{ and there exists }  
\alpha, \text{ so that }|\alpha|=q+1,~D^\alpha_x v(x,t,\tau)\neq 0 ~\}~.
\end{split}
\end{equation*}
By Proposition \ref{uniqcontB}, if $v(x,t,\tau)$
vanishes at infinite order in $x$, then $v(.,t,\tau)$ identically vanishes in
$\Omega$. By assumption, this is precluded. Thus, $M=\cup_{q\geq 1} M_q$.
And,  without loss of generality, we can replace $M$ by $M_q$ in Property 1) of
Theorem \ref{th-nodalset}. 

Let $q\geq 1$ and $(x_0,t_0,\tau_0)\in \Omega\times I\times J$. Let us first 
prove that there exists $\rho_{0,q}>0$ such that Property 1) of Theorem
\ref{th-nodalset} holds with $\Omega\times I\times J$
replaced by the ball $B((x_0,t_0,\tau_0),\rho_{0,q})$ and $M$ replaced by
$M_q$. 
Assume that $(x_0,t_0,\tau_0)\in M_q$ (otherwise the property is trivial).
There exists a multi-index $\beta$ with $|\beta| = q-1$ such that
$\hbox{Hess}(D^{\beta}_x v(x_0,t_0,\tau_0)) \neq 0$. In particular, there
exist $i, j$, $1 \leq i, j \leq d$, such that the derivative $D^2_{x_ix_j}
(D^{\beta}v(x_0,t_0,\tau_0)) \neq 0$. We next consider the $D^{\beta}_x$
derivative of the equation \eqref{eq-nodalset}. Since $v$ vanishes at order
$|\beta|+1$ at $(x_0,t_0,\tau_0)$, we obtain the equality
$$  \frac{d}{dt}D^{\beta}v(x_0,t_0,\tau_0) = \Delta_x
(D^{\beta}v(x_0,t_0,\tau_0))~.$$
Now two cases can occur:
\begin{itemize}
\item  Either $\frac{d}{dt}D^{\beta}v(x_0,t_0,\tau_0) =0$ and thus
$\sum_{k=1}^{d} \frac{\partial^2}{\partial x_k^2}
(D^{\beta}v(x_0,t_0,\tau_0))=0$. In this case, if $\frac{\partial^2}{\partial 
x_k^2}
(D^{\beta}v(x_0,t_0,\tau_0))$ $=0$ for all $k$, then
there exist $i \ne j$, such that
$D^2_{x_ix_j}(D^{\beta}v(x_0,t_0,\tau_0)) \ne 0$. By considering their
$i^{\text{th}}$ and $j^{\text{th}}$ components, we see that $\nabla_x
D_{x_i}(D^{\beta}v(x_0,t_0,\tau_0))$ and $\nabla_x
D_{x_j}(D^{\beta}v(x_0,t_0,\tau_0))$ are linearly independent.
If, on the contrary, there exists $i$ such that $\frac{\partial^2}{\partial
x_{i}^2} (D^{\beta}v(x_0,t_0,\tau_0)) \neq 0$, then there also exists
$j \neq i$ such that 
$$ \frac{\partial^2}{\partial
x_{i}^2} (D^{\beta}v(x_0,t_0,\tau_0)) \times
\frac{\partial^2}{\partial x_{j}^2} (D^{\beta}v(x_0,t_0,\tau_0))
< 0~.$$ 
By considering their ${i}^{\text{th}}$ and ${j}^{\text{th}}$ components, we
notice again that the vectors $\nabla_x 
D_{x_i}(D^{\beta}v(x_0,t_0,\tau_0))$ and $\nabla_x
D_{x_j}(D^{\beta}v(x_0,t_0,\tau_0))$ are
linearly independent. To summarize, in all the cases, there
exist $i$ and $j$, such that the vectors $\nabla_x
D_{x_i}(D^{\beta}v(x_0,t_0,\tau_0))$ and $\nabla_x
D_{x_j}(D^{\beta}v(x_0,t_0,\tau_0))$ are
linearly independent. 
This implies that there exists $\rho_{0,q}>0$ such that
$$ 
B((x_0,t_0, \tau_0),\rho_{0,q}) \cap (D_{x_i} D^{\beta}v)^{-1}(0)
\cap (D_{x_j}D^{\beta}v)^{-1}(0) 
$$ 
is an embedded $\Cc^r-$submanifold
$M_q(x_0,t_0,\tau_0)$ in $\RR^{d+2}$ of dimension $d$ which
contains all of $B((x_0,t_0, \tau_0),\rho_{0,q}) \cap M_q$.
This submanifold can be written as 
\begin{align*}
M_q(x_0,t_0,\tau_0)= & \left\{(x,t,\tau)\in B((x_0,t_0,
\tau_0),\rho_{0,q})\text{ such that }\right.\\
&~~~~~~~~\left. (x_i,x_j)=(\Phi_i ((x_k)_{k\neq i,j},t,
\tau), \Phi_j((x_k)_{k\neq i,j},t, \tau)) \right\}.
\end{align*}
\item Or $\frac{d}{dt}D^{\beta}v(x_0,t_0,\tau_0) \neq 0$, then
there exists $i$ such that $D_{x_i}^2D^{\beta}v(x_0,t_0,\tau_0) \neq 0$.
Notice that, since $D_{x_i}D^\beta v(x_0,t_0,\tau_0)=0$, 
$(D_{x_i},D_t)D^{\beta}v(x_0,t_0,\tau_0)$ and
$(D_{x_i},D_t)(D_{x_i} D^{\beta} v ( x_0, t_0,\tau_0))$ are linearly
independent. Thus, there exists $\rho_{0,q}>0$ such that
$$ 
B((x_0,t_0, \tau_0),\rho_{0,q}) \cap
(D_{x_i} D^{\beta}v)^{-1}(0) \cap (D^{\beta}v)^{-1}(0) 
$$
is an embedded $\Cc^r-$submanifold $M_q(x_0,t_0,\tau_0)$ in $\RR^{d+2}$
of dimension
$d$, which contains all of 
$B((x_0,t_0, \tau_0), \rho_{0,q}) \cap M_q$. This submanifold can be written as
\begin{align*}
M_q (x_0,t_0,\tau_0) =\big\{(x,t,\tau)&\in B((x_0,t_0,
\tau_0),\rho_{0,q})\text{ such that }\\
&(x_i,t)=(\Phi_i ((x_k)_{k\neq i},\tau), \Phi((x_k)_{k\neq i},\tau))\big\}~.
\end{align*}
\end{itemize}

To finish the proof of the first part of
Theorem \ref{th-nodalset}, notice that, since $\Omega\times I\times J$
is separable, for any $q\geq 1$, we can find a countable number of points
$(x_{n,q},t_{n,q},\tau_{n,q})_{n\geq 1}$ such that $\Omega\times I\times
J=\cup_n B((x_{n,q},t_{n,q},\tau_{n,q}),\rho_{n,q})$ and therefore we have
$M\subset \cup_{q\geq 1}\cup_{n\geq 1} M_{n,q}$ with
$M_{n,q}=M_q (x_{n,q},t_{n,q},\tau_{n,q})$.

Let $P:(x,t,\tau)\mapsto (x,t)$ be the canonical projection. Obviously,
$(TNS)$ is the complementary of $PM$. To prove the second part of  Theorem
\ref{th-nodalset}, it is thus sufficient to show that the projections of the
manifolds $M_{n,q}$ obtained above have an image
which is contained in a closed set of empty interior. For any $n$ and
$q$, $P_{|M_{n,q}}$ is a $\Cc^r-$ (and a fortiori a $\Cc^1-$) map defined
from a smooth manifold of dimension $d$ into $\Omega\times I\subset \RR^{d+1}$.
By the Sard theorem (see for example \cite[page 41]{AbRobbin}), the set of 
regular values
of this map is an open dense subset of $\Omega\times I$ (without loss of
generality, we may restrict the size of
$B((x_{n,q},t_{n,q},\tau_{n,q}),\rho_{n,q})$ in order to prove the openness
property). Obviously, the
derivative of $P_{|M_{n,q}}$ is never surjective and
thus the regular values of this projection map are not in its image. Hence,
$P(M_{n,q})$ is contained in a closed set of empty
interior, and property 2) of Theorem \ref{th-nodalset} follows from the
inclusion $M\subset \cup_{q\geq 1}\cup_{n\geq 1}M_{n,q}$.
\end{demo}

\begin{coro}\label{coro-nodalset}
Assume that the hypotheses of Theorem \ref{th-nodalset} hold. Assume
moreover that $a$ and $b$ and $v$ do not depend on $\tau$. Then the set
\begin{equation*}
(NS)= \{ x_0 \in \Omega \, | 
\hbox{ there does not exist
} t  \in I \hbox{ such that } (v(x_0,t),\nabla v(x_0,t))=(0,0)\}
\end{equation*}
is generic in $\Omega$.
\end{coro}
\begin{demo}
Since the problem is now independent of $\tau$, Property 1) of Theorem
\ref{th-nodalset} becomes: $M= \{(x,t) \in \Omega\times I \, | \,
v(x,t)=0~,~\nabla_x v(x,t) =0\}$ is contained in a countable union of manifolds
of dimension $d-1$, either parametrized by $t$ and $d-2$
components of $x$, or parametrized by $d-1$ components of $x$. Then, Corollary
\ref{coro-nodalset} follows from a use of the Sard theorem like in the proof of
Theorem \ref{th-nodalset}.
\end{demo}


\section{One-to-one properties for global solutions}\label{section-1to1}

In this section, we use the properties of the singular nodal sets of the
linearized
equation \eqref{eq-nodalset} of Section \ref{section-nodalset} in order  to 
prove
one-to-one properties for bounded complete solutions of the parabolic equation
\eqref{eq}. We recall that, in Section 2.4, we had deduced the backward 
uniqueness 
property of \eqref{eq} from the backward uniqueness property of the linearized
parabolic equation
\eqref{lineariz} with coefficients $a$ and $b$ given respectively by 
\eqref{eq-ab1} and \eqref{eq-ab2},
where  $u_1$ and $u_2$ are two solutions of \eqref{eq} (see the proposition 
\ref{prop-backward-uniqueness} and the corollary \ref{coro-backward-unique-S}).

Our first result concerns the periodic orbits $p$. It states that, for almost
every point $(x_0,t_0)\in \Omega\times\RR$, the value $(x_0,p(x_0,t_0),\nabla
p(x_0,t_0))$ is not taken twice during a period. Notice that if $\Omega$ is the
circle $S^1$, this property holds for all the points
$(x_0,t_0)$, see \cite{JR}. 
\begin{proposition}\label{prop-1to1-per}
Let $f\in\Cc^\infty(\overline\Omega\times \RR\times\RR^d,\RR)$. Let $p(t)$ be a
periodic solution of \eqref{eq} with minimal period $\omega>0$. Then there 
exists a
dense open set of points $(x_0,t_0)\in\Omega\times\RR$ such that
\begin{align*}
\text{i)}~~& (p_t(x_0,t_0),\nabla p_t(x_0,t_0))\neq (0,0)\\
\text{ii)}~~& (p(x_0,t_0),\nabla p(x_0,t_0)) \neq (p(x_0,t),\nabla p(x_0,t))
~\text{ if }~t\not\in t_0+\ZZ \omega 
\end{align*}
\end{proposition}
\begin{demo}
First, since $f$ is of class $\Cc^\infty$ and $p$ is a  bounded complete 
solution,
Proposition \ref{Pr-Cinfini} implies that 
$p\in\Cc^\infty(\Omega\times\RR,\RR)$. 
We already noticed that $p_t$ satisfies \eqref{lineariz} with coefficients $a$
and $b$ given by \eqref{eq-ab1}. Since $f$ and $p$ are of class  $\Cc^\infty$, 
the coefficients $a$ and $b$ are also of class  $\Cc^\infty$. Moreover, by 
Proposition 
\ref{prop-backward-uniqueness}, there exists no time $s$ such that $p_t(s) =0$. 
Thus, Corollary
\ref{coro-nodalset} implies that there is a generic set of points
$x_0 \in\Omega$ such that $(p_t(x_0,t),\nabla p_t(x_0,t))\neq (0,0)$, for any 
$t \in \RR$.

Next, we set $v(x,t,\tau)=p(x,t)-p(x,t+\tau)$, which solves
\eqref{lineariz} with coefficients given by \eqref{eq-ab2}. Again, we notice 
that 
$v$, $a$ and $b$ are infinitely differentiable with respect to $x$, $t$ and
$\tau$. Moreover, if there exist $t_1 \in\RR$ and $0 <\tau_1< \omega$ so that
$v(., t_1, \tau_1) \equiv 0$, then by the backward uniqueness property of 
Corollary
\ref{coro-backward-unique-S}, $v(., t,\tau_1) \equiv 0$, which means that 
$p(t)$ is periodic of period 
$\tau_1 < \omega$ and contradicts the fact that $\omega$ is the minimal period.
Thus, we can apply Theorem \ref{th-nodalset} to $v$ with $I=\RR$ and
$J=(0,\omega)$ to obtain a generic set of points $(x_0,t_0)\in\Omega\times\RR$ 
such
that the condition ii) holds. Therefore, both conditions i) and  ii) are
satisfied in a generic, and a fortiori dense, subset of $\Omega\times\RR$.

It remains to prove the openness. We consider the variable $t$ modulo the period
$\omega$, that is we work on $S=\RR/(\ZZ \omega)$. Let 
$(x_0,t_0)\in\Omega\times S$
satisfying i) and ii). There is an open neighborhood $\Uc$ of $(x_0,t_0)$ in
which i)
holds everywhere in $\Uc$. Moreover, since i) holds, we may assume that for any
$(x,t)$ and $(x,t')$ in $\Uc$, $t \ne t'$,  $(p(x,t),\nabla p(x,t)) \neq 
(p(x,t'),\nabla
p(x,t'))$.
The set of values $\{(p(x_0,t),\nabla p(x_0,t)), (x_0,t)\not\in \Uc\}$ is
compact and does not contain $(p(x_0,t_0),\nabla p(x_0,t_0))$ due to property
ii). Hence, this set of values is at positive distance of the value
$(p(x_0,t_0),\nabla p(x_0,t_0))$. Therefore, there exists a neighborhood
$\Vc\subset\Uc$ of $(x_0,t_0)$ such that, for any
$(x_1,t_1)\in \Vc$, $(p(x_1,t_1),\nabla p(x_1,t_1))$ is not contained in
$\{(p(x_1,t),\nabla p(x_1,t)), (x_1,t)\not\in \Uc\}$. This shows that ii) holds
in $\Vc$ and concludes the proof of the proposition.
\end{demo}

We also need to separate a periodic orbit from any other (bounded) complete 
solution.
\begin{proposition}\label{prop-1to1-up}
Let $f\in\Cc^\infty(\overline\Omega\times \RR\times\RR^d,\RR)$. Let $p(t)$ be
a periodic orbit of \eqref{eq} of minimal period $\omega$ or an equilibrium 
point, in which case we adopt the convention that   $p$ is a periodic solution 
with minimal period $\omega = 0$. 
Let $u(t)$ be a bounded complete solution of \eqref{eq}, such that, $p(t) \ne 
u(s)$, for any $(t,s) 
\in \RR^2$. Then there exists a
dense open set of points $(x_0,t_0)\in\Omega\times\RR$ such that
$(u(x_0,t_0),\nabla u(x_0,t_0)) \neq (p(x_0,t),\nabla p(x_0,t))$ for all
$t\in\RR$.
\end{proposition}
\begin{demo}
The proof is very similar to the one of Proposition \ref{prop-1to1-per} and thus
the details are left to the reader.  We emphasize only a few arguments. Since
$f$ is of class $\Cc^\infty$ and $u$, $p$ are bounded complete solutions,
Proposition \ref{Pr-Cinfini} implies that $p$ and $u$ belong to the space 
$\Cc^\infty(\Omega\times\RR,\RR)$. 
To prove the genericity of the points $(x_0,t_0)\in\Omega\times\RR$ such that
$(u(x_0,t_0),\nabla u(x_0,t_0)) \neq (p(x_0,t),\nabla p(x_0,t))$ for all
$t\in\RR$, we apply  Theorem \ref{th-nodalset} to 
$v(x,t,\tau)=u(x,t)-p(x,t+\tau)$, with
$I=J= \RR$. The function $v$ satisfies the hypotheses of  Theorem 
\ref{th-nodalset} and, in particular,  
due to the assumption of the proposition, there are no times $t$ and $\tau$ 
such that 
$v(., t,\tau) \equiv 0$. To show the openness of the set of the points 
$(x_0,t_0)\in\Omega\times\RR$ such that $(u(x_0,t_0),\nabla u(x_0,t_0)) \neq 
(p(x_0,t),\nabla p(x_0,t))$ for all
$t\in\RR$, one proceeds like in the proof of Proposition \ref{prop-1to1-per} by 
using the compactness of the set $\{(p(x_0,t),\nabla p(x_0,t)), t \in \RR\}$
 (but here the proof is even simpler, since we do not need to introduce the 
quotient $S$)
\end{demo}

As a particular case of the previous proposition, notice that we obtain the following 
result of separation of periodic orbits. In the case where $\Omega$ is the 
circle $S^1$, the arguments of \cite{Czaja-Rocha} show that this property holds 
for all the points $(x_0,t_0)$ (and not only for a dense open subset). The 
generalization to higher dimension is as follows.
\begin{proposition}\label{prop-1to1-CR}
Let $f\in\Cc^\infty(\overline\Omega\times \RR\times\RR^d,\RR)$. Let $p_1(t)$ and
$p_2(t)$ be two periodic solutions of \eqref{eq} of minimal periods $\omega_1$ 
and $\omega_2$. Assume that they do not correspond to the same 
periodic orbit, that is that $p_1(t)\neq p_2(s)$ for all $(t,s)\in\RR^2$. Then 
there exists
a dense open set of points $(x_0,t_0)\in\Omega\times\RR$ such that
$(p_1(x_0,t_0),\nabla p_1(x_0,t_0)) \neq (p_2(x_0,t),\nabla p_2(x_0,t))$ for all
$t\in\RR$.
\end{proposition}

The main dynamical result of this paper concerns heteroclinic and homoclinic 
orbits. We will need the following result.
\begin{proposition}\label{prop-1to1-hetero}
Let $f\in\Cc^\infty(\Omega\times \RR\times\RR^d,\RR)$. Let $p_-(t)$
and $p_+(t)$ be two periodic solutions of \eqref{eq} of minimal periods 
$\omega_-$ and
$\omega_+$ respectively. These periodic solutions may coincide or  each one may 
 be reduced to an equilibrium point, in which case we adopt the convention that 
the minimal period $\omega$ is equal to $0$.  Let $u(t)$ be
a global solution of \eqref{eq} connecting $p_-(t)$ and $p_+(t)$, that is, 
$$
u(t)-p_\pm(t)\xrightarrow[~t\longrightarrow \pm\infty~]{}0~.
$$
Then there exists a
dense open set of points $(x_0,t_0)\in\Omega\times\RR$ such that
\begin{align*}
\text{i)}~~& (\partial_t u(x_0,t_0),\nabla \partial_t u(x_0,t_0))\neq (0,0)\\
\text{ii)}~~& (u(x_0,t_0),\nabla u(x_0,t_0)) \neq (u(x_0,t),\nabla u(x_0,t))
~\forall~t\neq t_0\\
\text{iii)}~~&  (u(x_0,t_0),\nabla u(x_0,t_0)) \neq (p_\pm(x_0,t),\nabla
p_\pm(x_0,t)) ~\forall~t\in\RR
\end{align*}
\end{proposition}
\begin{demo}
Once again, the proof is very similar to the one of Proposition
\ref{prop-1to1-per}. We apply Theorem \ref{th-nodalset} to
$v(x,t,\tau)=u(x,t)-u(x,t+\tau)$ with $\tau<0$ and $\tau>0$ to prove the density of 
Property ii); and to $v(x,t,\tau)=u(x,t)-p_\pm(x,t+\tau)$ for the density of Property 
iii).
To prove the openness of Properties ii) and iii), we fix a point $(x_0,t_0)$ such that
i)-iii) hold. Due to i), there exists a neighborhood $\Uc=B(x_0,\rho)\times (t_0-\delta,t_0+\delta)$
of $(x_0,t_0)$ such that $(u(x,t),\nabla u(x,t))$ is injective in $\Uc$. Then we use the
compactness of $\{(u(x_0,t), \nabla u(x_0,t)) , t\in(-\infty,t_0-\delta]\cup [t_0+\delta,+\infty)\}
\cup \{(p_-(x_0,t), \nabla p_-(x_0,t)) ,t\in\RR\} \cup \{(p_+(x_0,t), \nabla p_+(x_0,t)), 
t\in\RR\}$ with arguments similar to the ones of the proof of Proposition
\ref{prop-1to1-per}.
\end{demo}


\section{Generic transversality of connecting orbits}\label{section-transverse}

To obtain the transversality of a connecting orbit as stated in Theorem 
\ref{th-main}, we need to show that we can perturb any parabolic semiflow 
$S_f(t)$ to another one, for which the considered stable and unstable 
manifolds intersect transversally. The construction of a
suitable perturbation $f+\varepsilon g$ of $f$ is the main difficulty in this 
task. Indeed, the global dynamical framework is classical and well understood 
in finite dimension. In Section \ref{section-dyn}, we have seen that the 
infinite dimension of $X^\alpha$ does not really affect this framework. The 
main novelty in this paper lies in the construction of a suitable perturbation 
$f+\varepsilon g$ of $f$ because we will need all the accurate PDE results 
proved in Sections \ref{section-nodalset} and \ref{section-1to1}.

\subsection{A perturbation to make an orbit transverse}

The first step consists in constructing a suitable perturbation $g$, which 
acts on a heteroclinic or homoclinic orbit $u(t)$ in a localized time 
interval only. In the following result, the one-to-one properties proved in Section 
\ref{section-1to1} are crucial. 
\begin{proposition}\label{prop-constr2}
Let $f \in \Cc^\infty(\overline\Omega\times\RR\times\RR^d,\RR)$ and let $u(t)$ 
be a bounded complete solution connecting $p_-(t)$ to $p_+(t)$ where $p_\pm(t)$ 
are two periodic solutions of minimal periods  $\omega_{\pm}$. 
Notice that $p_-=p_+$ is possible and that $p_\pm$ could be equilibrium points 
in which case we use the convention $\omega_\pm=0$.
Let $E$ be a compact subset of $\overline\Omega\times
\RR \times\RR^d$ with non empty interior, let $\Uc$ be an open subset of
$\Omega\times\RR$ and let $\psi\in\Cc^0(\Uc,\RR)$. Assume that there exists 
$(x_0,t_0)\in\Uc$ such that $(x_0,u(x_0,t_0),\nabla u(x_0,t_0))$ belongs to the 
interior of $E$ and $\psi(x_0,t_0)\neq 0$. 

Then, there exists a function
$h\in\Cc^\infty(\overline\Omega\times\RR\times\RR^d,\RR)$
such that 
\begin{enumerate}[(i)]
\item the function $h\,:\,\overline\Omega\times\RR\times\RR^d\rightarrow \RR$ 
has a compact support contained in $E$,
\item the function $h\circ u\,:\,(x,t)\in\overline\Omega\times\RR \longmapsto 
h(x,u(x,t),\grad u(x,t))\in\RR$ has a support contained in $\Uc$,
\item we have 
$\int_{\Omega\times\RR} \psi(x,t)h(x,u(x,t),\grad u(x,t))\,dxdt\,\neq\, 0$.
\end{enumerate}
\end{proposition}
\begin{demo}
Since $\psi(x_0,t_0)\neq 0$ and $(x_0,t_0)\in\Uc$, without loss of generality, 
by choosing $\Uc$ smaller, we may assume that $\psi$ does not vanish in $\Uc$. We set 
\begin{align*}
K=&\{(x,u(x,t),\grad u(x,t)),~(x,t)\not\in\Uc \}\cup \{(x,p_-
(x,t),\grad p_-(x,t)),~(x,t)\in\overline\Omega\times\RR \}\\
&~~\cup \{(x,p_+ (x,t),\grad p_+(x,t)),~(x,t)\in \overline\Omega\times\RR \}~.
\end{align*}
Proposition \ref{prop-1to1-hetero} shows that there is a dense open set of points 
$(\tilde x,\tilde t)\in\Uc$
such that $(\tilde x,u(\tilde x,\tilde t),\grad u(\tilde x,\tilde t))$ 
does not belong to $K$. Up to perturbing our reference point,
we can thus assume in addition that $(x_0,u(x_0,t_0),\nabla u(x_0,t_0))$ does not belongs to $K$.
Notice that $(x_0,u(x_0,t_0),\nabla u(x_0,t_0))$ still belongs to the interior of $E$ if our 
perturbation is small enough. Since $K$ is compact, $(x_0,u(x_0,t_0),\grad u(x_0,t_0))$ is in 
the interior of $E\setminus K$. Hence, we claim that it is sufficient to choose $h$ non-negative, 
with compact support in $E\setminus K$ and such that $h(x_0,u(x_0,t_0),\grad u(x_0,t_0))>0$. 

Property (i) holds by construction. For all $(x,t)\notin \Uc$,
$(x,u(x,t),\grad u(x,t))\in K$ and thus $h(x,u(x,t),\grad u(x,t))\equiv 0$, showing (ii).
Moreover, $\psi(x,t)h(x,u,\grad u)$ is not zero at $(x_0,t_0)$ and its sign
is constant in $\Uc$. These properties together with (ii) show that (iii) holds. 
\end{demo}

Using this perturbation $g$, we are able to perturb a non-transversal 
connecting orbit to a transversal one.
 \begin{proposition}\label{prop-transverse-dense}
Let $f_0  \in \Cc^\infty(\overline\Omega \times\RR\times\RR^d,\RR)$ 
 and let $\mathcal{N}_0$ be any small open neighborhood of $f_0$ in the 
$\Cg^r$-Whitney topology ($r\geq 2$). 
Let $\Gamma_{\pm} \equiv \{p_{\pm}(t) \, | \, t \in
[0,\omega_{\pm})\}$
 be two hyperbolic periodic orbits of minimal periods $\omega_{\pm} \geq 0$ of
$S_{f_0}(t)$, which may be not distinct and may be equilibrium points if 
$\omega_\pm=0$. 

Then there exists a function $f\in \mathcal{N}_0$ such that 
$\Gamma_-$ and $\Gamma_+$ are still hyperbolic periodic orbits for $S_{f}(t)$
and the unstable manifold $W^u(\Gamma_{-}, f)$ of $\Gamma_-$ intersects
transversally the local stable manifold $W^s_{loc}(\Gamma_+, f) \equiv 
W^s_{loc}(\Gamma_+, f_0)$ of $\Gamma_+$.
\end{proposition}
\begin{demo}
We will prove the  existence of a function $f \in \mathcal{N}_0$ satisfying the 
properties of Proposition \ref{prop-transverse-dense} by applying the 
transversal density Theorem \ref{th-AbRob} in Appendix 
\ref{section-Sard-Smale}. 

First, notice that the larger the regularity $r$ is, the more difficult is the 
result. Thus, without loss of generality we assume 
 $r > dim\, W^u(\Gamma_-) -codim\, W^s(\Gamma_+)$ in
the remaining part of the proof.

In what follows, $E$ will be a regular compact
subset of $\overline\Omega\times\RR\times\RR^d$ with non-empty interior. We
denote by $\Cc^r_0(E)$ the subset of
functions $g \in\Cc^r(\overline\Omega\times\RR\times\RR^d,\RR)$, which
identically vanish outside $E$; in fact, we identify $\Cc^r_0(E)$ with  the
space of functions in $\Cc^r(E,\RR)$, for which the first $r$ derivatives vanish
on $\partial E$. We recall that the topology induced in $\Cc^r_0(E)$ by the
Whitney topology coincides with the classical $\Cc^r$ topology and thus that
$\Cc^r_0(E)$ is actually a Banach space.

The proof splits in several steps.
\vskip 2mm

\noindent {\em First step: construction of particular
neighborhoods}\\
By theorems \ref{th-eWuL-WsL} and \ref{th-Wuper} and the remarks 
following both theorems, there exist two neighborhoods $\tilde\Nc_\pm$ of
$\Gamma_\pm$, for which the local stable and local unstable manifolds 
$W^s(\Gamma_{\pm}, \tilde \Nc_{\pm},f_0)$ and 
$W^u(\Gamma_{\pm}, \tilde \Nc_{\pm},f_0)$ of $\Gamma_\pm$ are well defined and
such that $\overline{\tilde\Nc_{-}} \cap \overline{\tilde\Nc_{+}} = \emptyset$ if $\Gamma_{+} \ne
\Gamma_{-}$. In the case where $\Gamma_{+} = \Gamma_{-}$, 
$\tilde \Nc_{+}= \tilde \Nc_{-}$ can be chosen so that $\overline
{W^u(\Gamma_{+},\tilde \Nc_{+},f_0)}\cap \overline{W^s(\Gamma_{+}, \tilde 
\Nc_{+},f_0)} = \Gamma_{+}$. 

We would like to perturb $f_0$ to deform the global unstable manifold $W^u(\Gamma_-,f_0)$ without changing the dynamics in $\tilde \Nc_\pm$. By construction, the part of $W^u(\Gamma_-,f_0)$ outside $\overline {\tilde \Nc_-} \cup \overline {\tilde \Nc_+}$ is a non-empty open subset of $W^u(\Gamma_-,f_0)$. The difficulty is that the nonlinearity $f$ sees the phase space $X^\alpha$ only through the projections by the evaluation map  
\begin{equation}\label{eq-evaluationmap}
{\rm Ev}~:~(x,\varphi) \in \Omega\times X^\alpha~\longmapsto~(x,\varphi(x),\nabla \varphi (x))\in\Omega\times\RR\times\RR^d~.
\end{equation}
We need to be sure that for all $u(t)$ connecting $\Gamma_-$ to $\Gamma_+$, not only $u(t)$ goes outside $\tilde\Nc_-\cup\tilde\Nc_+$ but also ${\rm Ev}(u(t))$ goes outside ${\rm Ev}(\tilde\Nc_-\cup\tilde\Nc_+)$.

The local unstable manifold $W^u(\Gamma_-, \tilde \Nc_-,f_0)$ is an embedded finite dimensional manifold and its boundary $\Sigma^u_{-}=\partial W^u(\Gamma_-, \tilde \Nc_-,f_0)$ is a compact set such
that, for all trajectory $\tilde u(t)$ belonging to the global unstable manifold
$W^u(\Gamma_-,f_0)\setminus\Gamma_-$, there exists a time $\tilde t_0\in\RR$ such that
$\tilde u(\tilde t_0)\in\Sigma^u_{-}$. Let $\sigma\in \Sigma^u_{-}$ and consider the trajectory $u_\sigma(t)=S_{f_0}(t)\sigma$, solution of \eqref{eq} with initial data $u_\sigma(t=0)=\sigma$ and nonlinearity $f=f_0$. For all $t<0$, $u_\sigma(t)$ belongs to the local unstable manifold $W^u(\Gamma_-,\tilde \Nc_-,f_0)$. Moreover, due
to Proposition \ref{prop-1to1-up}, there exists $(x_\sigma,t_\sigma)\in \Omega\times \RR_+$ such that $(u_\sigma(x_\sigma,t_\sigma),\nabla u_\sigma(x_\sigma,t_\sigma)) \neq (p_\pm(x_\sigma,t),\nabla p_\pm(x_\sigma,t))$ for all $t\in\RR$, or equivalently
$$\{(x_\sigma,u_\sigma(x_\sigma,t_\sigma),\grad u_\sigma(x_\sigma,t_\sigma))\}~\cap~ {\rm Ev}\big(\{x_\sigma\}\times (\Gamma_-\cup \Gamma_+)\big)~ = ~\emptyset~.$$
Since $\{(x_\sigma,u_\sigma(x_\sigma,t_\sigma),\grad u_\sigma(x_\sigma,t_\sigma)\}$ and $\{x_\sigma\}\times (\Gamma_-\cup \Gamma_+)$ are compact sets and since ${\rm Ev}$ is continuous because $X^\alpha$ is continuously embedded in $\Cc^1(\Omega)$, we can find $r_\sigma>0$ and $\rho_\sigma>0$ and neighborhoods $\Nc_{\sigma,\pm} \subset \tilde \Nc_{\pm}$ of $\Gamma_\pm$ in $X^\alpha$ such that 
$$\Uc_\sigma ~:=~B_{\Omega}(x_\sigma,r_\sigma)\times B_{\RR^{d+1}}((u_\sigma(x_\sigma,t_\sigma),\grad u_\sigma(x_\sigma,t_\sigma)),\rho_\sigma)$$
and $\Nc_{\sigma,\pm}$ satisfy 
$$\min~\big\{ \|\xi_1-\xi_2\|_{\Omega\times\RR^{d+1}}~|~\xi_1\in \Uc_\sigma ~\text{and}~ \xi_2\in {\rm Ev}\big(B_{\Omega}(x_\sigma,r_\sigma)\times (\Nc_{\sigma,-}\cup \Nc_{\sigma,+})\big\}~ > ~0~.$$
By continuity of ${\rm Ev}$ and of the flow $S_{f}(t)$ with respect to the initial
data and with respect to $f \equiv f_0 +g$, there are a neighborhood
$\Vc_\sigma$ of $\sigma$ in $X^{\alpha}$ and a neighborhood $\Wc_\sigma$ of $0$ in
$\Cg^r$ such that for any $\sigma' \in \Vc_\sigma$ and $g \in \Wc_\sigma$, the trajectory $S_{f_0 +g}(t)\sigma'$ has a projection ${\rm Ev}(\{x_0\}\times S_{f_0 +g}(t)\sigma')$ contained in $\Uc_\sigma$ for a non-empty open lapse of time. 

We can proceed as above for any point $\sigma \in \Sigma^u_-$. By compactness of $\Sigma^u_-$, it can be covered by a finite collection $\Vc_{\sigma_1}$,\ldots, $\Vc_{\sigma_N}$ of neighborhoods of points $\sigma_1$,\ldots, $\sigma_N$. We set $\Nc_\pm=\cap_n \Nc_{\sigma_n,\pm}$ and $E=\cup_n \overline{\Uc_{\sigma_n}}$. Notice that $E$ is a finite union of closed balls. Thus, $\Cc^r_0(E)$ is a well-defined Banach subspace of $\Cg^r$ and we set 
$\Wc^r=\cap_n \Wc_{u_n} \cap \Cc^r_0(E)$.\\[2mm]
To summarize, our construction satisfies the following properties (see Figure \ref{fig1}):
\begin{enumerate}
\item The neighborhoods $\Nc_\pm$ are small enough such that the local stable 
and local unstable manifolds $W^s(\Gamma_{+}, \Nc_{+}, f_0)$ and 
$W^u(\Gamma_{-}, \Nc_{-}, f_0)$ are well defined. Moreover, these local 
manifolds do not intersect if $\Gamma_{+} \ne \Gamma_{-}$, or have an 
intersection reduced to $\Gamma$ if $\Gamma_{+}=\Gamma_{-}=\Gamma$.
\item For any $f= f_0 +g$ where $g \in \Wc^r$ (in particular $g$ is supported in the set $E$), 
the flow $S_{f_0 +g}(t)$ is equal to the flow of
$S_{f_0}(t)$ in $\Nc_{\pm}$. In particular, we have $W^s(\Gamma_{+}, \Nc_{+},f_0) = 
W^s(\Gamma_{+}, \Nc_{+},f_0+g)$ and $W^u(\Gamma_{-}, \Nc_{-},f_0) = 
W^u(\Gamma_{-}, \Nc_{-},f_0+g)$ and the properties of 1. still hold when 
$f_0$ is perturbed to $f=f_0+g$.
\item For any $f= f_0 +g$ where $g \in 
\Wc^r$, for any global trajectory $u(t)=S_{f_0+g}(t)u(0)$
of the unstable manifold of $\Gamma_-$ ($\Gamma_-$ excluded), there exists $(x_0,t_0)\in\Omega\times\RR$ and $r>0$ such that for all $(x,t)\in B_{\Omega\times\RR}((x_0,t_0),r)$, $(x,u(x,t),\nabla 
u(x,t))$ belongs to the interior of $E$ (which is the set where 
the perturbations $g\in \Wc^r$ can be constructed) and not in ${\rm Ev}(\{x\}\times \Nc_\pm)$.
\end{enumerate}

\begin{figure}[htp]
\resizebox{15cm}{!}{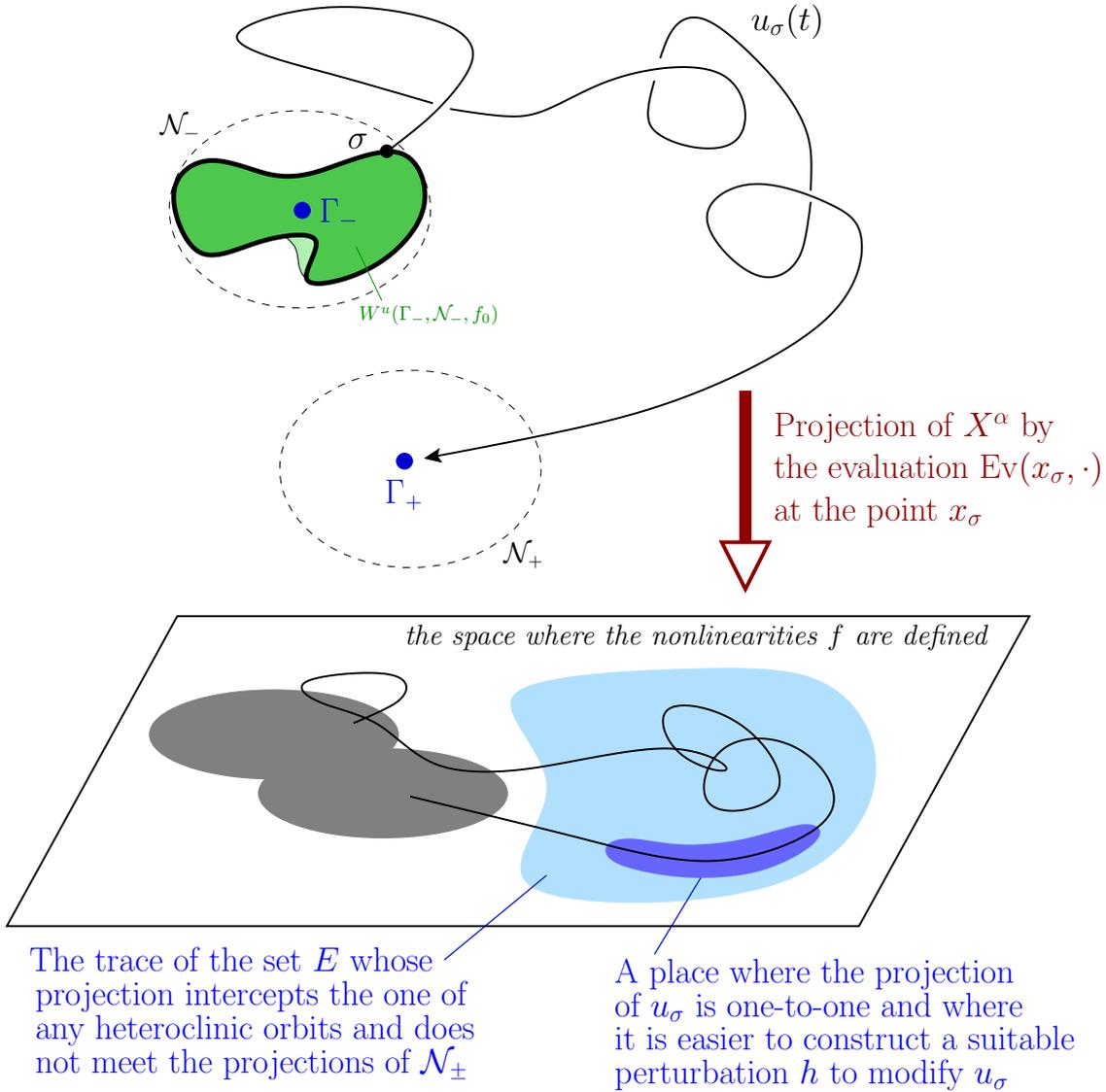}
\caption{\it A figure illustrating the proof of Proposition 
\ref{prop-transverse-dense}. In the phase space, $\Nc_\pm$ are small enough to 
define local dynamics and are disjoints in the heteroclinic case. The 
nonlinearity sees the dynamics only via the projections 
${\rm Ev}(x_\sigma,\cdot)$ by evaluating $(u_\sigma,\grad u_\sigma)$ at a point $x_\sigma$. 
In the first step, we construct a set $E$ whose projections do not meet the ones of the neighborhoods $\Nc_\pm$
of the closed orbits and such that, for all connecting orbit $u_\sigma$, there is a point $x_\sigma$ such that the evaluation of $u_\sigma(t)$ at this point enters in $E$ for an open lapse of times. The perturbation $g$ of the nonlinearity $f_0$ will be supported on this set $E$ to be able to modify any connecting orbits without modifying the closed orbits. Moreover, in the final step of our proof, we will also localize the perturbation in the place where the projection of $u_\sigma(t)$ has no self-intersection and where the modification of $u_\sigma(t)$ by a perturbation of the nonlinearity is easier to understand.   
\label{fig1}}
\end{figure}

\vspace{5mm}

\noindent {\em Second step: Application of the  Sard-Smale transversality 
Theorem \ref{th-AbRob}}\\
If $f= f_0 +g$, where $g$ is close to $0$ in $\Cc^r_0(E)$, then $f$ is close to 
$f_0$ in 
$\Cg^r$ (equipped with the Whitney topology). Moreover, by construction, for any
$f=f_0 + g$ with $g \in  \Wc^r$, $S_f(t)$ has the same dynamics as $S_{f_0}(t)$
in the neighborhoods $\Nc_{\pm}$ of $\Gamma_{\pm}$. Therefore, Proposition
\ref{prop-transverse-dense} holds if we can find a function $g \in \Wc^r \subset
\Cc^r_0(E)$ as close to $0$ as wanted such that $W^u(\Gamma_-, f_0 +g)$
intersects $W^s(\Gamma_{+}, \Nc_{+},f_0)$ transversally.

We recall that we did not assume global existence of solutions and thus the 
solutions in the unstable manifold may blow up. To overcome this technical 
problem, for all $m\geq 1$, we introduce the sets
$$\Nc^m_-=\{u_0\in\Nc_{-}~/~\forall g\in\Wc^r~,~S_{f_0+g}(t)u_0\text{ is well
defined for all }t\in[0,m]\}~.$$
The global orbit $\Gamma_-$ is obviously contained in $\Nc^m_-$ and we recall that ii) of Proposition \ref{Pr-exist} implies that $\Nc^m_-$ is open, in other words $\Nc^m_-$ is a neighborhood of $\Gamma_-$ contained in $\Nc_-$. Moreover, we have 
$$\forall g\in\Wc^r~,~~W^u(\Gamma_-,f_0+g)=\cup_{m \in\NN} S_f(m)W^u(\Gamma_-, 
\Nc^m_{-}, f_0+g)~.$$
To prove Proposition \ref{prop-transverse-dense}, it is sufficient to show that for any $m \in\NN$,
there exists a generic subset of functions $g\in\Wc^r$ such that
$S_{f_0+g}(m)W^u(\Gamma_-, \Nc^m_{-}, f_0)$ intersects $W^s(\Gamma_+, \Nc_{+},
f_0)$ transversally. Indeed the intersection of all these generic subsets is
generic and hence dense in $\Wc^r$ and consists in functions $f = f_0 +g$
such that $W^u(\Gamma_-, f_0 +g)$ intersects $W^s_{loc}(\Gamma_+,f_0 +g)$
transversally. 

To show this property, we are going to use the Sard-Smale transversality
theorem \ref{th-AbRob} in Appendix as follows. 
Let $m \geq 1$, let $\Mc=W^u(\Gamma_-,\Nc^m_-,f_0)$, $Y=X^\alpha$ and
$\Wc=W^s(\Gamma_+,\Nc_+,f_0)$. Let $\Lambda=\Wc^r$ and 
$\hat\Lambda=\Cc^\infty_0(E)\cap\Wc^r$. We define the mapping
$$\Phi~:~\left(\begin{array}{ccc} \Mc \times \Lambda & 
\longrightarrow & Y\\
(u_0,g) & \longmapsto & S_{f_0 +g}(m)u_0 \end{array}\right)$$
Notice that $S_{f_0 +g}(m)W^u(\Gamma_-,\Nc^m_-,f_0)$ intersects
$W^s(\Gamma_+,\Nc_+,f_0)$ transversally if and only if $\Phi(.,g)$ intersects
$W^s(\Gamma_+,\Nc_+,f_0)$ transversally. Thus, due to the above discussions, 
the conclusion of Theorem \ref{th-AbRob} in this framework will complete the 
proof of Proposition \ref{prop-transverse-dense}. Hypothesis i) of Theorem 
\ref{th-AbRob} is a consequence of the assumption $r > dim\, W^u(\Gamma_-) 
-codim\, W^s(\Gamma_+)$ made at the beginning of this proof and of the 
regularity of the parabolic flow with respect to the parameters. Thus, 
Hypothesis ii) is the only assumption which remains to be verified.

\vspace{5mm}

\noindent {\em Third step: checking Hypothesis ii) of Theorem \ref{th-AbRob}}\\
Let $u_0\in W^u(\Gamma_-,\Nc^m_-,f_0) \setminus \Gamma_{-}$ and 
$f= f_0 + g$, where $g \in\Wc^r$. If $S_f(m)u_0$
does not belong to $W^s(\Gamma_+,\Nc_+,f_0)$, then ii) is trivially satisfied.
If $S_f(m)u_0$ belongs to $W^s(\Gamma_+,\Nc_+,f_0)$, we set
$u(t)=S_f(t)u_0$ and we remark that, since $W^s(\Gamma_+,\Nc_+,f_0)=W^s(\Gamma_+,\Nc_+,f)$, $u(t)$ is a global solution and $u(t)\in W^s(\Gamma_+,\Nc_+,f)$ for all $t\geq m$. 

It remains to show that $\Phi$ is transversal to $\Wc$ in $X^\alpha$ at the point $u_0$, we have to compute 
$$D\Phi(u_0,g).(v_0,h)~=~ D_{u}\Phi(u_0,g).v_0~+~D_g \Phi(u_0,g).h~.$$
Let us consider the second term and let $v(t)$ be the derivative of $u(t)$ with respect to a variation $h$ of the nonlinearity $g$. By differentiating Equation \eqref{eq}, we have that $v$ solves
$$\partial_t v = \Delta v + h(x,u,\grad u) + f'_u(x,u,\grad u).v+f'_{\grad u}(x,u,\grad u).\grad v$$
with $v(t=0)=0$. 
We denote by $U(t,s)$ the family of evolution operators generated by the
equation \eqref{lineariz} with coefficients given by \eqref{eq-ab1}, which is
the linearization of the nonlinear equation along the trajectory $u(t)$. Using the variation of constants formula,
we get
\begin{equation}\label{calcul-derivee}
D_{g}\Phi(u_0,g).h = \int_0^{m} U(m,s)h(.,u(.,s),\grad u(.,s))~ds~.
\end{equation}
In a similar way, we obtain that $D_u\Phi(u_0,g).v_0= U(m,0)v_0$ whose range is the tangent space $T_{u(m)}W^u(\Gamma_-,f)$. 

We claim that the image of $D_g\Phi(u_0,g)$ is dense in $X^\alpha$ and we postpone the proof of this density in a final step below. Assuming this property, let us check Hypothesis ii) of Theorem \ref{th-AbRob} using Definition \ref{defi-trans-f}. First notice that $T_{u(m)}\Wc=T_{u(m)}W^s(\Gamma_+,\Nc_+,f)$ is a closed subspace with finite codimension (see Theorem \ref{th-eWuL-WsL}). To show that the image of $D\Phi(u_0,g)$ contains a closed complementary subspace of $T_{u(m)}\Wc$ in $X^\alpha$, it is sufficient to reach a given finite number of independent vectors $\phi_1$,\ldots, $\phi_p$ outside $T_{u(m)}\Wc$. This is obviously implied by the density of the image of $D_g\Phi(u_0,g)$ in $X^\alpha$. Since ${\rm span}(\phi_1,\ldots,\phi_p)\oplus T_{u(m)}\Wc=X^\alpha$, we have that $T_{u_0,g}\Mc\times \Lambda=D\Phi(u_0,g)^{-1}(T_{u(m)}\Wc)\oplus {\rm span}(\psi_1,\ldots,\psi_p)$ where $D\Phi(u_0,g).\psi_j=\phi_j$. By continuity, we directly have that $D\Phi(u_0,g)^{-1}(T_{u(m)}\Wc)$ is closed and its complementary space is also closed because of its finite-dimensionality.

\vspace{5mm}

\noindent {\em Fourth step: the image of $D_g\Phi(u_0,g)$ is dense in $X^\alpha$}\\
The operator $(-\Delta_D)^\alpha$ is a homeomorphism from $X^\alpha$ into $X$. Hence, it is sufficient to show that for any non-zero $\psi_{m}\in X^*$, there exists $h\in \Cc^\infty_0(E)$ such that
$$\langle\, \psi_{m} \,|\, (-\Delta_D)^\alpha D_{g} \Phi(u_0,g).h
\,\rangle_{X^*,X}\neq 0~.$$
Hence, using the expression of $D_{g} \Phi(u_0,g).h$ given by \eqref{calcul-derivee}, we have to find a function $h\in\Cc^\infty_0(E)$ such that
$$
\int_0^{m} \langle\, U(m,s)^* ((-\Delta_D)^\alpha)^*\psi_{m} \,|\,
h(.,u(.,s),\grad u(.,s)) \,\rangle_{X^{\alpha,*},X^\alpha}~ds \neq 0~.
$$
Now, we use Proposition \ref{prop-adjoint}: $\psi(s)=U(m,s)^*
((-\Delta_D)^\alpha)^* \psi_{m}$ is well defined in $X^*$ and is a solution
in $\Cc^0((0,m),\Cc^1(\overline\Omega))$ of \eqref{eq-adjoint} with $a$ and
$b$ as in \eqref{eq-ab1}. In particular, $\psi$ satisfies the unique
continuation property stated in Proposition \ref{prop-Fabre}: in any open set of $\Omega\times(0,m)$, there exists $(x,t)$ such that $\psi(x,t)\neq 0$. 

By considering the constructions made during the first step (see the third of the properties recalled at the end), we know that there exists a non-empty open set $\Uc\subset\Omega\times\RR$ such that for all $(x_0,t_0)\in\Uc$, $(x_0,u(x_0,t_0),\grad u(x_0,t_0))$ belongs to the interior of the set $E$ and is not in ${\rm Ev}(\{x\}\times \Nc_\pm)$. In particular, $u(x_0,t_0)$ cannot belongs to $\Nc_\pm$ and thus $t_0\in (0,m)$ because we have already noticed that $u(t)\in W^s(\Gamma_+,\Nc_+,f)$ for all $t\geq m$ and because $u(t)\in W^u(\Gamma_-,\Nc_-,f)$ for all $t\leq 0$ by definition of $\Phi$ and $u$. We now apply Proposition \ref{prop-constr2}, noticing that the unique continuation property for $\psi$ yields the existence of $(x_0,t_0)\in\Uc$ such that $\psi(x_0,t_0)\neq 0$. We obtain a function $h \in\Cc^\infty_0(E)$ such that
$$\int_\RR \int_\Omega \psi(x,s)h(x,u(x,s),\grad u(x,s)) \,dxds~\neq~ 0~.$$
It remains to notice Proposition \ref{prop-constr2} guarantees that $h\circ u$ is supported in $\Uc$ and that the above discussion shows that $\Uc\subset \Omega\times (0,m)$. Thus, for any $\psi_m\in X^*$, we may replace the domain $\RR\times\Omega$ by $[0,m]\times\Omega$ in the above integral and, in conclusion, we have obtained $h$ such that
$$\langle\, \psi_{m} \,|\, (-\Delta_D)^\alpha D_{g} \Phi(u_0,g).h
\,\rangle_{X^*,X}~=~\int_0^m \langle\, \psi(s) \,|\, h(.,u(.,s),\grad u(.,s)) \,\rangle_{X^{*},X}\,ds~\neq~ 0~.$$
which implies that the image of $D_g\Phi(u_0,g)$ is dense in $X^\alpha$.
\end{demo}

\subsection{Proof of Theorem \ref{th-main}}
The proof of our main theorem easily follows from the perturbation result of 
Proposition \ref{prop-transverse-dense}.

Let $f_0\in\Cg^r$ be given and let $\Cc_0^\pm$ be two hyperbolic critical 
elements. By Theorems \ref{th-stab-hyperbo}, \ref{th-eWuL-WsL} and 
\ref{th-pWuL-WsL}, there exists a neighborhood $\Oc$ of $f_0$ such that 
$\Cc_0^\pm$ are associated with two families $\Cc^\pm(f)$ of hyperbolic critical 
elements depending smoothly on $f$. Moreover, the corresponding local stable 
and unstable manifolds $W^u_{loc}(\Cc^-(f))$ and $W^s_{loc}(\Cc^+(f))$ also 
depend smoothly on $f$. 

Let $m\in\NN$ be given and let
\begin{align*}
W^u_m(\Cc^-(f))=&\{u\in X^\alpha\text{ such that }\|u\|_{X^\alpha}< m\text{ and 
there exists } \\ &  t\in [0,m]\text{ and }~u_0\in W^u_{loc}(\Cc^-(f))\text{ 
such that }u=S_f(t)u_0\}~.
\end{align*}
The set $W^u_m(\Cc^-(f))$ is a bounded open subset of the global unstable 
manifold $W^u(\Cc^-(f))$ and an immersed manifold of $X^\alpha$. Also notice 
that $W^u_m(\Cc^-(f))$ depends smoothly on $f$. We consider the sets
$$\Gg_m=\{ f \in \Oc~|~W^u_m(\Cc^-(f))\pitchfork W^s_{loc}(\Cc^+(f))\}~.$$
The smooth dependences yield that $\Gg_m$ are open subsets of $\Oc$ (see 
Appendix \ref{section-whitney} to understand what these smooth dependences 
mean with respect to the Whitney topology). We claim that the sets $\Gg_m$ are 
also 
dense. Indeed, $X^\alpha$ is embedded in $\Cc^1$ and so its ball 
$\{u~|~\|u\|_{X^\alpha}\leq m\}$ provides values $(x,u(x),\grad u(x))$ uniformly 
bounded by some constant $C(m)$. For any $f\in\Oc$, we may perturb $f$ to 
$\tilde f$ such that $\tilde f$ is of class $\Cc^\infty$ in the ball of radius 
$C(m)$ and equal to $f$ outside the ball of radius $C(m)+1$. In this way, 
$\tilde f$ is as close as wanted to $f$ in the $\Cg^r$ Whitney topology. Moreover,  
any solution $u$ in $W^u_m(\Cc^-(\tilde f))$ stays in the place where $\tilde f$ 
is a $\Cc^\infty-$non-linearity. 
Applying Proposition \ref{prop-transverse-dense}, we may perturb $\tilde f$ to 
obtain a non-linearity in $\Gg_m$.

Since the sets $\Gg_m$ are open and dense in $\Oc$, by setting $\Gg=\cap_m 
\Gg_m$, we obtain the generic set of Theorem \ref{th-main}.


\section{Further generalizations of the generic transversality stated in Theorem
\ref{th-main}}\label{section-beyond}

Our above arguments are not exactly specific to Equation \eqref{eq}. We may 
easily check the following generalizations.

\vspace{5mm}

\noindent {\bf Other geometries}\\
Dirichlet boundary conditions are not mandatory, we may choose Neumann ones or 
Robin ones. We may also consider other flat geometries such as $\Omega$ being a torus 
or a cylinder.

We may also add coefficients to the Laplacian operator $\Delta$, typically 
considering the Laplace-Beltrami operator $\frac 1 {\sqrt{g}} 
\text{div}(\sqrt{g}g_{ij} \grad\cdot)$ associated to a metric $g$. However, 
notice that part of our results, e.g. Theorem \ref{th-nodalset}, require smooth 
coefficients and thus $g$ needs to be smooth. Thus, we may generalize Theorem 
\ref{th-main} in the case where $\Omega$ is a bounded $\Cc^\infty-$submanifold 
of $\RR^n$, as a sphere for example.

\vspace{5mm}

\noindent {\bf Systems of parabolic equations}\\
Instead of considering the scalar parabolic equation \eqref{eq}, we consider a
system of $n$ parabolic equations as follows.  
We keep the same space $X =L^p(\Omega)$, $p>d$,  and the same $\Delta_D$
Laplacian operator with homogeneous Dirichlet boundary conditions. Like in the
introduction, we keep $\alpha \in (1/2 + d/2p, 1)$, so that $X^{\alpha}=
D((-\Delta_D)^{\alpha} \hookrightarrow
W^{2\alpha,p}(\Omega)$ is compactly embedded in $\Cc^1(\overline{\Omega})$.  Let
$n \in \NN$, $n \geq 1$. We consider the system of parabolic equations
\begin{equation}
\label{eqSn}
\left\{\begin{array}{ll}U_t(x,t)=\Delta U(x,t)+ F(x,U(x,t),\nabla
U(x,t)),&\quad (x,t)\in \Omega\times (0,+\infty)\\ 
U(x,t)=0,& \quad (x,t)\in\partial\Omega\times (0,+\infty)\\ 
U(x,0)=U_0(x) \in X_n^\alpha \equiv (X^\alpha)^n  , 
\end{array}\right.
\end{equation}
where $F \equiv (f_1, f_2, \ldots,f_n)
\in\Cc^r(\overline\Omega\times\RR^n\times\RR^{nd},\RR^n)$,
$r\geq 2$, and where $U \equiv (u_1(x,t), u_2(x,t),\ldots, u_n(x,t))$ belongs to 
$\RR^n$. As in the case $n=1$, the system 
\eqref{eqSn} generates a local dynamical system $S_n(t) \equiv S_{n,F}(t)$ on
$X_n^{\alpha}$.  This (local) dynamical system $S_{n,F}(t)$ satisfies all the
smoothing properties of Section \ref{section-basics} as well as the dynamical 
systems properties given in Section \ref{section-dyn}. The strong unique 
continuation property of Proposition \ref{uniqcontB} still holds and is proved 
in \cite[Theorem 2.2]{Chen} (see also \cite{HanLin}). 
The singular nodal sets properties as given in Theorem \ref{th-nodalset} and its
Corollary \ref{coro-nodalset} are still true and are proved with the same 
arguments (see also \cite[Theorem 2.3]{Chen}). These facts allow us to 
generalize Theorem \ref{th-main} to the system \eqref{eqSn}.

\vspace{5mm}

\noindent {\bf Genericity for other topologies}\\
We have chosen here to consider the genericity in $\Cg^r$ by endowing 
$\Cc^r(\overline\Omega\times\RR\times\RR^d,\RR)$ with the Whitney topology (see 
the 
precise definition in Appendix \ref{section-whitney}). Indeed this topology 
seems to be the most usual one for this kind of question concerning generic 
dynamics. Moreover, it also seems to be the most delicate topology since it has only a
few nice properties (for example the closed sets are not the sequentially closed 
sets and in particular $\Cg^r$ is not a metric space). However, Theorem 
\ref{th-main} also holds if we endow 
$\Cc^r(\overline\Omega\times\RR\times\RR^d,\RR)$ with other reasonable topology. 
We may for example consider $\Cc^r_b(\overline\Omega\times\RR\times\RR^d,\RR)$, 
the set of bounded $\Cc^r-$functions on $\overline\Omega\times\RR\times\RR^d$ 
endowed with the supremum $\Cc^r$-norm. We may also extend the previous metric 
by considering unbounded $\Cc^r-$functions but defining their neighborhoods with 
bounded perturbation only (in other words, we may say that if $f-g$ or one 
of its $r$ first derivatives is unbounded, then $f$ and $g$ are at 
infinite distance). In any case, the conclusions of Theorem \ref{th-main} remain 
valid since, in the proofs, we in fact only consider non-linearities via a 
bounded set of $\overline\Omega\times\RR\times\RR^d$, where all these 
topologies are equivalent (see Appendix \ref{section-whitney}).

\vspace{5mm}

\noindent {\bf Some open problems}\\
To conclude, let us mention cases where the generalization is not 
straightforward and remains an open problem.

We may wonder if Theorem \ref{th-main} is still true for systems of parabolic 
equations if, instead of considering mappings $F(x,U,\nabla U)$ in the set
$\Cc^r(\overline\Omega\times\RR^n\times\RR^{nd},\RR^n)$, one considers only
mappings $F(x,U)\in \Cc^r(\overline\Omega\times\RR^n,\RR^n)$ depending only on 
$x$ and of the value of $U$.
Since the Hausdorff dimension of the nodal set is, in general, larger by $1$
than the dimension of the singular nodal set (see a simple example in 
\cite[Section 9]{Chen}), the one-to-one properties of
global trajectories, as given in Section \ref{section-1to1}, can be false if 
$F=F(x,U)$ and are no longer consequences of 
Theorem \ref{th-nodalset} (see \cite[Section 9]{Chen}).

We can also wonder if one can extend Theorem
\ref{th-main} to the case where the Laplacian operator is replaced by a 
$2m$-th order homogeneous elliptic linear operator. In this case, in Equation 
\eqref{eq}, we replace the non-linearity $f(x,u,\nabla u)$ by a non-linearity 
$f^*(x,u, D_x u, D_x ^2u, ...., D_x^{2m-1}u)$ depending on the values of 
$u$, $D_xu$, ....,$D_x^{2m-1}u$. 
If the strong unique continuation property of Proposition \ref{uniqcontB} holds,
then, arguing exactly as in the proof of Theorem \ref{th-nodalset}, one shows
that the statement of this theorem is still true provided we
replace the singular nodal set by 
\begin{equation*}
\begin{split}
(TNS)= \{ (x_0,t_0) \in \Omega \times I \, | 
&\hbox{ there does not exist
} \tau  \in J \hbox{ such that }\cr
& (v,D_x v,D_x^2v,...., D_x^{2m-1}v)(x_0,t_0,\tau)=(0,0,....,0)\}
\end{split}
\end{equation*}
Unfortunately, the strong unique continuation property for the 
parabolic equation with higher order elliptic operators is not always true 
(concerning the elliptic equation, see \cite{HanHardtLin} and \cite{Plis} for example). 
For this reason, we cannot state here a generalization of Theorem \ref{th-main} for higher-order 
parabolic equation.


\setcounter{section}{0}
\renewcommand{\thesection}{{\Alph{section}}}

\section{Appendix: The Whitney topology}\label{section-whitney}

If we want to prove generic properties for the parabolic equation
\eqref{eq} with respect to the non-linearity $f$, we need to equip the space
of nonlinear functions $f$ with a topology. Let $E\subset \RR^n$, $n\geq 1$,
by $f\in\Cc^r(E,\RR)$, we mean that $f$ is $r$ times differentiable
in the set $E$ and that these derivatives are continuous. We do 
not a priori endow $\Cc^r(E,\RR)$ with any topology and we do not assume that 
$f$ or its derivatives are bounded. 

In this article, we consider $E=\overline\Omega\times\RR\times\RR^d$ which is 
unbounded in $\RR^{2d+1}$. Since we do not want to exclude unbounded 
non-linearities, we cannot equip $\Cc^r(E,\RR)$ with the classical 
$\Cc^r$-topology. 
\begin{defi}
For any $r\in\NN$, we denote by $\Cg^r \equiv
\Cg^r(E,\RR)$ the space $\Cc^r(E,\RR)$ endowed 
with the Whitney topology, that is the topology generated by the
neighborhoods
\begin{equation*}
\{g\in \Cc^r(E,\RR)~|~|D^i f(y)-D^i
g(y)|\leq\delta(y)~,~~\forall i\in\{0,1,\ldots,r\}~,~~\forall
y\in E\}~,
\end{equation*}
where $f$ is any function in $\Cc^r(E,\RR)$ and $\delta$ is any positive 
continuous function. 
\end{defi}
 
We emphasize that, if $E$ is bounded, then the Whitney topology 
coincides with the classical $\Cc^r$-topology and thus 
$\Cg^r(E,\RR)$ is a Banach space equipped with the 
classical norm $\|f\| = \sup_{i=0, 1,\ldots,r}\| f^{(i)}\|_{L^{\infty}}$. 
However, if $E=\overline\Omega\times \RR \times \RR^d$, the neighborhoods of a 
function $f$ in the Whitney topology cannot be
generated by a countable number of them. As a consequence, this topology is not
metrizable and open or closed sets cannot be characterized by sequences. In 
order to give an idea about the uncountable conditions imposed by the Whitney 
topology, we recall that a sequence of
functions $(f_n)$ converges to a function $f$ in the Whitney topology if
and only if there is a compact set $K \subset E$ such that $f_n 
\equiv f$ in $E\setminus K$ for any $n\in\NN$, but for a finite number of
them, and such that $(f_n)$ converges to $f$ in  the space $\Cc^r (K,\RR)$,
equipped with the classical topology of uniform convergence of the functions
together with their derivatives up to order $r$. This means that the Whitney 
topology imposes an uncountable number of conditions of proximity outside 
compact sets and thus a sequence has to be constant there in order to 
be convergent. 

As already written in Section \ref{section-beyond}, we could have chosen a simpler topology, 
but the Whitney topology seems to be the most usual one. In order 
to overcome several technical problems due to this topology, we make more precise some 
arguments in this appendix. We omit the corresponding problems during the main 
proofs of this paper to avoid too heavy proofs. However, if all the technical details are written, 
the interested reader will notice that we easily deal with the fact that the Whitney topology 
does not generate a Banach space as follows.

\vspace{5mm}

\noindent {\bf Genericity and Baire property:~}The main purpose of this paper is 
to obtain the genericity of the transversality of heteroclinic and homoclinic 
orbits. The notion of generic sets, that are sets containing a countable 
intersection of dense open sets, is important because it provides a nice notion 
of large subset. However, the acceptance of this notion is mainly related to the Baire 
property, that is the fact that the countable intersection of generic sets is 
generic. A space satisfying the Baire property is called a Baire space. Complete 
spaces, and in particular Banach ones, are Baire spaces. But when $E$ is 
unbounded, $\Cg^r(E,\RR)$ with its Whitney topology is even not metrizable. 
Thus, it is important to emphasize that it is at least a Baire space, implying 
that the genericity is still a meaningful concept (see \cite{GG} or 
\cite{Hirsch} for example). 

\vspace{5mm}

\noindent {\bf Smooth dependences, open or dense subsets and other abuses of 
notations:} When $E$ is unbounded, since $\Cg^r(E,\RR)$ is not metrizable, we 
can speak about continuous dependence on $f\in\Cg^r(E,\RR)$ but not about 
smooth dependence, even not about derivatives with respect to $f$. We sometimes 
use the following abuse of notation. Consider $K$ a compact subset of $E$ and 
define $P$ as the canonical projection from $\Cg^r(E,\RR)$ onto $\Cg^r(K,\RR)$, 
that is $Pf:=f_{|K}$ is the restriction of $f$ to $K$. Now, as already noticed, 
$\Cg^r(K,\RR)$ endowed with the Whitney topology is equivalent to the Banach 
space 
$\Cc^r(K,\RR)$ endowed with the classical $\Cc^r-$norm. Consider a function 
$\Phi$ depending on $f$ via the values in $K$ only. We may thus associate 
with $\Phi$ defined in $\Cg^r(E,\RR)$ a function $\tilde\Phi$ defined in 
$\Cg^r(K,\RR)$ and then it is relevant to say that $\tilde\Phi$ depends 
smoothly on $Pf$. In this case, we may use an abuse of notations by saying that 
$\Phi$ depends smoothly on $f$ instead of saying that $\tilde\Phi$ depends 
smoothly on $Pf$ (notice that, rigorously, we should not even say that $Pf$ 
depends smoothly on $f$). 

At this point, it is important to notice that, the restriction operator 
$$P:\Cg^r(E,\RR)\rightarrow \Cg^r(K,\RR)~\text{ with }K\subset E \text{ 
compact and }E\subset\RR^n$$ 
is continuous, open and surjective. Continuity is clear and surjectivity 
follows from the Whitney extension theorem (see \cite{AbRobbin}), or 
a simpler result if $r=0$ or $K$ is a regular subdomain for which the 
extension is easily constructed. Openness follows from the following argument: 
consider $g\in\Cg^r(K,\RR)$ close to $0$, extend $g$ to $f\in\Cg^r(E,\RR)$ and 
truncate $f$ by multiplying it by a smooth function $\chi$ with $0\leq \chi \leq 
1$, $\chi_{|K}\equiv 1$ and $\chi\equiv 0$ outside a small neighborhood of $K$. 
This provides a function $\chi f\in \Cg^r(E,\RR)$ with $P(\chi f)=g$ and $\chi 
f$ as close to $0$ in $\Cg^r(E,\RR)$ as wanted as soon as $g$ is small enough. 
Thus, the image by $P$ of any neighborhood of $0$ contains a neighborhood of 
$0$.

The surjectivity of $P$ enables to define the above functional $\tilde\Phi$ in 
$\Cg^r(K,\RR)$ because to each function $g\in\Cg^r(K,\RR)$ indeed corresponds 
a class of equivalence of functions $f\in\Cg^r(E,\RR)$ with $Pf=g$. The 
openness is useful to show that a property is open in $\Cg^r(E,\RR)$ if this 
property depends on the value of $f$ in $K$ only: if the property is open in 
$\Cg^r(K,\RR)$ with the above abuse of notation, then it is open in  
$\Cg^r(E,\RR)$. Together, these properties show that, with the abuse of 
notation, if a property is open and dense (resp. generic) in $\Cg^r(K,\RR)$ 
then it is open and dense (resp. generic) in $\Cg^r(E,\RR)$.

\vspace{5mm}

Notice that the above tricks have already been widely used in
previous articles (see \cite{Bruno-Pola} for instance).
Finally, for a further study of the Whitney topology and the comparison with
the weak topology, we refer the reader to \cite{GG} or  \cite{Hirsch} for
example.


\section{Appendix: Sard Theorem and Sard-Smale trans\-ver\-sa\-li\-ty 
theorems}\label{section-Sard-Smale}
The Sard theorem (\cite{Sard}) and the transversality theory (which goes back
to Thom \cite{Thom}) are very useful tools for proving the genericity of a
given property in finite dimension. In \cite{Smale-trans}, Smale has shown how
to use Fredholm theory to generalize the transversality theorems to
infinite-dimensional Banach spaces. There exist different version
of this kind of transversality theorems (often called Sard-Smale theorems or 
Thom theorems) with slight changes in the hypotheses, depending on the  
framework, in which they are used. We recall here the general framework 
and the version used in this paper.

Let $\Mc$ and $\Nc$ be two differentiable Banach manifolds and let
$f:\Mc\longrightarrow \Nc$ be a differentiable map. We say that 
$x\in\Mc$ is a regular point of $f$ if $Df(x):T_x\Mc \rightarrow T_{f(x)}\Nc$ is
surjective and its kernel splits (that is, has a closed complement in
$T_x\Mc$). A point $y\in \Nc$ is a
regular value of $f$ if any $x\in \Mc$ such that $f(x)=y$ is a regular point of
$f$. The points of $\Nc$ which are not regular values are said critical values.
The classical theorem of Sard is as follows.
\begin{theorem}\label{th-sard}
If $U$ is an open set of $\RR^p$ and if
$f:U\longrightarrow \RR^q$ is of class $\Cc^s$ with $s>max(p-q,0)$,
then, the set of critical values of $f$ in $\RR^q$ is of Lebesgue
measure zero.
\end{theorem}

Using Fredholm operators and a Lyapounov-Schmidt method,
Smale has generalized Sard Theorem to infinite-dimensional spaces (for 
introduction to Fredholm operators, see \cite{Bonic} for example). As a 
consequence of Smale theorem in \cite{Smale-trans}, many versions of Sard-Smale 
theorems can be obtained, see \cite{AbRobbin} and \cite{Henry05} for examples. 
The versions involving a functional formulation have been used since the pioneer work of 
Robbin \cite{Robbin} and are very useful in the PDE context where the 
geometrical arguments may be too difficult to perform, see Theorem \ref{th-sard-smale} below
and \cite{Bruno-Pola,Bruno-Raugel,Joly,JR,JR2}.
In this article, the transversality of connecting orbits may be proved with a 
more geometrical version of Sard-Smale theorems. Indeed, we only need to 
perturb an unstable manifold, which is finite-dimensional, and we may do it far 
from the periodic orbit, so that the basic framework does not depend on the 
parameter (see Section \ref{section-transverse}). This kind of geometrical 
setting is more difficult to use if we want to prove generic hyperbolicity as 
discussed in Appendix \ref{section-hyperbolicity} below.

We recall the following definition (see \cite{AbRobbin} for more details).
\begin{defi}\label{defi-trans-f}
Let $\Mc$ and $\Nc$ be two $\Cc^1$ Banach manifolds and let
$f\in\Cc^1(\Mc,\Nc)$. Let $\Wc$ be a $\Cc^1$ submanifold of $\Nc$. The function
$f$ is said to be {\bf transversal} to $\Wc$ at a point $x\in \Mc$ if either
$f(x)\not\in\Wc$ or $f(x)\in\Nc$ and
\begin{enum-i}
\item $D_xf^{-1}(T_{f(x)}\Wc)$ is a closed subspace of $T_x\Mc$ which admits a
closed complementary space,
\item $D_x f(T_x \Mc)$ contains a closed complement to $T_{f(x)}\Wc$ in
$T_{f(x)}\Nc$.
\end{enum-i}
\end{defi}

\noindent We need in this article a slight improvement of Theorem 19.1 of \cite{AbRobbin}.
The idea of replacing the condition on $\Lambda$ by a condition on a dense subset $\hat\Lambda$ 
only has been already used in \cite{Bruno-Pola,Bruno-Raugel,Joly} for example.
\begin{theorem}\label{th-AbRob}
Let $r\geq 1$. Let $\Mc$ be a $\Cc^r$ separable manifold of dimension $n$. Let 
$\Wc$ be a $\Cc^r$ manifold of codimension $m$ in a Banach space $Y$. Let
$\Lambda$ be an open subset of a separable Banach space and let $\hat\Lambda$
be a dense subset of $\Lambda$. Let $\Phi \in\Cc^r(\Mc\times\Lambda,Y)$. Assume
that
\begin{enum-i}
\item $r>n-m$,
\item $\Phi$ is transversal to $\Wc$ at any point $(x,\lambda)\in
\Mc\times\hat\Lambda$.
\end{enum-i}
Then, there is a generic set of parameters $\lambda\in\Lambda$ such that the
map $x\mapsto \Phi(x,\lambda)$ is everywhere transversal to $\Wc$.
\end{theorem}
\begin{demo}
Theorem \ref{th-AbRob} is proved as Theorem 19.1 of
\cite{AbRobbin}. The only difference is that hypothesis ii) is assumed here only
for a dense set of parameters $\lambda$. To obtain this improvement from the
classical version where ii) is assumed everywhere, we argue as follows. Since
$\Mc$ is separable and finite dimensional, we can find a countable sequence of
open subsets $(\Mc_k)$ such that $\Mc=\cup \Mc_k$ and $\overline\Mc_k$ is
contained in $\Mc$ and is compact. Let $\lambda_0\in\hat\Lambda$. Let
$(\lambda_p)$ be a sequence converging to $\lambda_0$. Assume that there is a
point $x_p\in \overline\Mc_k$ such that $\Phi$ is not transversal to $\Wc$ at
$(x_p,\lambda_p)$. By the compactness property, one may assume that $(x_p)$ 
converges to
$x_0\in\overline\Mc_k$. Since $\Phi$ is $\Cc^1$, $\Phi$ is not transversal to
$\Wc$ at $(x_0,\lambda_0)$ which is absurd. Thus, there exists a neighborhood
$\Uc$ of $\lambda_0$ such that ii) holds for any $(x,\lambda)\in \Mc_k
\times\Uc$. By applying \cite[Theorem 19.1]{AbRobbin}, we obtain a generic
subset $\Uc_k\subset \Uc$ such that for any $\lambda\in\Uc_k$, the map $x\mapsto
\Phi(x,\lambda)$ is transversal to $\Wc$ for any $x\in \Mc_k$. Since
$\hat\Lambda$ is dense in $\Lambda$, we have a generic
subset $\tilde\Uc_k\subset \Lambda$ such that for any $\lambda\in\tilde\Uc_k$,
the map $x\mapsto \Phi(x,\lambda)$ is transversal to $\Wc$ for any $x\in
\Mc_k$. The generic set of parameters appearing in the conclusion of Theorem
\ref{th-AbRob} is then $\cap_k \tilde\Uc_k$.
\end{demo}

\begin{figure}[ht]
\begin{center}
\resizebox{12cm}{!}{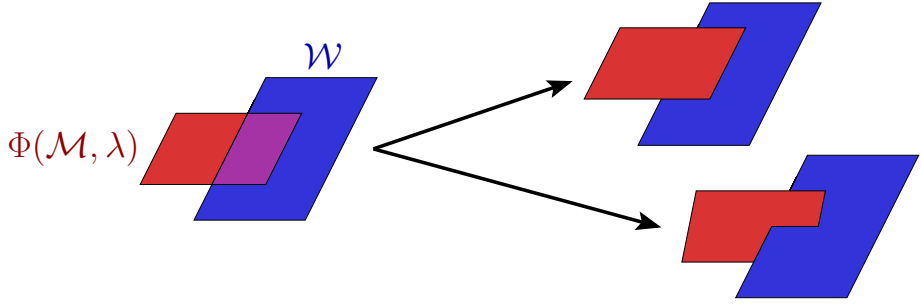}
\end{center}
\caption{\it The geometric idea behind Sard-Smale theorems as Theorem \ref{th-AbRob}: 
if perturbing the parameter $\lambda$ provides enough freedom, 
a non-transversal intersection between $\Phi(\Mc,\lambda)$ and $\Wc$ 
is generically perturbed into either an empty, and thus transversal, 
intersection or a non-empty transversal intersection.\hfill} \label{fig-sard-smale}
\end{figure}

For brief discussions in Appendix \ref{section-hyperbolicity} and for the curious reader,
we finish by a brief recall of one of the simplest version of Sard-Smale theorem with a functional
formulation (see for example \cite{Henry05} for other versions or proofs). Let us recall that a 
continuous linear map $f:E\longrightarrow F$ between two Banach spaces is a Fredholm map 
if its image is closed and if the dimension of its kernel and the codimension of its image are finite.
\begin{theorem}\label{th-sard-smale}
Let $k\geq 1$ and let $\Mc$, $\Nc$ and $\Lambda$ be three $\Cc^k$ Banach manifolds.
Let $y\in \Nc$ and let
$\Phi\in \Cc^k(\Mc\times\Lambda,\Nc)$. Assume that:
\begin{enum-i}
\item for any $(x,\lambda)\in \Phi^{-1}(\{y\})$, $D_x
\Phi(x,\lambda):T_x\Mc\rightarrow T_y\Nc$ is a Fredholm map of index $i$
strictly less than $k$,
\item for any $(x,\lambda)\in \Phi^{-1}(\{y\})$, $D\Phi(x,\lambda):T_x\Mc\times
T_\lambda \Lambda \rightarrow T_y\Nc$ is surjective,
\item $\Mc$ is separable.
\end{enum-i}
Then, there is a generic set of parameters $\lambda\in\Lambda$ such that 
for all $x\in\Mc$ such that $(x,\lambda)\in \Phi^{-1}(\{y\})$, $D_x\Phi(x,\lambda)$ is
surjective. 
\end{theorem}
As in Theorem \ref{th-AbRob}, a similar result holds if $\Lambda$ is replaced by a dense subset 
$\hat \Lambda \subset \Lambda$ and if $\Lambda$ is separable (see \cite{Bruno-Pola}).


\section{Appendix: discussion about proving the generic 
hyperbolicity of periodic orbits}\label{section-hyperbolicity}

The purpose of this section is unusual. To obtain the genericity of the Kupka-Smale 
property for the parabolic equation \eqref{eq}, it remains to prove the 
genericity of hyperbolicity of equilibrium points and periodic orbits. The generic 
hyperbolicity of equilibrium points is proved in \cite{JR}. We tried to obtain the 
generic hyperbolicity of periodic orbits but failed to get a complete proof. In this section, we 
would like to present some ideas and to point out where there is still a gap in the proof. 
Maybe this discussion could inspire a motivated reader.

The first proofs of generic hyperbolicity of periodic orbits appeared in 
\cite{Kupka,Smale63}. Peixoto in \cite{Peix66} introduced a nice recursion argument, 
which has been modified in \cite{AbRobbin} and \cite{Mallet}. Basically, the 
recursion is as follows. We introduce the sets
$$
{\cal G}_1(K)=\{f \in \Cg^r \, | \,\hbox{ any
equilibrium point }e\hbox{ of \eqref{eq} with
}\|e\|_{X^\alpha}\leq
K \hbox{ is hyperbolic}\}
$$
\begin{equation*}
\begin{split}
{\cal G}_{3/2}&(A,K) = \{f \in {\cal G}_1(K)\, | \, \hbox{ any non-constant
periodic solution } p(t) \hbox{ of \eqref{eq} } \cr
& \text{with period }T \in (0,A]
\hbox{ such that sup}_{t\in\RR} \|p(t)\|_{X^\alpha} \leq K
 \hbox{ is non-degenerate} \}~.
\end{split}
\end{equation*}
and 
\begin{equation*}
\begin{split}
{\cal G}_{2}(A,K) =&\{f \in {\cal G}_1(K) \, | \, \hbox{ any non-constant
periodic
solution } p(t) \hbox{ of \eqref{eq}} \cr
& \text{with period } T \in (0,A]
\hbox{ such that sup}_{t\in\RR} \|p(t)\|_{X^\alpha} \leq K
\hbox{ is hyperbolic} \}~.
\end{split}
\end{equation*}
The slightly strange above notation comes from the fact that ${\cal G}_1$ and
${\cal G}_2$ are the sets originally introduced by Peixoto, whereas the set 
${\cal G}_{3/2}$ has been introduced later.

We know from the arguments of the second part of Section 3 of \cite{JR} that 
$\Gc_1(K)$ is a dense open subset of $\Cg^r$. The idea of the recursion 
argument is that there exists $\varepsilon>0$ small enough, such that ${\cal 
G}_{2}(\varepsilon,K)={\cal G}_1(K)$ due to the
absence of periodic orbits of small period. Then, the method of Peixoto would 
consist in proving, like in \cite{Mallet}, that 
$\mathcal{G}_2(A,K) \cap \mathcal{G}_{3/2}(3A/2,K)$ is dense in
$\mathcal{G}_2(A,K)$ and that $\mathcal{G}_{2}(3A/2,K)$ is dense in
$\mathcal{G}_{3/2}(3A/2,K)$. 
By this way, we obtain a chain of dense inclusions\\[2mm]
\resizebox{\textwidth}{!}{\makebox{
$\ldots \mathcal{G}_{2}(9\varepsilon/4,K)
\underset{\text{\footnotesize dense}}{\text{\large $\subset$}}
\mathcal{G}_{3/2}({9\varepsilon}/4,K)
\underset{\text{\footnotesize dense}}{\text{\large $\subset$}}
\mathcal{G}_{2}({3\varepsilon}/2,K)
\underset{\text{\footnotesize dense}}{\text{\large $\subset$}}
\mathcal{G}_{3/2}({3\varepsilon}/2,K)
\underset{\text{\footnotesize dense}}{\text{\large $\subset$}}
{\cal G}_{2}(\varepsilon,K)={\cal G}_1(K)$}}\\[2mm]
which shows the density of the hyperbolicity of periodic orbits in 
${\cal G}_1$. The openness of these sets is rather simple
and similar to the finite-dimensional case considered in \cite{Peix66}. 
This scheme of proof has been exactly performed in \cite{Mallet} and in
\cite{AbRobbin}. The difficulty lies in the proofs of density.

\vspace{5mm}

We claim that the following density holds.
\begin{proposition} \label{generic3A/2A}
For any positive $A$ and $K$, ${\cal G}_{3/2}(3A/2,K)\cap {\cal G}_{2}(A,K)$ is
dense in ${\cal G}_{2}(A,K)$.
\end{proposition}
\begin{demo}
We give here very brief arguments since this proposition is only an 
auxiliary result in the whole proof of generic hyperbolicity, which is 
unfortunately not yet completed.

The proof of Proposition \ref{generic3A/2A} is very similar to the one of 
Proposition \ref{prop-transverse-dense}. We apply a suitable version of 
Sard-Smale theorem (similar to Theorem \ref{th-sard-smale}) to the map
$$ \Phi~:~(T,u_0,g) ~ \longmapsto ~ S_{f_0+g}(T)u_0 - u_0~.$$
As usual, the main difficulty is to obtain a surjectivity 
as required by Hypothesis ii) of Theorem \ref{th-sard-smale}. We skip the details, but 
simply notice that checking this property is very similar to the end of the 
proof of Proposition \ref{prop-transverse-dense}: we have to find for any 
solution $\varphi^*$ of the adjoint equation along a periodic orbit $p$, a 
perturbation $g$ of $f$ such that
\begin{equation*}
\int_{\Omega}\int_{0}^{T} g(x, p(x,s),\grad p(x,s)) \varphi^*(x,s) ds dx  \ne 
0~.
\end{equation*}
This is achieved by constructing a function as in Proposition 
\ref{prop-constr2} by using Proposition \ref{prop-1to1-per}.
\end{demo}

The proof of the genericity of the Kupka-Smale property would be obtained if we 
could prove the following result.

\begin{conj} \label{genericAA}
For any $A>0$ and $K$, ${\cal G}_{2}(3A/2,K)$ is
dense in ${\cal G}_{3/2}(3A/2,K)\cap {\cal G}_{2}(A,K)$.
\end{conj}
To prove this conjecture, we only need to know how to make hyperbolic a given 
simple periodic orbit in the following sense.
\begin{conj}\label{prop-hyperbo}
Let $f  \in \Cc^\infty(\Omega\times\RR\times\RR^d,\RR)$ 
and let $\Nc$ be any small open neighborhood of $f$ in $\Cc^r$.
Let $p$ be a simple periodic solution of \eqref{eq} with minimal period $\omega
>0$ and such that $\sup_{t \in [0,\omega]} \| p(t)\|_{X^\alpha}
\leq \tilde K$, where $\tilde K>0$. Then, there exists a function $\tilde f \in\Nc$ 
such that $p$ is a hyperbolic periodic solution of \eqref{eq} with non-linearity
$\tilde f$. 
\end{conj}
Once again, the usual strategy would be to apply a Sard-Smale theorem 
(similar to Theorem \ref{th-sard-smale}) to an 
appropriate functional $\Phi$ and then to check a surjectivity hypothesis as 
ii) of Theorem \ref{th-AbRob}. If we try the most natural way, we will have to 
find a perturbation $g$ of $f$ satisfying
\begin{equation}\label{eq-appendix-conj}
\Re \int_0^\omega \int_\Omega 
\,(D_u g,D_{\grad u} g)(x,p(x,t),\grad 
p(x,t))\,.\,\psi^*(x,t)(\phi,\grad \phi)(x,t)~dx dt
~\neq~ 0
\end{equation}
where $p$ is the considered simple periodic orbit, $\phi$ a solution of the 
linearized equation associated to an eigenvalue $\lambda$ with modulus 
$|\lambda|=1$ and $\psi^*$ a solution of the adjoint equation. 
Notice in \eqref{eq-appendix-conj} the presence of the real part $\Re$ since 
the spectrum of a periodic orbit has complex eigenvalues.
To obtain this perturbation $g$, we may use a construction as follows.
\begin{proposition}\label{prop-constr}
Let $f\in\Cc^\infty(\overline\Omega\times\RR\times\RR^d,\RR)$ and let 
$p \in \Cc^\infty(\Omega\times\RR,\RR)$
be a periodic solution of \eqref{eq} with minimal period $\omega$. Let
$V\in\Cc^\infty(\Omega\times [0,\omega],\RR^{d+1})$ be a function, which
 is not everywhere colinear to $(p_t(x,t),\grad p_t(x,t))$. Then, there exists 
a function $g\in\Cc^\infty(\overline\Omega\times\RR\times\RR^d,\RR)$ such that
\begin{align*}
\text{i)}~~& g(x,p(x,t),\grad p(x,t))=0~~~\forall (x,t)\in\Omega\times\RR,\\
\text{ii)}~~& \int_0^\omega \int_\Omega (D_ug,D_{\nabla u}g)(x,p(x,t),\grad
p(x,t)).V(x,t)~dxdt\neq 0~.
\end{align*}
\end{proposition}
\begin{demo}
To simplify the notations, we denote by $U$ the variable $(u,\grad
u)\in\RR^{d+1}$ and we set $P(x,t)=(p(x,t),\grad p(x,t))\in\RR^{d+1}$.

By assumption, there is an open set $\Uc$ with 
$\overline{\Uc} \subset \Omega\times (0,\omega)$ such that
$V$ is never colinear to $P_t$ on $\Uc$. Notice that, in particular
$P_t(x,t)\neq 0$ for all $(x,t) \in \Uc$. Due to Proposition
\ref{prop-1to1-per},  restricting
$\Uc$, we can assume that, for all $(x_0,t_0)\in\Uc$, the map
$(x,t) \in \Omega\times [0,\omega) \mapsto (x,P(x,t))\in\Omega\times\RR^{d+1}$ 
reaches
the value $(x_0,P(x_0,t_0))$ at $(x_0,t_0)$ only.

Let $(x_0,t_0)\in\Uc$. We complete $(P_t,V)$ to a basis of
$\RR^{d+1}$: let $W_1$,...,$W_{d-1}$ be $d-1$ vectors of
$\RR^{d+1}$ such that ($P_t(x_0,t_0)$, $V(x_0,t_0)$, $W_1$, $\ldots$,
$W_{d-1}$) is a basis of $\RR^{d+1}$. Restricting again $\Uc$, we can
assume that ($P_t(x,t)$, $V(x,t)$, $W_1$, $\ldots$, $W_{d-1}$) is a
basis of $\RR^{d+1}$ for all $(x,t)\in\Uc$. Let $\Vc=\Uc\times\Wc$
where $\Wc\subset\RR^d$ is a neighborhood of $0$. We define $h:\Vc \rightarrow
\Omega\times\RR^{d+1}$ by
$$h(x,t,\tau,s_1,\ldots,,s_{d-1})= \left(x,P(x,t) + \tau V(x,t) +
s_1W_1+...+s_{d-1}W_{d-1}\right)~.$$
Up to choosing $\Vc$ smaller, the local inversion theorem shows that $h$
is a $\Cc^\infty$-diffeomorphism into its
image. We recall that for all $(x_0,t_0)\in\Uc$, the map $\Omega\times
[0,\omega)\ni (x,t)\mapsto (x,P(x,t))\in\Omega\times \RR^{d+1}$ takes
the value $(x_0,P(x_0,t_0))$ at $(x_0,t_0)$ only. Due to the
compactness of the graph of this map, we can restrict $\Wc$ such
that $(x,P(x,t))$ belongs to $h(\Vc)$ if and only if $(x,t)$
belongs to $\Uc$. Let
$\chi\in\Cc^\infty(\Omega\times\RR^{d+1},\RR)$ be a function with compact 
support in 
 $\Vc$, which will be made more precise later.
We set 
$\theta(x,t,\tau,s_1,\ldots,,s_{d-1})=\chi(x,t,\tau,s_1,\ldots,,s_{d-1})\tau$. 
We define the
function $g: h(\Vc)\rightarrow \RR$ by $g(x,u,\grad u)=g(x,U)=\theta\circ
h^{-1}(x,U)$. We can extend $g$ by $0$ outside $h(\Vc)$ to obtain a function in
$\Cc^\infty(\overline\Omega\times\RR^{d+1})$. By construction, for all $(x,t) 
\notin \Uc$, 
$g(x,P(x,t))=0$ and $D_U g(x,P(x,t))=0$. Moreover, for all $(x,t)\in\Uc$,
$g(x,P(x,t))=\theta(x,t,0,0,...,0)=0$ and
\begin{align*}
\partial_U g(x,P(x,t)).V(x,t)&=D\theta(h^{-1}(x,P(x,t))).\left(\partial_U
h^{-1}(x,P(x,t)).V(x,t)\right)\\
&=D\theta(x,t,0,...,0).\left(\partial_U h^{-1}(h(x,t,0,...,0)).\partial_\tau
h(x,t,0,...,0)\right)\\
&=D\theta(x,t,0,...,0).\partial_\tau (h^{-1}\circ h) (x,t,0,...,0)\\
&=\partial_\tau \theta(x,t,0,...,0)\\
&=\chi(x,t,0,...,0)
\end{align*}
Thus, Property i) of Proposition \ref{prop-constr} holds and moreover
$$\int_0^\omega\int_\Omega \partial_U g(x,P(x,t)).V(x,t)~dxdt=\int_\Uc
\chi(x,t,0,...,0)~dxdt~.$$
Therefore, we can easily choose $\chi$ such that Property ii) of Proposition
\ref{prop-constr} also holds. \end{demo}

The final problem lies in checking that the real part of  
$\psi^*(x,t)(\phi,\grad \phi)$ in \eqref{eq-appendix-conj} 
is not everywhere colinear to $(p_t,\grad p_t)$. This is true if we only 
consider real functions (see Proposition \ref{prop-non-colinear} 
below), but we consider here complex solutions $\psi^*$ and $\phi$ and thus the 
real part of $\psi^*(\phi,\grad \phi)$ correspond to a combination of two real 
solutions of the linearized equation: the real and the imaginary parts of 
$\phi$. Even if 
this colinearity would be very strange and holds surely in very rare cases only 
(remember that we may break potential symmetries by perturbing $f$), we found no 
rigorous argument to avoid it. 

We finish with a statement of non-colinearity which could be inspiring.
\begin{proposition}\label{prop-non-colinear}
Let $I$ be an open interval of $\RR$ and $\Omega$ and open subset of $\RR^d$. 
Let $a\in\Cc^\infty(\Omega\times I,\RR)$ and $b\in\Cc^\infty(\Omega\times
I,\RR^d)$ be bounded coefficients. Let $v_1$ and 
$v_2$ be two solutions of the real equation
\begin{equation}\label{eq-nodalset2}
\partial_t v(x,t)=\Delta v(x,t) + a(x,t)v(x,t)+b(x,t).
\nabla_x v(x,t)~.
\end{equation}
Assume that $(v_1,\grad v_1)$ is colinear to $(v_2,\grad v_2)$ at each points 
$(x,t)$, meaning that there exists real values $\alpha(x,t)$ and $\beta(x,t)$ 
such that for all $(x,t)\in\Omega\times I$,
\begin{equation}\label{eq-nodalset3}
\alpha(x,t)(v_1,\grad v_1)(x,t)+\beta(x,t)(v_2,\grad v_2)(x,t)=0.
\end{equation}
Then $v_1$ and $v_2$ are colinear to $v_2$ as solutions, that is that 
\eqref{eq-nodalset3} holds with real constants $\alpha$ and $\beta$.
\end{proposition}
\begin{demo}
If $v_i\equiv 0$ for $i=1$ or $i=2$ the conclusion is trivial. By the unique 
continuation properties of Section \ref{section-basics}, up to choose $I$ and 
$\Omega$ smaller, we may thus assume that $(v_i,\grad v_i)$ are not zero and thus 
that $\alpha(x,t)$ and $\beta(x,t)$ are smooth non-zero functions. Moreover, we 
may fix the normalization $\alpha^2(x,t)+\beta^2(x,t)=1$. Fix 
$(x_0,t_0)$ and set 
$(\tilde\alpha,\tilde\beta)=(\alpha(x_0,t_0),\beta(x_0,t_0))$. 
We notice that the value $(\tilde\alpha,\tilde\beta)$ is taken by 
$(\alpha(x,t),\beta(x,t))$ in a submanifold $\Mc$ of dimension $d'\geq d$ of 
$\Omega\times I$ because the possible values of the function lie in the circle 
$S^1$ which is one-dimensional. The function $w=\tilde\alpha v_1+\tilde\beta 
v_2$ is also a solution of 
\eqref{eq-nodalset2} and by construction $(w,\grad w)$ vanishes in the 
submanifold $\Mc$ of dimension $d'$. We now apply Theorem \ref{th-nodalset} with 
families independent of $\tau\in J=\RR$. The singular nodal set of 
$w(x,t,\tau)$ is $\Mc\times J$ of dimension $d'+1\geq d+1$. Thus $w\equiv 0$ which 
concludes the proof. 
\end{demo}


\end{document}

%% file: transverse.pdf_t
\begin{picture}(0,0)%
\includegraphics{transverse.pdf}%
\end{picture}%
\setlength{\unitlength}{2072sp}%
\begingroup\makeatletter\ifx\SetFigFont\undefined%
\gdef\SetFigFont#1#2#3#4#5{%
  \reset@font\fontsize{#1}{#2pt}%
  \fontfamily{#3}\fontseries{#4}\fontshape{#5}%
  \selectfont}%
\fi\endgroup%
\begin{picture}(10313,3711)(1898,-5483)
\put(3108,-3369){\makebox(0,0)[lb]{\smash{{\SetFigFont{17}{20.4}{\rmdefault}{\mddefault}{\updefault}{\color[rgb]{0,0,0}$\Cc^-$}%
}}}}
\put(9618,-2544){\makebox(0,0)[lb]{\smash{{\SetFigFont{14}{16.8}{\rmdefault}{\mddefault}{\updefault}{\color[rgb]{0,0,0}$u(t)$}%
}}}}
\put(6123,-4884){\makebox(0,0)[lb]{\smash{{\SetFigFont{14}{16.8}{\rmdefault}{\mddefault}{\updefault}{\color[rgb]{.82,0,0}$W^s(\Cc^-)$}%
}}}}
\put(9283,-4681){\makebox(0,0)[lb]{\smash{{\SetFigFont{17}{20.4}{\rmdefault}{\mddefault}{\updefault}{\color[rgb]{0,0,0}$\Cc^+$}%
}}}}
\put(11060,-5338){\makebox(0,0)[lb]{\smash{{\SetFigFont{14}{16.8}{\rmdefault}{\mddefault}{\updefault}{\color[rgb]{.82,0,0}$W^s(\Cc^+)$}%
}}}}
\put(4636,-2266){\makebox(0,0)[lb]{\smash{{\SetFigFont{14}{16.8}{\rmdefault}{\mddefault}{\updefault}{\color[rgb]{0,.69,0}$W^u(\Cc^-)$}%
}}}}
\put(11386,-4156){\makebox(0,0)[lb]{\smash{{\SetFigFont{14}{16.8}{\rmdefault}{\mddefault}{\updefault}{\color[rgb]{0,.69,0}$W^u(\Cc^+)$}%
}}}}
\end{picture}%

%% file: preuve.pdf_t
\begin{picture}(0,0)%
\includegraphics{preuve.pdf}%
\end{picture}%
\setlength{\unitlength}{2072sp}%
\begingroup\makeatletter\ifx\SetFigFont\undefined%
\gdef\SetFigFont#1#2#3#4#5{%
  \reset@font\fontsize{#1}{#2pt}%
  \fontfamily{#3}\fontseries{#4}\fontshape{#5}%
  \selectfont}%
\fi\endgroup%
\begin{picture}(17554,17487)(-1832,-18067)
\put(9181,-17251){\makebox(0,0)[lb]{\smash{{\SetFigFont{17}{20.4}{\rmdefault}{\mddefault}{\updefault}{\color[rgb]{0,0,1} }%
}}}}
\put(-1349,-17161){\makebox(0,0)[lb]{\smash{{\SetFigFont{17}{20.4}{\rmdefault}{\mddefault}{\updefault}{\color[rgb]{0,0,1}any heteroclinic orbits and does}%
}}}}
\put(-1349,-17701){\makebox(0,0)[lb]{\smash{{\SetFigFont{17}{20.4}{\rmdefault}{\mddefault}{\updefault}{\color[rgb]{0,0,1}not meet the projections of $\Nc_\pm$}%
}}}}
\put(7831,-16261){\makebox(0,0)[lb]{\smash{{\SetFigFont{17}{20.4}{\rmdefault}{\mddefault}{\updefault}{\color[rgb]{0,0,1}A place where the projection}%
}}}}
\put(7831,-16801){\makebox(0,0)[lb]{\smash{{\SetFigFont{17}{20.4}{\rmdefault}{\mddefault}{\updefault}{\color[rgb]{0,0,1}of $u_\sigma$ is one-to-one and where}%
}}}}
\put(7831,-17341){\makebox(0,0)[lb]{\smash{{\SetFigFont{17}{20.4}{\rmdefault}{\mddefault}{\updefault}{\color[rgb]{0,0,1}it is easier to construct a suitable}%
}}}}
\put(7831,-17881){\makebox(0,0)[lb]{\smash{{\SetFigFont{17}{20.4}{\rmdefault}{\mddefault}{\updefault}{\color[rgb]{0,0,1}perturbation $h$ to modify $u_\sigma$  }%
}}}}
\put(-1349,-16621){\makebox(0,0)[lb]{\smash{{\SetFigFont{17}{20.4}{\rmdefault}{\mddefault}{\updefault}{\color[rgb]{0,0,1}projection intercepts the one of}%
}}}}
\put(-1439,-16081){\makebox(0,0)[lb]{\smash{{\SetFigFont{17}{20.4}{\rmdefault}{\mddefault}{\updefault}{\color[rgb]{0,0,1}The trace of the set $E$ whose}%
}}}}
\put(4501,-10861){\makebox(0,0)[lb]{\smash{{\SetFigFont{14}{16.8}{\rmdefault}{\mddefault}{\itdefault}{\color[rgb]{0,0,0}the space where the nonlinearities $f$ are defined}%
}}}}
\put(631,-2716){\makebox(0,0)[lb]{\smash{{\SetFigFont{14}{16.8}{\rmdefault}{\mddefault}{\updefault}{\color[rgb]{0,0,0}$\Nc_-$}%
}}}}
\put(3781,-5686){\makebox(0,0)[lb]{\smash{{\SetFigFont{10}{12.0}{\rmdefault}{\mddefault}{\updefault}{\color[rgb]{0,.56,0}$W^u(\Gamma_-,\Nc_-,f_0)$}%
}}}}
\put(4186,-8611){\makebox(0,0)[lb]{\smash{{\SetFigFont{17}{20.4}{\rmdefault}{\mddefault}{\updefault}{\color[rgb]{0,0,.82}$\Gamma_+$}%
}}}}
\put(6076,-9556){\makebox(0,0)[lb]{\smash{{\SetFigFont{14}{16.8}{\rmdefault}{\mddefault}{\updefault}{\color[rgb]{0,0,0}$\Nc_+$}%
}}}}
\put(3151,-4111){\makebox(0,0)[lb]{\smash{{\SetFigFont{17}{20.4}{\rmdefault}{\mddefault}{\updefault}{\color[rgb]{0,0,.82}$\Gamma_-$}%
}}}}
\put(9991,-1051){\makebox(0,0)[lb]{\smash{{\SetFigFont{17}{20.4}{\rmdefault}{\mddefault}{\updefault}{\color[rgb]{0,0,0}$u_\sigma(t)$}%
}}}}
\put(3601,-2941){\makebox(0,0)[lb]{\smash{{\SetFigFont{17}{20.4}{\rmdefault}{\mddefault}{\updefault}{\color[rgb]{0,0,0}$\sigma$}%
}}}}
\put(10351,-7531){\makebox(0,0)[lb]{\smash{{\SetFigFont{17}{20.4}{\rmdefault}{\mddefault}{\updefault}{\color[rgb]{.56,0,0}Projection of $X^ \alpha$ by}%
}}}}
\put(10351,-8206){\makebox(0,0)[lb]{\smash{{\SetFigFont{17}{20.4}{\rmdefault}{\mddefault}{\updefault}{\color[rgb]{.56,0,0}the evaluation ${\rm Ev}(x_\sigma,\cdot)$}%
}}}}
\put(10351,-8881){\makebox(0,0)[lb]{\smash{{\SetFigFont{17}{20.4}{\rmdefault}{\mddefault}{\updefault}{\color[rgb]{.56,0,0}at the point $x_\sigma$}%
}}}}
\end{picture}%

%% file: Sard-Smale.pdf_t
\begin{picture}(0,0)%
\includegraphics{Sard-Smale.pdf}%
\end{picture}%
\setlength{\unitlength}{1243sp}%
\begingroup\makeatletter\ifx\SetFigFont\undefined%
\gdef\SetFigFont#1#2#3#4#5{%
  \reset@font\fontsize{#1}{#2pt}%
  \fontfamily{#3}\fontseries{#4}\fontshape{#5}%
  \selectfont}%
\fi\endgroup%
\begin{picture}(22877,7474)(-1089,-6983)
\put(6351,-1171){\makebox(0,0)[lb]{\smash{{\SetFigFont{17}{20.4}{\rmdefault}{\mddefault}{\updefault}{\color[rgb]{0,0,.69}$\Wc$}%
}}}}
\put(-1074,-3421){\makebox(0,0)[lb]{\smash{{\SetFigFont{17}{20.4}{\rmdefault}{\mddefault}{\updefault}{\color[rgb]{.56,0,0}$\Phi(\Mc,\lambda)$}%
}}}}
\end{picture}%